\documentclass[10pt,twoside]{amsart}

\usepackage{amsmath,amssymb,mathrsfs}
\usepackage{oldgerm}

\input xypic
\xyoption{all}
\CompileMatrices 

\newdir{ >}{{}*!/-10pt/\dir{>}}

\def\smono[#1]{\ar@{->}[#1]|@{|}}
\newcommand{\xymono}{\ar@{ >->}}
\newcommand{\xyepi}{\ar@{->>}}

\newcommand{\colim}{\operatorname{colim}}

\renewcommand{\lim}{\operatorname{lim}}

\newcommand{\Fun}{\operatorname{Fun}}

\newcommand{\Ar}{\operatorname{Ar}}

\newcommand{\Ch}{\operatorname{Ch}}

\newcommand{\Mor}{\operatorname{Mor}}

\newcommand{\Hom}{\operatorname{Hom}}

\newcommand{\can}{\operatorname{can}}
\newcommand{\quis}{\operatorname{quis}}
\newcommand{\Quis}{\operatorname{-\, quis}}

\newcommand{\Qcoh}{\operatorname{Qcoh}}
\newcommand{\str}{\operatorname{str}}

\newcommand{\lax}{\operatorname{lax}}
\newcommand{\Ac}{\operatorname{Ac}}
\newcommand{\const}{\operatorname{\Gamma}}
\newcommand{\Lmod}{\operatorname{-Mod}}
\newcommand{\Rmod}{\operatorname{Mod\phantom{,}-}}
\newcommand{\Bimod}{\operatorname{-Bimod}}
\newcommand{\LRmod}{\operatorname{-Mod\hspace{.2ex}-}\hspace{-.3ex}}

\newcommand{\Vect}{\operatorname{Vect}}
\newcommand{\sPerf}{\operatorname{sPerf}}

\newcommand{\Tot}{\operatorname{Tot}}
\newcommand{\bmu}{\mu}

\newcommand{\sptL}{{\Bbb L}}

\newcommand{\ffi}{\varphi}

\newcommand{\sCx}{\, s{\rm Cx}}

\newcommand{\Z}{{{\Bbb Z}}}

\newcommand{\N}{{{\Bbb N}}}

\renewcommand{\1}{{{1\!\!1}}}

\newcommand{\ie}{{\it i.e.,\,}}

\newcommand{\im}{{\rm im}\, }
\newcommand{\A}{\mathcal{A}}
\newcommand{\B}{\mathcal{B}}
\newcommand{\C}{\mathcal{C}}
\newcommand{\D}{\mathcal{D}}
\newcommand{\E}{\mathcal{E}}
\renewcommand{\H}{\mathcal{H}}

\newcommand{\scK}{\mathcal{K}}
\renewcommand{\L}{{L}}
\newcommand{\M}{\mathcal{M}}
\newcommand{\cN}{\mathcal{N}}
\newcommand{\cP}{\mathcal{P}}

\newcommand{\scR}{\mathcal{R}}
\newcommand{\scS}{\mathcal{S}}
\newcommand{\T}{\mathcal{T}}
\newcommand{\V}{\mathcal{V}}

\newcommand{\coker}{\operatorname{coker}}

\newcommand{\U}{\mathcal{U}}

\newcommand{\eps}{\varepsilon}

\newcommand{\Qed}{\hfil\penalty-5\hspace*{\fill}$\square$}

\newcommand{\GW}{{{\Bbb G}W}}

\newcommand{\meps}{\epsilon}

\title[The Mayer-Vietoris principle for Grothendieck-Witt groups]{The
  Mayer-Vietoris principle for Grothendieck-Witt groups of schemes}
\author{Marco Schlichting}
\address{Marco Schlichting 
\newline Department of Mathematics, Louisiana State University, Baton Rouge, LA 70803-4918
\newline
\phantom{Dep}and
\newline Max-Planck-Institut f\"ur Mathematik, Postfach 7280, D-53072 Bonn, Germany
}

\thanks{The author was partially supported by NSF grant DMS-0604583}

\email{mschlich@math.lsu.edu {\it and} mschlich@mpim-bonn.mpg.de}



\begin{document}

\begin{abstract}
We prove localization and Zariski-Mayer-Vietoris for higher Gro-thendieck-Witt
groups, 
alias hermitian $K$-groups, of schemes admitting an ample family of
line-bundles. 
No assumption on the characteristic is needed, and our schemes can be
singular.
Along the way, we prove additivity, fibration and approximation theorems for
the hermitian $K$-theory of exact categories with weak equivalences and
duality. 
\end{abstract}

\maketitle

\bibliographystyle{alpha}

\vspace{3ex}

\tableofcontents

\section{Introduction}

A classical invariant of a scheme $X$ is its Grothendieck-Witt group $GW_0(X)$
of symmetric bilinear spaces over $X$.
According to Knebusch \cite[\S4]{knebusch:queensLect}, this is
the abelian group generated by isometry classes $[\V,\ffi]$ of
vector bundles $\V$ over $X$ equipped with a non-singular symmetric bilinear
form $\ffi: \V \otimes_{O_X}\V \to O_X$ modulo the relations
$[(\V,\ffi)\perp (\V',\ffi')]= [\V,\ffi]+[\V',\ffi']$ and $[\M,\ffi]=[\H(\cN)]$
for every metabolic space $(\M,\ffi)$ with Lagrangian subbundle $\cN =
\cN^{\perp} \subset \M$ and associated hyperbolic space 
$\H(\cN)$.
Grothendieck-Witt groups naturally occur in ${\Bbb A}^1$-homotopy theory
\cite{Morel:pi0} and are to oriented Chow groups what
algebraic $K$-theory is to ordinary Chow groups, see \cite{bargeMorel},
\cite{faselSrinivas}, \cite{hornbostel:bloch}. 

Using a hermitian version of Quillen's $Q$-construction, 
we have defined in \cite{myGWex} the higher
Grothendieck-Witt groups 
$GW_i(X)$, $i\in \N$, of a scheme $X$, generalizing the
group $GW_0(X)$. 
The purpose of this article is
to prove the following Mayer-Vietoris principle for open covers.

\subsection{Theorem}
\label{thm:MVintro}
{\it
Let $X=U\cup V$ be a scheme with an ample family of line-bundles
(e.g., quasi-projective over an affine scheme, or regular separated noetherian)
which is covered by two open quasi-compact subschemes $U,V \subset X$.
Then there is a long exact sequence, $i \in \Z$,

$$\cdots GW_{i+1}(U\cap V) \to GW_i(X) \to GW_i(U) \oplus GW_i(V) \to
GW_i(U\cap V)  \to GW_{i-1}(X) \cdots$$
}

\noindent
This is a
special case of our theorem \ref{cor:fullMV}
which also includes versions of theorem \ref{thm:MVintro}
for skew-symmetric forms, for forms with
coefficients in line-bundles other than $O_X$ and
for certain
non-commutative schemes.
Note that we don't need the common assumption $\frac{1}{2}\in \Gamma(X,O_X)$,
and $X$ can be singular!

Theorem \ref{thm:MVintro} is a consequence of two theorems, ``Localization'' and ``Zariski-excision''.
To explain the implication, 
let $X$ be a scheme, $L$ a line bundle on $X$, $n \in \Z$ an integer, and $Z
\subset X$ a closed subscheme with open complement $U$. 
With this set of data, we associate in definition \ref{dfn:GWnXonZ} a
topological space 
$GW^n(X\phantom{i}on\phantom{i}Z,\phantom{i}L)$
which,
for $Z=X$, $n=0$ and $L=O_X$, yields the Grothendieck-Witt space $GW(X)$
introduced in \cite{myGWex} whose homotopy groups are the higher
Grothendieck-Witt groups $GW_i(X)$ in theorem \ref{thm:MVintro}, see corollary
\ref{cor:sPerf=exVect}. 
The space $GW^n(X\phantom{i}on\phantom{i}Z,\phantom{i}L)$
is the Grothendieck-Witt space (as defined in
\ref{dfn:GWSpace}) of an exact category with weak equivalences and duality,
namely, 
the exact category of bounded chain complexes of vector bundles on $X$
which are (cohomologically) 
supported in $Z$, equipped with the set of quasi-isomorphisms as weak
equivalences and duality $E\mapsto Hom(E,L[n])$, where 
$L[n]$ denotes the complex which is $L$ in degree $-n$. 
If $Z=X$, we write $GW^n(X,L)$ for
$GW^n(X\phantom{i}on\phantom{i}Z,\phantom{i}L)$.
The non-negative part of theorem \ref{thm:MVintro} is a consequence of the
following two theorems (proved in theorems \ref{thm:Locn1} and
\ref{thm:ZarExc1}).
They are extended to negative Grothendieck-Witt groups in
\S\ref{sec:NegGW} (theorems \ref{thm:LocnSp1} and \ref{thm:ZarExcSp1}).

\subsection{Theorem \rm (Localization)}
\label{thm:LocnIntro}
{\it
Let $X$ be a scheme with an ample family of line-bundles, let $U \subset X$ be
a quasi-compact open subscheme with closed complement $Z=X-U$.
Let $L$ be a line bundle on $X$, and $n\in \Z$ an integer.
Then there is a homotopy fibration 
$$GW^n(X\phantom{i}on\phantom{i}Z,\phantom{i}L) \longrightarrow
GW^n(X,\phantom{i}L) \longrightarrow GW^n(U,\phantom{i}j^*L).$$ 
}

\subsection{Theorem \rm (Zariski excision)}
\label{thm:ZerDescintro}
{\it
Let $j:U\subset X$ be quasi-compact open subscheme of a scheme $X$ which has an
ample family of line-bundles. 
Let $Z \subset X$ be a closed subset such that $Z \subset U$.
Then restriction of vector-bundles induces a homotopy equivalence for all
$n\in\Z$ and all line bundles $L$ on $X$
$$GW^n(X\phantom{i}on\phantom{i}Z,\phantom{i}L)
\stackrel{\sim}{\longrightarrow}
GW^n(U\phantom{i}on\phantom{i}Z,\phantom{i}j^*L).$$
}

Theorems \ref{thm:MVintro} -- \ref{thm:ZerDescintro} have well-known analogs
in algebraic $K$-theory proved by 
Thomason in \cite{TT} based on the work of Waldhausen \cite{wald:spaces} and
Grothendieck et al. \cite{SGA6}.
In fact, our theorems \ref{thm:Locn1}, \ref{thm:LocnSp1}, \ref{thm:ZarExc1},
\ref{thm:ZarExcSp1} and \ref{cor:fullMV} -- special cases of which are
theorems \ref{thm:LocnIntro}, \ref{thm:ZerDescintro} and \ref{thm:MVintro} --
are generalizations of the corresponding theorems in Thomason's work.
More recently, Balmer \cite{balmer:gersten} and Hornbostel
\cite{hornbostel:reps} proved results reminiscent of our theorems
\ref{thm:MVintro} -- \ref{thm:ZerDescintro}.
Both need $X$ to be regular noetherian and separated
and they need $2$ to be a unit in the ring of regular functions on $X$. 
Balmer works with (triangular) Witt-groups instead of Grothendieck-Witt
groups, and 
Hornbostel works with Karoubi's hermitian $K$-groups of rings extended to
regular separated schemes using Jouanolou's device of replacing such a
scheme by an affine vector-bundle torsor.

Neither Balmer's nor Hornbostel's methods can be generalized to cover our
theorems \ref{thm:LocnIntro}, \ref{thm:ZerDescintro} and \ref{thm:MVintro}.
This is because the assumption $\frac{1}{2} \in \Gamma(X,O_X)$ is ubiquitous
in their work, the analog of theorem \ref{thm:MVintro} for
Balmer's triangular Witt groups fails to hold for singular quasi-projective
schemes (see \cite{myWstable} for a counter example even with $\frac{1}{2}\in
\Gamma(X,O_X)$),  
Hornbostel imposes homotopy invariance which doesn't hold for singular
schemes, and his proof uses Karoubi's fundamental theorem
\cite{karoubi:annals} which fails to hold for higher Grothendieck-Witt groups
when $\frac{1}{2} \notin \Gamma(X,O_X)$ (see \cite{myGWder} for a counter
example). 
Instead, we generalize Thomason's work \cite{TT}.
His proofs of the $K$-theory analogs of theorems \ref{thm:LocnIntro},
\ref{thm:ZerDescintro} and \ref{thm:MVintro} are based on a fibration theorem
of Waldhausen \cite[1.6.4]{wald:spaces} and on ``invariance of $K$-theory
under derived equivalences'' \cite[Theorem 1.9.8]{TT} which itself is a
consequence of Waldhausen's approximation theorem \cite[1.6.7]{wald:spaces}.
We prove in theorem \ref{thm:ChgOfWkEq} the analog of Waldhausen's fibration
theorem for higher
Grothendieck-Witt groups.
Its proof, however, is not a formal consequence of
``additivity'' (proved for higher Grothendieck-Witt groups in
\S\ref{sec:Addthms}), contrary to the $K$-theory situation.
Our proof relies on the author's cone construction in \cite{myGWex}.
``Invariance under derived equivalences'' as well as the naive generalization
of Waldhausen's approximation theorem fail to hold
for higher 
Grothendieck-Witt groups when ``2 is not a unit'' (see \cite{myGWder} for a
counter example).
We prove in theorems \ref{thm:ApprFunctorialCalcFrac}
and \ref{thm:ApprCalcOfFrac2}
versions of Waldhausen's approximation theorem for higher
Grothendieck-Witt groups.
Though not as general as one might wish, they are
enough to show theorems \ref{thm:LocnIntro},
\ref{thm:ZerDescintro} and \ref{thm:MVintro} and their generalizations in
\S\ref{sec:LocnZar} and \S\ref{sec:NegGW}.

\subsection*{Prerequisites}
The article can be read independently from \cite{TT} and \cite{wald:spaces},
though, of course, much of our inspiration derives from these two papers.
The reader is advised to have some background in homotopy theory in the form
of \cite[I-IV]{goerssJardine} and in the
theory of triangulated categories in the form of \cite{keller:uses},
\cite{neeman:locn}, \cite[\S1-2]{neeman:grothendieck}.
Also, we will frequently use results from \cite{myGWex}.

\subsection*{Acknowledgments}
Part of the results of this article were obtained and written down while I was
visiting I.H.E.S. in Paris and Max-Planck-Institut in Bonn. 
I would like to express my gratitude for their hospitality.

\section{The Grothendieck-Witt space}

In this section we introduce the Grothendieck-Witt group and the
Grothendieck-Witt space of an exact category with weak equivalences and
duality (Definitions \ref{dfn:GW0} and \ref{dfn:GWSpace}), and we show in
proposition \ref{Qh=Seh} that the Grothendieck-Witt space of an exact category
with duality defined here is equivalent to the one defined in \cite{myGWex}.
We start with recalling definitions from \cite{myGWex}.
Note that our terminology (for ``category with duality'', ``duality preserving
functor'') sometimes differs from standard terminology
as in \cite{scharlau:book}, \cite{knus:book}.

\subsection{Categories with duality, $\C_h$  and  form functors}
\label{subsec:castWduals}
A {\it category
with duality} is a triple $(\C,*,\eta)$ with $\C$ a category, $*: \C^{op} \to
\C$ a functor, $\eta: 1 \to **$ a natural transformation, called {\em double
  dual identification}, such that 
$1_{A^{*}}=\eta_A^{*}\circ \eta_{A^{*}}$ for all objects $A$ in $\C$.
If $\eta$ is a natural isomorphism, we say that the duality is {\it strong}.
In case $\eta$ is the identity (in which case $**=id$), we call the duality
{\it strict}.

A {\it symmetric form}  in a category with duality
$(\C,*,\eta)$ is a pair $(X,\ffi)$ where $\ffi:X \to
X^{*}$ is a morphism in $\C$ satisfying $\ffi^{*}\eta_X=\ffi$.
A map of symmetric forms $(X,\ffi) \to (Y,\psi)$ is a map $f:X \to Y$ in $\C$
such that $\ffi=f^*\circ\psi\circ f$.
Composition of such maps is composition in $\C$.
For a category with duality $(\C,*,\eta)$, we denote by
$\C_h$ the {\em category of symmetric forms in $\C$}.
It has objects the symmetric forms in $\C$ and maps the maps between symmetric
forms. 

A {\it form functor}  from a category with duality
$(\A,*,\alpha)$ to another such category $(\B,*,\beta)$  is a pair $(F,\ffi)$
with 
$F:\A \to \B$ a functor and $\ffi: F* \to *F$ a natural transformation, called
{\em duality compatibility morphism},
such that $\ffi_A^*\beta_{FA}=\ffi_{A^*}F(\alpha_A)$ for every object $A$ of
$\A$.
There is an evident definition of composition of form functors, see
\cite[3.2]{myGWex}. 
The category $\Fun(\A,\B)$ of functors $\A \to \B$ is a category with duality,
where the dual
$F^{\sharp}$ of a functor $F$ is $*F*$, and double dual identification 
$\eta_F: F \to F^{\sharp\sharp}$ at an object $A$ of $\A$ is the map 
$\beta_{F(A^{**})}\circ F(\alpha_A) = F(\alpha_A)^{**}\circ \beta_{FA}$.
To give a form functor $(F,\ffi)$ is the same as to give a symmetric form
$(F,\hat{\ffi})$ in the category with duality $\Fun(\A,\B)$
in view of the formulas ${\ffi}_A=F(\alpha_A)^*\circ \hat{\ffi}_{A^*}$ and
$\hat{\ffi}_A={\ffi}_{A^*}\circ F(\alpha_A)$.
A natural transformation $(F,\ffi) \to (G,\psi)$ of form functors is a map
$(F,\hat{\ffi}) \to (G,\hat{\psi})$ of symmetric forms in $\Fun(\A,\B)$.

A {\em duality preserving functor} between categories with duality
$(\A,*,\alpha)$ and $(\B,*,\beta)$ is a functor $F: \A \to \B$ which commutes
with dualities and double dual identifications, that is, we have
$F*=*F$ and $F(\alpha)=\beta_{F}$.
In this case, $(F,id)$ is a form functor.
We will consider duality preserving
functors $F$ as form functors $(F,id)$.
Note that our use of the phrase ``duality preserving functor'' may differ from
its use by other authors! 
\vspace{2ex}

Recall that an {\em exact category} is 
 an additive category $\E$ equipped with a family of 
sequences of maps in $\E$, called {\it conflations} (or {\em admissible short
  exact sequences}, or simply {\em exact sequences}), 
$$X \stackrel{i}{\to} Y \stackrel{p}{\to} Z$$
satisfying a list of axioms, see 
\cite{quillen:higherI}, \cite[\S 4]{keller:uses}, \cite[2.1]{myGWex}.
The map $i$ in an exact sequence is called inflation (or admissible
monomorphism) and may be depicted as $\rightarrowtail$, and the 
map $p$ is called {\it deflation} (or admissible 
epimorphism) and may be depicted as $\twoheadrightarrow$ in diagrams.
Unless otherwise stated, all exact categories in this article will be
(essentially) small.

\subsection{Exact categories with weak equivalences}
\label{subsec:ExCatWeEq}
An {\it exact category with weak equivalences} is a pair $(\E,w)$ 
with $\E$ an exact category and $w \subset \Mor\E$ a set of
morphisms, called weak equivalences, which contains all identity morphisms, is
closed under isomorphisms, retracts, push-outs along inflations, pull-backs
along deflations, composition and the $2$ out of three property for
composition (if $2$ of the $3$ maps among $a$, $b$, $ab$ are in $w$ then so is
the third).
A weak equivalence is usually depicted as $\stackrel{_{\sim}}{\to}$ in diagrams.
A functor $F:\A \to \B$ between exact categories with weak equivalences 
$(\A,w)$ and $(\B,w)$ is called {\em exact} if it sends conflations to
conflations and weak equivalences to weak equivalences.

For an exact category with weak equivalences $(\E,w)$, we will write $w\E$ for
the subcategory of weak equivalences in $\E$. 
Its objects are the objects
of $\E$ and its maps the maps in $w$.
Also, we will regard an exact category $\E$ (without specifying weak
equivalences) as an exact category with weak
equivalences $(\E,i)$ where $i$ is the set of isomorphisms in $\E$.

\subsection{Exact categories with weak equivalences and duality}
An {\it exact category with weak equivalences and duality} is a quadruple
$(\E,w,*,\eta)$
with $(\E,w)$ an exact category with weak equivalences and
$(\E,*,\eta)$ a category with duality such that 
$*:(\E^{op},w) \to (\E,w)$ is an exact functor (in particular, $*(w) \subset
w$) and 
$\eta: id \to **$ is a natural weak equivalence, that is, $\eta_X \in w$ for
all objects $X$ in $\E$. 
We may simply say $\E$ or $(\E,w)$ is an exact category with weak equivalences
and duality if the remaining data are understood.
Note that if $\E$ is an exact category with weak equivalences and duality, the
category $w\E$ of weak equivalences in $\E$ is a category with duality.

A symmetric form $(X,\ffi)$ in $(\E,w,*,\eta)$ is called {\em non-singular} if
$\ffi$ is a weak equivalence.
In this case, we call the pair $(X,\ffi)$ a {\em symmetric space} in
$(\E,w,*,\eta)$.
A form functor $(F,\ffi): (\A,w,*,\eta)\to (\B,w,*,\eta)$ is called {\em
  exact} if $F: (\A,w) \to (\B,w)$ is exact.
It is called {\em non-singular}, if the duality compatibility morphism $\ffi:
F* \to *F$ is a natural weak equivalence. 

An {\it exact category with duality} is an exact category with weak
equivalences and duality where the set of weak equivalences is the set of
isomorphisms. 
In particular, the double dual identification has to be a natural isomorphism.

\subsection{Definition}
\label{dfn:GW0}
The {\em Grothendieck-Witt group} 
$$GW_0(\E,w,*,\eta)$$ 
of an exact category with weak equivalences and duality $(\E,w,*,\eta)$
is the free abelian group generated by isomorphism classes $[X,\ffi]$ of
symmetric spaces $(X,\ffi)$ in $(\E,w,*,\eta)$, subject to the following
relations 
\begin{enumerate}
\item
\label{cor:itm0:Kh_0}
$[X,\ffi]+[Y,\psi] = [X \oplus Y,\ffi \oplus \psi]$
\item
\label{cor:itm1:Kh_0}
if $g:X \to Y$ is a weak equivalence, then $[Y,\psi]=[X,g^*\psi g]$,
and
\item
\label{cor:itm2:Kh_0}
if $(E_{\bullet},\ffi_{\bullet})$ is a symmetric space in the category
of exact sequences in $\E$, that is, a map 
$$\xymatrix{
\hspace{2ex}E_{\bullet} \ar[d]^{\ffi_{\bullet}}_{\wr}:  & E_{-1}
\xymono[r]^i\ar[d]^{\ffi_{-1}}_{\wr} & E_0 \xyepi[r]^p\ar[d]^{\ffi_{0}}_{\wr}
&  E_1 \ar[d]^{\ffi_{1}}_{\wr} \\
\hspace{2ex}E_{\bullet}^*: & E_1^* \xymono[r]_{p^*} & E_0^* \xyepi[r]_{i^*} &
E_{-1}^*}$$ 
of exact sequences with
$(\ffi_{-1},\ffi_0,\ffi_1)=(\ffi_1^*\eta,\ffi_0^*\eta,\ffi_{-1}^*\eta) $ a 
weak equivalence, then $$[E_0,\ffi_0] = \left[ E_{-1}\oplus E_1,
\left(\begin{smallmatrix} 0 & \ffi_1\\ \ffi_{-1} &
    0\end{smallmatrix}\right)\right].
$$
\end{enumerate}

\subsection{Remark}
If in definition \ref{dfn:GW0}, the set of weak equivalences is the set of
isomorphisms, then we recover the classical Grothendieck-Witt group of an
exact category with duality, see for instance \cite[2.9]{myGWex}.
In this case, relation \ref{cor:itm2:Kh_0} says that the class $[E_0,\ffi_0]$
of a metabolic space
$(E_0,\ffi_0)$ with Lagrangian $i:E_{-1}\rightarrowtail E_0$ is equivalent in
the Grothendieck-Witt group to the class of the hyperbolic space $\H(E_{-1})$
of the Lagrangian $E_{-1}$.
In particular, if $\E$ is the category $\Vect(X)$ of vector bundles on $X$,
$*$ is the duality functor $E \mapsto Hom(E,O_X)$ and $\eta$ is the usual
canonical double dual identification, the group $GW_0(\E,i,*,\eta)$ is
Knebusch's Grothendieck-Witt group $GW_0(X)$ of a scheme $X$, denoted $L(X)$
in \cite{knebusch:queensLect}.

\subsection{Definition}
The {\em Witt} group 
$$W_0(\E,w,*,\eta)$$ 
of an exact category with weak equivalences and duality $(\E,w,*,\eta)$ 
is the free abelian group generated by isomorphism classes $[X,\ffi]$ of
symmetric spaces $(X,\ffi)$ in $(\E,w,*,\eta)$, subject to the
relations \ref{dfn:GW0} \ref{cor:itm0:Kh_0}, \ref{cor:itm1:Kh_0} and
\begin{itemize}
\item[(c')]
if $(E_{\bullet},\ffi_{\bullet})$ is a symmetric space in the category
of exact sequences in $(\E,w,*,\eta)$, then $[E_0,\ffi_0]=0$.
\end{itemize}
\vspace{2ex}

The hermitian $S_{\bullet}$-construction of \cite{yaoShapiro},
\cite[1.5]{jensme}, which gives rise to the Grothendieck-Witt space to be
defined in \ref{dfn:GWSpace}, is the edgewise subdivision of
Waldhausen's $S_{\bullet}$-construction \cite{wald:spaces}.
We review the relevant definitions and start with the edgewise subdivision of a
simplicial object, see
\cite[1.9 Appendix]{wald:spaces}, \cite[Appendix
1]{segal:configuration}.

\subsection{Edgewise subdivision}
\label{subsec:edgewise}
Let $\Delta$ be the category with objects $[n]$ 
the totally ordered sets $[n]=\{0<1<...<n\}$, $n\in \N$,
and morphisms the order preserving maps.
Let $\underline{n}$ be the totally ordered set 
$$\underline{n} = \{n' < (n-1)' < ... <0'<0<...<n\}.$$
It is (uniquely) isomorphic to $[2n+1]$.
The assignment $T: [n]\mapsto \underline{n}$ defines a functor $\Delta \to
\Delta$
where a map $\theta:[n] \to [m]$ goes to the map
$T(\theta):\underline{n} \to \underline{m}: p \mapsto \theta(p), p'\mapsto
\theta(p)'$. 
For a simplicial object $X_{\bullet}$, 
the {\em edge-wise subdivision} $X_{\bullet}^e$  of $X_{\bullet}$ is 
the simplicial object $X_{\bullet}\circ T$.
The inclusion $[n] \hookrightarrow \underline{n}:i \mapsto i$ defines a map
$X^e \to X$ of simplicial objects.
It is known \cite[Appendix
1]{segal:configuration} that for a simplicial set $X_{\bullet}$, the topological
realization of $X_{\bullet}$ and of its edge-wise subdivision $X^e_{\bullet}$
are homeomorphic.
We need the following (well-known) variant.

\subsection{Lemma}
\label{lem:Xe->X}
{\it
For any simplicial set $X_{\bullet}$, the map $X^e_{\bullet} \to X_{\bullet}$
is a homotopy equivalence.}

\subsection*{\it Proof}
Let $X_{\bullet}$ be a simplicial set.
For a small category $\C$, write $X_{\C}$ for the set 
$\Hom(N_*\C,X_{\bullet})$ of simplicial maps from the nerve
$N_*\C$ of $\C$ to $X_{\bullet}$.
Note that $X^e$ is the simplicial set $[n]\mapsto X_{\underline{n}}$.
We define bisimplicial sets $X^e_{\bullet\bullet}$ and $X_{\bullet\bullet}$ by 
the formulas $X^e_{m,n} = X_{\underline{m}\times [n]}$ and
$X_{m,n} = X_{[m]\times [n]}$.
Consider the following diagram of bisimplicial sets
$$\xymatrix{X^e_{0\bullet} \ar[r]^{\sim} \ar[d]^{\wr} & X^e_{\bullet\bullet}
  \ar[d] & X^e_{\bullet 0} \ar[l]_{\sim} \ar[d]\\
X_{0\bullet} \ar[r]^{\sim} & X_{\bullet\bullet}
   & X_{\bullet 0} \ar[l]_{\sim}
}$$
in which the horizontal maps are the canonical inclusions of horizontally
respectively vertically constant bisimplicial sets, and the vertical maps are
induced by the inclusions $[m]\subset \underline{m}$.
Once we show that all arrows labeled $\stackrel{\sim}{\to}$ are homotopy
equivalences, we are done, because the right vertical map can be identified
with $X_{\bullet}^e \to X_{\bullet}$.

The left vertical map $X^e_{0\bullet} \to X_{0\bullet}$ is a homotopy
equivalence since it
can be identified with the map $X_{\bullet}^I \to X_{\bullet}$ which is
evaluation at $0$, where $I=N_*\underline{0}\cong N_*[1]$ is the standard
simplicial interval. 
In order to see that the upper right horizontal map $X^e_{\bullet 0} \to
X^e_{\bullet\bullet}$ is a 
homotopy equivalence, it suffices to prove that for every $n$, the map 
$X^e_{\bullet 0} \to X^e_{\bullet n}$ is a homotopy equivalence of simplicial
sets.
Since the map $[0] \to [n]:0 \mapsto 0$ induces a retraction $X^e_{\bullet n} \to
X^e_{\bullet 0}$,
we have to show that the composition $X^e_{\bullet n} \to
X^e_{\bullet 0} \to X^e_{\bullet n}$ is homotopic to the identity.
The (unique) natural transformation from the constant functor $[n] \to [n]:i
\mapsto 0$ to the identity functor $[n] \to [n]$ defines a functor 
$h:[1]\times [n] \to [n]$ such that the restrictions to $\{i\} \times [n] \to
[n]$, $i=0,1$ are the constant respectively the identity functor.
We have a map of simplicial sets
\begin{equation}
\label{eqn:Xe}
I^e\times X^e_{\bullet n} \to X^e_{\bullet n}
\end{equation}
which in degree $m$ sends the pair $(\xi,f) \in I_{\underline{m}}\times
X_{\underline{m}\times [n]}$ to the composition
$$N_*(\underline{m} \times [n]) 
\stackrel{(1,\xi)\times 1}{\longrightarrow} 
N_*(\underline{m}\times [1] \times [n])
\stackrel{(1\times h)}{\longrightarrow} 
N_*(\underline{m} \times [n])
\stackrel{f}{\longrightarrow} X_{\bullet}
$$
Since the two points $\{0,1\}=\{0,1\}^e\subset I^e$ are path connected in
$I^e$, the map (\ref{eqn:Xe}) defines the desired homotopy.
The other horizontal homotopy equivalences in the diagram are similar, and we
omit the details.
\Qed 
\vspace{2ex}

Now we recall Waldhausen's $S_{\bullet}$-construction \cite[\S1.3]{wald:spaces}.

\subsection{Waldhausen's $S_{\bullet}$-construction}
Let $\Ar[n]$ denote the category whose objects are the arrows of
$[n]=\{0<1<...<n\}$ and 
whose morphisms are the commutative squares in $[n]$.
For an exact category with weak equivalences $(\E,w)$, 
Waldhausen defines $S_n\E \subset \Fun(\Ar[n],\E)$ as the full subcategory of
the category $\Fun(\Ar[n],\E)$ of
functors $$A: \Ar[n] \to \E:(p\leq q) \mapsto A_{p,q}$$
for which $A_{p,p}=0$ and $A_{p,q} \rightarrowtail A_{p,r} \twoheadrightarrow
A_{q,r}$ is a conflation 
whenever $p\leq q \leq r$, $p,q,r \in [n]$.
The category $S_n\E$  is an exact category with weak equivalences 
where a sequence $A \to B \to C$ of functors $\Ar[n]\to \E$ in $S_n\E$
is exact if 
$A_{p,q} \rightarrowtail B_{p,q} \twoheadrightarrow C_{p,q}$ is exact in $\E$,
and a map $A \to B$ of functors in $S_n\E$ is a weak equivalence if 
$A_{p,q} \to B_{p,q}$ is a weak equivalence in $\E$ for all $p\leq q \in [n]$.

The cosimplicial category $n \mapsto \Ar[n]$ makes the assignment 
$n \mapsto S_n\E$ into a simplicial exact category with weak equivalences.
According to \cite{wald:spaces}, \cite{TT}, the $K$-theory space $K(\E,w)$ of
an exact category with weak equivalences $(\E,w)$ is the space
$$K(\E,w)=\Omega |wS_{\bullet}\E|.$$

\subsection{The hermitian $S_{\bullet}$-construction}
\label{dfn:Se_dot}
The category $[n]$ has a unique structure of a category with strict duality
$[n]^{op} \to [n]: i \mapsto n-i$.
This induces a strict duality on the category $\Ar[n]$ of arrows in $[n]$.
For an exact category with weak equivalences and duality $(\E,w,*,\eta)$, the
category $\Fun(\Ar[n],\E)$ is therefore a category with duality (see
\ref{subsec:castWduals}).
This duality preserves the subcategory $S_n\E \subset \Fun(\Ar[n],\E)$, and 
makes $S_n\E $ into an exact category with weak equivalences and duality.
It turns out that the simplicial structure maps of $n \mapsto S_n\E$ are not
compatible with dualities.
However, its edgewise subdivision 
$$S^e_{\bullet}\E: n \mapsto S_n^e\E = S_{2n+1}\E$$
is a simplicial exact
category with weak equivalences and duality; the simplicial structure maps
being duality preserving.
Considering $S_n^e\E$ as a full subcategory of $\Fun(\Ar(\underline{n}),\E)$,
the dual $A^*$ of an object $A: \Ar(\underline{n}) \to \E$ satisfies
$(A^*)_{p,q}=A_{q',p'}^*$ for $p\leq q \in \underline{n}$ where $p''$ trickily
denotes $p$.
The double dual identification $A \to A^{**}$ at $(p\leq q)$ is
$\eta_{A_{p,q}}$.

\subsection{Definition}
\label{dfn:GWSpace}
Let $(\E,w,*,\eta)$ be an exact category with weak equivalences and duality.
By \ref{dfn:Se_dot}, the assignment $n \mapsto S^e_{n}\E$ defines a
simplicial exact category $S^e_{\bullet}\E$ with weak equivalences and duality.
The subcategories of weak equivalences define a simplicial category with
duality $n \mapsto wS^e_{n}\E$.
Taking associated categories of symmetric forms (see \ref{subsec:castWduals}),
we obtain a simplicial category $(wS^e_{\bullet}\E)_h$.

The composition $(wS^e_{\bullet}\E)_h \to wS^e_{\bullet}\E \to wS_{\bullet}\E$
of simplicial categories, in which 
the first arrow is the forgetful functor $(X,\ffi) \mapsto X$, 
and the second is the canonical map $X^e_{\bullet} \to X_{\bullet}$ of
simplicial objects (see \ref{subsec:edgewise}),
yields a map of classifying spaces
\begin{equation}
\label{eqn:GWdfn}
|(wS^e_{\bullet}\E)_h| \to |wS_{\bullet}\E|
\end{equation}
whose homotopy fibre (with respect to
a zero object  of $\E$ as base point of $wS_{\bullet}\E$) is defined to be the
{\it 
  Grothendieck-Witt space} 
$$GW(\E,w,*,\eta)$$ of $(\E,w,*,\eta)$.
If $(*,\eta)$ are understood, we may simply write $GW(\E,w)$ instead of
$GW(\E,w,*,\eta)$. 
We define the
{\em higher Grothendieck-Witt groups} of $(\E,w,*,\eta)$ as the homotopy
groups $$GW_i(\E,w,*,\eta) = \pi_i GW(\E,w,*,\eta), \phantom{1234}
i\geq 1,$$ 
and show in proposition \ref{prop:pi0GW=GW0} below that $\pi_0GW(\E,w,*,\eta)
\cong GW_0(\E,w,*,\eta)$, so that our definition here extends that in
\ref{dfn:GW0}. 

\subsection{Remark}
Orthogonal sum makes the spaces $(wS_{\bullet}^e\E)_h$ and $GW(\E,w,*,\eta)$
into commutative $H$-spaces.
Since the commutative monoid of connected components of these spaces are
groups (see proposition \ref{prop:pi0GW=GW0} and remark \ref{rem:pi0Sdot}
below), both spaces are actually commutative $H$-groups.

\subsection{Functoriality}
A non-singular exact form functor $(F,\ffi):(\A,w,*,\eta) \to (\B,w,*,\eta)$
between exact categories with weak equivalences and duality induces maps 
$$
\renewcommand\arraystretch{1.5} 
\begin{array}{rcl}
(F,\ffi): & (wS^e_{\bullet}\A)_h \to (wS^e_{\bullet}\B)_h:& (A,\alpha) \mapsto
(FA,\ffi_A F(\alpha)),\phantom{a}{\rm and}\\
F:& wS_{\bullet}\A \to wS_{\bullet}\B:& A \mapsto FA
\end{array}$$
of simplicial categories compatible with composition of form functors.
Taking homotopy fibres of $(wS^e_{\bullet})_h \to wS_{\bullet}$, we obtain an
induced map
$$GW(F,\ffi): GW(\A,w,*,\eta) \to GW(\B,w,*,\eta)$$
of associated Grothendieck-Witt spaces.
For the next lemma, recall that a natural transformation of form functors
$(F,\ffi) \to (G,\psi)$ is a map of associated symmetric forms in
$\Fun(\A,\B)$.
It is a natural weak equivalence if $FA \to GA$ is a weak equivalence for all
objects $A$ of $\A$.

\subsection{Lemma}
\label{lem:htpyOfFunctors}
{\it Let $(F,\ffi) \stackrel{\sim}{\to} (G,\psi)$ be a natural weak
  equivalence of non-singular exact form functors 
$(\A,w,*,\eta) \to (\B,w,*,\eta)$ between exact categories
  with weak equivalences and duality.
Then,  on associated Grothendieck-Witt spaces, $(F,\ffi)$ and
$(G,\psi)$ induce homotopic maps  
$GW(\A,w,*,\eta) \to GW(\B,w,*,\eta).$
}

\subsection*{\it Proof}
The natural weak equivalence $(F,\ffi) \stackrel{\sim}{\to} (G,\psi)$ induces
natural transformations of functors $(wS^e_{n}\A)_h \to (wS^e_{n}\B)_h$ and
$wS_{n}\A \to wS_{n}\B$.
These natural transformations define functors 
$[1] \times (wS^e_{n}\A)_h \to (wS^e_{n}\B)_h$ and $[1]\times wS_{n}\A \to
wS_{n}\B$ whose restrictions to $0, 1 \in [1]$ are the two given functors.
They are compatible with the simplicial structure and induce, after
topological realization, the homotopy between $GW(F,\ffi)$ and $GW(G,\psi)$. 
\vspace{-3ex}

\Qed
\vspace{2ex}

Next, we will associate to every exact category with weak equivalences
$(\E,w)$ a 
category with weak equivalences and duality $(\H\E,w)$ such that the
Grothendieck-Witt space of $(\H\E,w)$ is equivalent to the $K$-theory space of
$(\E,w)$. 
In this sense, Grothendieck-Witt theory is a generalization of algebraic
$K$-theory. 

\subsection{Hyperbolic categories}
Let $\C$ be a category.
Its hyperbolic category is the category with strict duality 
$\H\C=(\C\times \C^{op},*)$  
where $(X,Y)^* = (Y,X)$.
For any category with duality $\A$ there is a functor $\A_h \to \A:
(X,\ffi) \mapsto X$ that ``forgets the forms''.
We define the functor $(\H\C)_h \to \C$ as the composition 
of the functor
$(\H\C)_h \to \H\C$ and the projection $\H\C=\C \times \C^{op} \to \C$ onto the
first factor.

\subsection{Lemma}
\label{lem:HCh=C}
{\it
For any small category $\C$, the functor 
$(\H\C)_h \to \C$
is a homotopy equivalence.
}

\subsection*{\it Proof}
The category $(\H\C)_h$ of symmetric forms in $\H\C$ is isomorphic to the
category whose objects 
are maps $f:X \to Y$ in $\C$ and where a map from $f$ to $f':X' \to Y'$ is a
pair of maps $a:X \to X'$, $b:Y' \to Y$ in $\C$ such that $f=bf'a$.
Composition is composition in $\C$ of the $a$'s and $b$'s.
The functor $(\H\C)_h \to \C$ in the lemma sends the object $f:X \to Y$
to $X$ and the map $(a,b)$ to $a$.
Write $F$ for this functor, and let $A$ be an object of the target category
$\C$. 
We will show that the comma categories $(A\downarrow F)$ are contractible.
By Quillen's theorem A \cite[\S1]{quillen:higherI}, this implies the lemma.

The category $(A\downarrow F)$ has objects sequences $A \stackrel{x}{\to} X
\stackrel{f}{\to} Y$ of maps in $\C$.
A morphism from $(x,f)$ to $A \stackrel{x'}{\to} X'
\stackrel{f'}{\to} Y'$ is a pair $a:X \to X'$, $b:Y' \to Y$ of maps in $\C$
such that $x'=ax$ and $f=bf'a$.
In particular, the category $(A \downarrow F)$ is non-empty as $(1_A,1_A)$ is
one of its objects.
Let $\C_0 \subset (A \downarrow F)$ be the full subcategory of objects $(x,f)$
with $x=1_A$.
The inclusion has a right adjoint $(A \downarrow F) \to \C_0: (x,f) \mapsto
(1_A,fx)$ with counit of adjunction $(1_A,fx) \to (x,f)$ given by the pair of
maps $x:A \to X$ and $1:Y \to Y$.
It follows that the inclusion $\C_0 \subset (A \downarrow F)$ is a homotopy
equivalence. 
Since the category $\C_0$ has a terminal object, namely $(1_A,1_A)$, the
categories $\C_0$ and $(A \downarrow F)$ are contractible.
\Qed
\vspace{2ex}

If $(\E,w)$ is an exact category with weak equivalences, we make $\H\E$ into an
exact category with weak equivalences and (strict) duality by declaring a map
$(a,b):(X,Y) \to (X',Y')$ in $\H\E$ to be a weak equivalence if $a:X \to X'$
and $b: Y' \to Y$ are weak equivalences in $\E$.
Note that $w\H\E=\H w\E$ as categories with strict duality.

\subsection{Proposition}
\label{prop:GWH=K}
{\it Let $(\E,w)$ be an exact category with weak equivalences, then there is a
  natural homotopy equivalence 
$$GW(\H\E,w) \simeq K(\E,w).$$
}

\subsection*{\it Proof}
Consider the commutative diagram of simplicial categories
$$\xymatrix{
(wS^e_{\bullet}\H\E)_h \ar[r]^{=} \ar[d] & (\H wS^e_{\bullet}\E)_h \ar[d]
\ar[dr]^{\sim} &\\ 
wS^e_{\bullet}\H\E \ar[r]_{\hspace{-6ex} =} \ar[d]_{\wr} &
wS^e_{\bullet}\E\times 
(wS^e_{\bullet}\E)^{op} \ar[d]^{\wr} \ar[r]_{\hspace{6ex} p_1} &
wS^e_{\bullet}\E \ar[d]^{\wr}\\  
wS_{\bullet}\H\E \ar[r]^{\hspace{-6ex} =}  & wS_{\bullet}\E\times
(wS_{\bullet}\E)^{op}  \ar[r]^{\hspace{6ex} p_1} & wS_{\bullet}\E
}
$$
where the upper vertical maps are the functors that ``forget the forms''.
The 
lower vertical maps are induced by the inclusion $[n] \hookrightarrow
\underline{n}$ and are thus homotopy equivalences, by lemma \ref{lem:Xe->X}.
The diagonal map is a homotopy equivalence by lemma \ref{lem:HCh=C}.
By the ``octahedron axiom'' for homotopy fibres applied to the upper right
triangle, it follows that the homotopy fibre of the composition of the left
vertical 
maps is equivalent to the loop space of the fibre of 
$p_1:wS_{\bullet}\E\times (wS_{\bullet}\E)^{op}  \to wS_{\bullet}\E$
which is the loop of the fibre of 
$(wS_{\bullet}\E)^{op}  \to ({\rm point})$ which is $K(\E,w)$.
\Qed
\vspace{2ex}

\subsection{Remark}
\label{rmk:dfnGW2}
Let $(\E,w,*,\eta)$ be an exact category with weak equivalences and duality.
The functor
$$F: (\E,w,*,\eta) \to (\H\E,w): X \mapsto (X,X^*)$$
together with the duality compatibility morphism $(1,\eta_X): (X^*,X^{**}) \to
(X^*,X)$ is called {\em forgetful form functor}.
It is a non-singular exact
form functor between exact categories 
with weak equivalences and duality.
The map $(wS^e_{\bullet}\E)_h \to wS_{\bullet}\E$ defining the
Grothendieck-Witt space factors as
$$(wS^e_{\bullet}\E)_h \stackrel{F}{\to} (wS^e_{\bullet}\H\E)_h
\stackrel{\sim}{\to} wS_{\bullet}\E,$$ 
where 
the second map is the homotopy equivalence in the diagram of the proof of
proposition \ref{prop:GWH=K} (going right, diagonally and down).
It follows that the Grothendieck-Witt space
$GW(\E,w,*,\eta)$ is naturally homotopy equivalent to the homotopy fibre of
$$(wS^e_{\bullet}\E)_h \stackrel{F}{\to}(wS^e_{\bullet}\H\E)_h.$$ 
\vspace{2ex}

We finish the section with a comparison result between the definition of the
Grothendieck-Witt space of an exact category with duality $(\E,i,*,\eta)$ in
terms of the hermitian $S_{\bullet}$-construction and the definition given in
\cite[Definition 4.6]{myGWex} in terms of the hermitian $Q$-construction.
We recall the relevant definitions.

\subsection{The hermitian $Q$-construction}
Recall from \cite{quillen:higherI} that for an exact category $\E$ its
$Q$-construction is the category with objects the objects of $\E$ and maps $X
\to Y$ equivalence classes of diagrams
\begin{equation}
\label{eqn:quillenQdatum}
X \stackrel{p}{\twoheadleftarrow} U \stackrel{i}{\rightarrowtail} Y
\end{equation}
with $p$ a deflation and $i$ an inflation.
The datum $(U,i,p)$ is equivalent to $(U',i',p')$ if there is an
isomorphism $g:U \to U'$ in $\E$ such that $p=p'\circ g$ and $i=i'\circ g$.
The composition in $Q\E$ of maps $X \to Y$ and $Y \to Z$ represented by the
data $(U,i,p)$ and $(V,j,q)$ is given by the datum
$(U\times_YV,p\bar{q},j\bar{i})$ where $\bar{q}:U\times_YV \to U$ and
$\bar{i}:U\times_YV \to V$ are the canonical projections to $U$ and $V$,
respectively.

For an exact category with duality $(\E,*,\eta)$, the hermitian
$Q$-construction is the category $Q^h(\E,*,\eta)$
with objects the symmetric spaces $(X,\ffi)$ in $\E$.
 A map $(X,\ffi) \to (Y,\psi)$ is a map $X \to Y$ in Quillen's
 $Q$-construction, that is, an equivalence class of diagrams 
(\ref{eqn:quillenQdatum}), such that the square of maps
$p$, $i$, $p^*\ffi$ and $i^*\psi$ is commutative and bicartesian.
Composition of maps is as in Quillen's $Q$-construction.
For more details, we refer the reader to \cite[4.2]{myGWex} and the references
in \cite[Remark 4.3]{myGWex}.
\vspace{2ex}

In \cite[Definition 4.6]{myGWex}, we defined the Grothendieck-Witt space of 
$(\E,*,\eta)$ as the homotopy fibre of the forgetful functor 
$Q^h\E \to Q\E: (X,\ffi) \mapsto X$.
The following proposition reconciles this definition with the one given in 
\ref{dfn:GWSpace}.
This allows us to freely use the results proved in \cite{myGWex}.

\subsection{Proposition}
\label{Qh=Seh}
{\it
For an exact category with duality $(\E,*,\eta)$, there are natural homotopy
equivalence between $|Q^h\E|$ and $|(iS_{\bullet}^e\E)_h|$ and between the
homotopy fibre of the forgetful functor $|Q^h\E| \to |Q\E|$ and the
Grothendieck-Witt space as defined in \ref{dfn:GWSpace}.
}

\subsection*{\it Proof}
For the first homotopy equivalence, the proof is the same as in \cite[1.9
Appendix]{wald:spaces}. 
In detail, let $iQ^h_{\bullet}\E$ be the simplicial category which in degree
$n$ is the category $iQ^h_n\E$ whose objects are sequences $X_0 \to X_1 \to
... \to X_n$ of maps in $Q^h\E$ and a map in $iQ^h_n\E$ is an isomorphism of
such sequences. 
As $n$ varies, $iQ^h_n\E$ defines a simplicial category where face and
degeneracy maps are defined as in the usual nerve construction.
The nerve of $iQ^h_{\bullet}\E$ as a bisimplicial set is isomorphic to the
nerve of the simplicial category which in degree $m$ are the sequences $X_0
\to X_1 \to ... \to X_m$ of isomorphisms in $Q^h\E$ and where maps are maps of
sequences in $Q^h\E$ (which are not necessarily isomorphisms).
The latter simplicial category is degree-wise equivalent to $Q^h\E$ (via the
embedding of $Q^h\E$ as the constant sequences).
Thus the latter simplicial category (and therfore also $iQ^h_{\bullet}\E$) is
homotopy equivalent to $Q^h\E$.

Every object $(A_{p,q})_{p\leq q \ \in \ \underline{n}}$ in $(S^e_n\E)_h$ defines a
  string of maps 
$$A_{0',0} \to A_{1',1} \to ... \to A_{n',n}$$ in $Q^h\E$.
This defines a map 
$(iS^e_{\bullet}\E)_h \to iQ_{\bullet}^h\E$
of simplicial categories which is degree-wise an equivalence.
Therefore, this map defines a homotopy equivalence on topological
realizations.

For the second homotopy equivalence, consider the commutative diagram of
topological spaces
$$\xymatrix{
|(iS^e_{\bullet}\E)_h| \ar[d] & |(iS^e_{\bullet}\E)_h| \ar[l]_1 \ar[r]^{\sim}
\ar[d] & |iQ^h_{\bullet}\E| \ar[d] & |Q^h\E|\ar[l]_{\sim} \ar[d]\\
|iS_{\bullet}\E| & |iS^e_{\bullet}\E| \ar[l]_{\sim} \ar[r]^{\sim}
 & |iQ_{\bullet}\E|  & |Q\E|\ar[l]_{\sim} 
}$$
in which the lower right two horizontal maps are defined in a similar way as
their hermitian analogs above them (see \cite[1.9 Appendix]{wald:spaces}), the
three right vertical maps are ``forgetful'' functors and 
all maps labeled $\stackrel{\sim}{\to}$ are homotopy equivalences.
It follows that the homotopy fibre of the first vertical map is equivalent to
the homotopy fibre of the last vertical map.
\Qed

\section{Additivity theorems}
\label{sec:Addthms}

Additivity theorems are fundamental in algebraic $K$-theory.
They imply, for instance,  Waldhausen's fibration theorem
\cite[1.6.4]{wald:spaces} which is the basis for the $K$-theory version of
theorem \ref{thm:LocnIntro}.
In this section, we prove the analogs of additivity for higher
Grothendieck-Witt theory. 
In order to formulate them, we recall the
concept of ``admissible short complexes'' from \cite[\S7]{myGWex}.

\subsection{Admissible short complexes and their homology}
\label{sec:shrtCx}
Let $(\E,w,*,\eta)$ be an exact category with weak equivalences and duality.
A {\em short complex} in $\E$ is a complex 
\begin{equation}
\label{eqn:sCx}
A_{\bullet}:\phantom{123456} 0 \to A_{1} \stackrel{d_1}{\to} A_0
\stackrel{d_0}{\to} A_{-1} \to 0, \phantom{1234}(d_0\circ d_1=0)
\end{equation}
in $\E$ concentrated in degrees $-1,0,1$.
It is {\em admissible} if 
$d_1$ and $d_{0}$ are inflation and
deflation, respectively, and the map 
$A_1 \to \ker(d_0)$ (or equivalently $\coker(d_1) \to A_{-1}$) is an inflation
(deflation).
A sequence $A_{\bullet} \to B_{\bullet} \to C_{\bullet}$ of admissible
short complexes  is {\em exact} if $A_i \to B_i \to C_i$ is exact
in $\E$;
a map $A_{\bullet} \to B_{\bullet}$ is a weak equivalence if $A_i \to B_i$ is
a weak equivalence in $\E$, $i=-1,0,1$.
The dual of the complex (\ref{eqn:sCx}) is the (admissible short)
complex
$$A^*_{-1} \stackrel{d_0^*}{\to} A_0^* \stackrel{d_1^*}{\to} A_1^*,$$
and the double dual identification
$\eta_{\A_{\bullet}}:A_{\bullet} \to 
(A_{\bullet})^{**}$ is $\eta_{A_i}$ in degree $i=-1,0,1$.
We denote by $(\sCx(\E),w,*,\eta)$ the exact category with weak equivalences
and duality of admissible short complexes in $\E$.

If $(A_{\bullet},\alpha_{\bullet})$ is a symmetric form in $\sCx(\E)$, we write
$H_0(A_{\bullet},\alpha_{\bullet})$ for its zeroth homology symmetric form. 
Its underlying object is $H_0(A_{\bullet})=\ker(d_0)/\im(d_1)$, and it is
equipped with the form
$\bar{\alpha}$ which is the unique symmetric form such that
$\bar{\alpha}_{|\ker(d_0)} = \alpha_{|\ker(d_0)}$.
This makes $H_0:\sCx(\E) \to \E$ into a non-singular exact form functor for
every exact category with weak equivalences and duality $\E$.
For more details, we refer the reader to \cite[\S7]{myGWex}.
\vspace{2ex}

In the special case of an exact category with duality (where all weak
equivalences are isomorphisms), the following two
theorems were proved in \cite[Theorems 7.1, 7.4]{myGWex}.
The $K$-theory version is due to Waldhausen in \cite[Theorem 1.4.2 and
proposition 1.3.2]{wald:spaces},
in view of the equivalence $S^e_1\E \to \sCx(\E):A\mapsto (A_{1'0'}\to A_{1'1}
\to A_{01})$ of exact categories with weak equivalences and duality.

\subsection{Theorem \rm (Additivity for short complexes)}
\label{thm:addtysCx}
{\it Let $(\E,w,*,\eta)$ be an exact category with weak equivalences and
  duality. 
Then the non-singular exact form functor $H_0:\sCx(\E) \to \E$ together with
the exact functor $ev_1:\sCx(\E) \to \E:A_{\bullet} \mapsto A_1$ induce
homotopy equivalences $(H_0,ev_1):$
$$
\renewcommand\arraystretch{1.5} 
\begin{array}{rcl}
 (wS^e_{\bullet}\sCx \E)_h & \stackrel{\sim}{\longrightarrow} &
(wS^e_{\bullet} \E)_h \times wS_{\bullet}\E\\
 GW(\sCx \E, w,*,\eta) & \stackrel{\sim}{\longrightarrow} &
GW(\E,w,*,\eta) \times K(\E,w).
\end{array}$$  
}
\vspace{2ex}

For the next theorem, recall that a form functor $(\A,*,\eta) \to (\B*,\eta)$
between categories with duality is nothing else than a symmetric form 
$(F,\ffi)$ in the category with duality of functors
$(\Fun(\A,\B),\sharp,\eta)$, see \ref{subsec:castWduals}.

\subsection{Theorem \rm(Additivity for functors)}
\label{thm:addtyFun}
{\it
Let $(\A,w,*,\eta)$ and $(\B,w,*,\eta)$ be exact categories with weak
equivalences and duality.
Given a non-singular exact form functor $(F_{\bullet},\ffi_{\bullet}):
(\A,w,*,\eta) \to \sCx (\B,w,*,\eta)$, that is,
a commutative diagram of exact functors $F_i: (\A,w) \to (\B,w)$
$$\xymatrix{
\hspace{2ex}F_{\bullet} \ar[d]^{\ffi_{\bullet}}_{\wr}:  & F_{1}
\xymono[r]^{d_1}\ar[d]^{\ffi_{1}}_{\wr} & F_0 \xyepi[r]^{d_0}\ar[d]^{\ffi_{0}}_{\wr}
&  F_{-1} \ar[d]^{\ffi_{-1}}_{\wr} \\
\hspace{2ex}F_{\bullet}^{\sharp}: & F_{-1}^{\sharp} \xymono[r]_{d_0^{\sharp}} & F_0^{\sharp}
\xyepi[r]_{d_1^{\sharp}} & F_{1}^{\sharp}}$$ 
with $d_0d_1=0$, $F_1 \rightarrowtail \ker{d_0}$ an inflation, and 
$(\ffi_{1},\ffi_0,\ffi_{-1})=(\ffi_{-1}^{\sharp}\eta,\ffi_0^{\sharp}\eta,\ffi_{1}^{\sharp}\eta)$. 
Then the two non-singular exact form functors 
\begin{equation}
\label{eqn:AddtyFun}
(F_0,\ffi_0) \phantom{as}{\rm and}\phantom{wer} 
H_0(F_{\bullet},\ffi_{\bullet}) \perp \left(F_1\oplus F_{-1}, \left(\begin{smallmatrix} 0 & \ffi_{-1}\\ \ffi_{1} &
    0\end{smallmatrix}\right)\right)
\end{equation}
induce homotopic maps on Grothendieck-Witt spaces
$$GW(\A,w,*,\eta) \to GW(\B,w,*,\eta).$$
}
\vspace{2ex}

We will reduce the proofs of the additivity theorems \ref{thm:addtysCx} and
\ref{thm:addtyFun} to the case of exact
categories with dualities dealt with in \cite[\S7]{myGWex}.
This will be done with the help of the simplicial resolution lemma
\ref{lem:simplResLem} below.
For that, we need to replace an exact category with weak equivalences and
duality (where the double dual identification is a natural weak
equivalence) 
by one with a strong duality (where the double dual identification is a
natural isomorphism) without changing its hermitian $K$-theory.
This will be done in the strictification lemma \ref{lem:strictification}
below. 
\vspace{2ex}

We introduce notation for lemma \ref{lem:strictification}.
Let $ExWeDu$ be the category of small exact categories with weak equivalences
and duality; and non-singular exact form functors as maps.
Recall that a category with duality $(\A,*,\eta)$ has a strict duality if
$\eta=id$ (and in particular, $**=id$). 
Let $ExWeDu^{\str}$ be the category of small exact categories with weak
equivalences and strict duality; and duality preserving functors as maps.
Write $\lax:ExWeDu^{\str} \subset ExWeDu$ for the natural inclusion.

\subsection{Lemma \rm (Strictification lemma)}
\label{lem:strictification}
{\it There is a (strictification) functor 
 $$\str: ExWeDu \to ExWeDu^{str}: (\A,w,*,\eta) \mapsto (\A^{\str}_w,w,\sharp,id)$$
and natural transformations
$(\Sigma,\sigma):id \to \lax\circ\str$ and $(\Lambda,\lambda):\lax\circ\str
\to id$ such that the 
compositions $(\Sigma,\sigma)\circ (\Lambda,\lambda)$ and
$(\Lambda,\lambda)\circ (\Sigma,\sigma)$ are weakly equivalent to the
identity form functor.  
In particular, for any exact category with weak equivalences and duality
$(\A,w,*,\alpha)$, we have a homotopy equivalence of Grothendieck-Witt spaces
$$GW(\Sigma,\sigma):GW(\A,w,*,\alpha) \stackrel{\sim}{\longrightarrow}
GW(\A^{\str}_w,w,\sharp,id).$$ 
}

\subsection*{\it Proof}
The construction of the strictification functor is as follows.
Let $(\A,w,*,\alpha)$ be an exact category with weak equivalences and duality.
The objects of $\A^{\str}_w$ are triples $(A,B,f)$ with $A$, $B$ objects of
$\A$ and  $f:A \stackrel{\sim}{\to} B^*$ a weak equivalence in $\A$.   
A morphism from $(A_0, B_0,f_0)$ to $(A_1, B_1,f_1)$ 
  is a pair $(a,b)$ of morphisms $a:A_0 \to A_1$ and $b:B_1 \to B_0$ in $\A$
  such that $f_1 a = b^* f_0$.
Composition is composition of the $a$'s and $b$'s in $\A$.
A map $(a,b)$ is a weak equivalence if $a$ and $b$ are weak equivalences in
$\A$. 
A sequence $(A_0,B_0,f_0) \to
 (A_1,B_1, f_1) \to (A_2,B_2,f_2)$ is exact if $A_0 \to
 A_1\to A_2$ and $B_2 \to   B_1 \to B_0$ are exact in $\A$.
The duality $\sharp: (\A^{str}_w)^{op}\to \A^{str}_w$ is defined by
  $(A,B,f: A\to B^*)^{\sharp} = (B,A,f^* \alpha_{B})$ on objects, and 
  by $(a,b)^{\sharp}= (b,a)$ on morphisms.
The double dual identification is the identity natural transformation.
The category thus constructed $\A^{str}_w= (\A^{str}_w,w ,\sharp , id)$ is an
exact category with weak equivalences and strict duality.
We may write $\A^{\str}$, or $\A^{str}_w$ for  $(\A^{str}_w,w ,\sharp , id)$
if the remaining data are understood.
If $(F,\ffi):(A,w,*,\alpha) \to (\B,v,*,\beta)$ is a non-singular exact form
functor,
its image  under the strictification functor 
 is the functor $F^{str}: \A^{str}_w \to \B^{str}_v$ given by
 $F^{str}:(A,B,f) \mapsto 
(FA,FB,{\ffi}_B F(f))$ on objects, and by $(a,b)\mapsto (Fa,Fb)$ on morphisms. 
One checks that $F^{str}$ preserves composition and commutes with dualities.

The natural transformations $(\Sigma,\sigma):id \to \lax \circ
\str$ and $(\Lambda,\lambda):  \lax \circ \str \to id$ are defined as
follows.
The functor 
$\Sigma:\A \to \A^{str}_w$
sends an object $A$ to $(A,A^*,\alpha_A)$
and a morphism $a$ to $(a,a^*)$. 
The duality compatibility morphism $\sigma: \Sigma(A^*) \to \Sigma(A)^{\sharp}
$ is the map $(1,\alpha_A): (A^*,A^{**},\alpha_{A^*}) \to (A^{*},A,1)$.
The functor $\Lambda: \A^{str} \to \A$ sends the object $(A,B,f)$ to $A$ and a
morphism $(a,b)$ to $a$, and is equipped with the duality compatibility
morphism  $\lambda:\Lambda[(A,B,f)^{\sharp}] \to [\Lambda(A,B,f)]^*$ the map
$f^*\alpha_B:B \to A^*$.

The composition 
$(\Sigma,\sigma) \circ (\Lambda,\lambda)$ is the functor sending 
$(A,B,f)$ to $(A,A^*,\alpha_A)$, and the morphism $(a,b)$ to $(a,a^*)$.
It is equipped with the
duality compatibility morphism
$(f^*\alpha_B,f):(B,B^*,\alpha_B) \to (A^*,A,1)$.
There is a natural weak equivalence of form functors $(\Sigma,\sigma)\circ
(\Lambda,\lambda) \stackrel{\sim}{\to} id$ given by the map $(1,f^*\alpha_B): 
(A,A^*,\alpha_A) \to (A,B,f)$.
Regarding the other composition, we have
$(\Lambda,\lambda)\circ (\Sigma,\sigma)=id$.

By lemma \ref{lem:htpyOfFunctors}, the form functor $(\Sigma,\sigma): \A \to
\A^{\str}_w$ induces a homotopy equivalence of Grothendieck-Witt spaces
with homotopy inverse $(\Lambda,\lambda)$.
\vspace{-2.5ex}

\Qed

\subsection{Remark}
Here is a slight generalization of lemma \ref{lem:strictification}.
If $v$ is another set of weak equivalences for $(\A,*,\alpha)$
such that $w\subset v$, then
$(\A_w^{\str},v,\sharp,id)$ is also an exact category with weak equivalences
and strict duality.
The form functors $(\Lambda,\lambda): (\A_w^{\str},v) \to (\A,v)$ and
$(\Sigma,\sigma): (\A,v) \to (\A_w^{\str},v)$ are still exact and non-singular
with compositions that are naturally
weakly equivalent to the identity functors.
In particular, we have a homotopy equivalence
$$
GW(\Sigma,\sigma):GW(\A,v,*,\alpha) \stackrel{\sim}{\longrightarrow}
GW(\A^{\str}_w,v,\sharp,id).
$$
\vspace{2ex}

The next lemma
will sometimes allow us to replace an exact category with weak equivalences
and 
duality by a simplicial exact category with duality (where weak equivalences
and double dual identification are isomorphisms).

\subsection{Notation for lemma \ref{lem:simplResLem}}
Let $(\E,w,*,\eta)$ be an exact category with weak equivalences and duality,
and let $\D$ be an arbitrary (small) category.
Recall from \ref{subsec:castWduals} that the category $\Fun(\C,\E)$ of
functors $\D \to \E$ is a category with duality.
It is an exact category with weak equivalences and duality
if we declare maps $F \to G$ (sequences $F_{-1} \to F_0\to F_1$) of functors
$\D \to \E$ to be a weak equivalence (conflation) if $F(A) \to G(A)$ is a weak
equivalence ($F_{-1}(A) \to F_0(A) \to F_1(A)$ is a
conflation) in $\E$ for all objects $A$ of $\D$.
We write 
$\Fun_w(\D,\E)\subset \Fun(\D,\E)$ for the full subcategory of those
functors $F:\D \to \C$ for which the image $F(d)$ of all maps $d$ of $\D$ are
weak equivalences in $\E$: $F(d)\in w\E$.
The category $\Fun_w(\D,\E)$ inherits from $\Fun(\D,\E)$ the structure of an
exact category with weak equivalences and duality. 
In particular, for $n\in \N$ and $(\E,w,*,\eta)$ an exact category with weak
equivalences and strong duality,  
the category $\Fun_w(\underline{n},\E)$ is an exact category with weak
equivalences and strong duality.
It has objects strings of weak
equivalences and maps commutative diagrams in $\E$.
Varying $n$, the categories $\Fun_w(\underline{n},\E)$ define a simplicial
exact category with weak equivalences and strong duality.
Recall that the symbol $i$ stands for the set of isomorphisms in a category.

\subsection{Lemma \rm (Simplicial resolution lemma)}
\label{lem:simplResLem}
{\it
Let $(\E,w,*,\eta)$ be an exact category with weak equivalences and strong
duality. 
Then there are homotopy equivalences
$$
\renewcommand\arraystretch{1.5} 
\begin{array}{rcl}
(wS_{\bullet}^e\E)_h & \simeq & |n \mapsto
(iS_{\bullet}^e\Fun_w(\underline{n},\E))_h|,\\
GW(\E,w) & \simeq & | n \mapsto GW(\Fun_w(\underline{n},\E),i)|
\end{array}$$
which are functorial for exact form functors $(F,\ffi)$ 
for which $\ffi$ is an isomorphism.
}

\subsection*{\it Proof}
We start with some general remarks.
Let $\C$ be a category, and recall that $i$ stands for the set of isomorphisms
in $\C$.
Since $\C=\Fun_i([0],\C)$, inclusion of degree-zero simplices yields a map of
simplicial categories 
$\C \to (n \mapsto \Fun_i([n],\C))$ which is degree-wise an equivalence of
categories, and thus induces a homotopy equivalence after topological
realization. 
Using the equality of bisimplicial sets
$$p,q \phantom{qw} \mapsto \phantom{qwret} N_p\Fun_i([q],\C) = N_qi\Fun([p],\C),$$
where $N_{\bullet}$ stands for the nerve of a category, 
we obtain a homotopy equivalence 
$|\C| \stackrel{\sim}{\to} | n \mapsto i\Fun([n],\C)|$.
Since the topological realizations of $p\mapsto i\Fun([p],\C)$ and of
$p\mapsto i\Fun([p]^{op},\C)$ are isomorphic,
we have a homotopy equivalence 
\begin{equation}
\label{eqn:simplReal}
\C \stackrel{\sim}{\to} | n \mapsto i\Fun([n]^{op},\C)|.
\end{equation}
The homotopy equivalence is natural in the category $\C$.

Let $(\C,*,\eta)$ be a category with strong duality.
There is an equivalence of categories
$i\Fun([n]^{op},\C_h) \to (i\Fun(\underline{n},\C))_h$
which sends an object 
$$(X_n,\ffi_n) \stackrel{f_n}{\to} (X_{n-1},\ffi_{n-1})
\stackrel{f_{n-1}}{\to} \cdots \stackrel{f_{1}}{\to}(X_0,\ffi_0)$$ of the left
hand category to the object 
$$X_n \stackrel{f_n}{\to} X_{n-1}
\stackrel{f_{n-1}}{\to} \cdots \stackrel{f_{1}}{\to} X_0 \stackrel{\ffi_0}{\to}
X_0^* \stackrel{f_{1}^*}{\to} X_{1}^* \stackrel{f_{2}^*}{\to} \cdots
\stackrel{f_{n}^*}{\to} X_n^*$$
equipped with the form $(\eta_{X_n},...,\eta_{X_0},1,...,1)$.
A map $(g_n,...,g_0)$ (which is an isomorphism compatible with forms) is sent
to the map $(g_n,...,g_0,(g_0^*)^{-1},...,(g_n^*)^{-1}$).
The equivalence is functorial in $[n]\in \Delta$ and thus induces a
homotopy equivalence after topological realization.
Together with (\ref{eqn:simplReal}), we obtain a 
homotopy equivalence of topological spaces
\begin{equation}
\label{eqn:simplRealII}
|\C_h| \stackrel{\sim}{\longrightarrow} | n \mapsto
(i\Fun(\underline{n},\C))_h|
\end{equation}
which is natural for categories with strong duality $(\C,*,\eta)$ and
form functors $(F,\ffi)$ between them for which $\ffi$ is an isomorphism.

For an exact category with weak equivalences and strong duality 
$(\E,w,*,\eta)$, we apply the general homotopy equivalence
(\ref{eqn:simplRealII}) to the form functor  
$wS_p^e\E \to wS_p^e\H\E$ 
induced by the forgetful form functor $\E \to \H\E$ (see remark \ref{rmk:dfnGW2}).
Varying $p$, we obtain a map of homotopy equivalences after topological realization
$$\xymatrix{
(wS_{\bullet}^e\E)_h \ar[r]^{\hspace{-7ex} \sim} \ar[d] & |n \mapsto
(iS_{\bullet}^e\Fun_w(\underline{n},\E))_h| \ar[d]\\
(wS_{\bullet}^e\H\E)_h \ar[r]^{\hspace{-8ex} \sim} & |n \mapsto
(iS_{\bullet}^e\H\Fun_w(\underline{n},\E))_h|.
}$$
The top row gives the first homotopy equivalence of the lemma.
By remark \ref{rmk:dfnGW2}, the left vertical homotopy fibre of the diagram
is $GW(\E,w,*,\eta)$.
In view of Bousfield-Friedlander's theorem \cite[B4]{BF}, \cite[Theorem IV
4.9]{goerssJardine}, the homotopy fibre 
of the right vertical map is the simplicial realization of the degree-wise
homotopy fibres.
By remark \ref{rmk:dfnGW2}, this is
$$|n \mapsto GW( \Fun_w(\underline{n},\E),i,*,\eta)|.$$
\Qed
\vspace{2ex}

Before proving the additivity theorems, we give a first application of the
simplicial resolution lemma and show that 
$\pi_0GW(\E,w,*,\eta)$ is the Grothendieck-Witt group as defined in
\ref{dfn:GW0}.

\subsection{Proposition \rm (Presentation of $GW_0$)}
\label{prop:pi0GW=GW0}
{\it
Let $(\E,w,*,\eta)$ be an exact category with weak equivalences and duality.
There is a natural isomorphism
$$GW_0(\E,w,*,\eta) \cong \pi_0GW(\E,w,*,\eta).$$
}

\subsection*{\it Proof}
In view of relation \ref{dfn:GW0} \ref{cor:itm1:Kh_0} and lemma
\ref{lem:htpyOfFunctors}, weakly equivalent non-singular 
exact form functors $(F,\ffi) \stackrel{\sim}{\longrightarrow} (G,\psi)$
 induce the same map on $GW_0$ and on $\pi_0GW$.
Therefore, we can replace $(\E,w)$ by its strictification $(\E^{\str}_w,w)$
from lemma \ref{lem:simplResLem} which has a strong (in fact strict) duality.
So we can assume the duality on $(\E,w)$ to be strong which will allow us to
use the simplicial resolution lemma \ref{lem:simplResLem}.
For a bisimplicial set $X_{\bullet\bullet}$, there is a co-equalizer diagram 
\begin{equation}
\label{eqn:pf:GW0}
\xymatrix{
\pi_0\, |X_{1\bullet}| \ar@<.7ex>[r]^{d_1}  \ar@<-.7ex>[r]_{d_0}
& \pi_0\, |X_{0\bullet}| \ar[r] & \pi_0\, |X_{\bullet\bullet}|.
}
\end{equation}
This is well known and follows from an examination of the usual skeletal
filtration of $|X_{\bullet\bullet}|= |n\mapsto X_{n,n}| = |n \mapsto
X_{n\bullet}|$ -- which the reader can find in  
\cite[Diagram IV.1 (1.6)]{goerssJardine}, for instance --
using the fact that the functor
$\pi_0:\Delta^{op}\operatorname{Sets} \to \operatorname{Sets}$ preserves
push-out diagrams as it is left adjoint to the inclusion functor 
$\operatorname{Sets} \to \Delta^{op}\operatorname{Sets}$.
By the simplicial resolution lemma \ref{lem:simplResLem} and the fact (proven
in \cite[Proposition 4.14]{myGWex} together with proposition \ref{Qh=Seh})
that the proposition 
holds when the set of weak equivalences is the set of isomorphisms, we deduce
that $\pi_0GW(\E,w)$ is the co-equalizer of the diagram
$$\xymatrix{
GW_0\left\lgroup\Fun_w(\underline{1},\E)\right\rgroup \ar@<.7ex>[r]^{d_1}  \ar@<-.7ex>[r]_{d_0}
&  GW_0\left\lgroup\Fun_w(\underline{0},\E)\right\rgroup
}$$ 
of Grothendieck-Witt groups of exact categories with duality.
It suffices therefore to display $GW_0(\E,w,*,\eta)$ as the co-equalizer of the
same diagram.
Evaluation at the object $0'$ of $\underline{0}$ defines a non-singular exact
form functor 
$F:(\Fun_w(\underline{0},\E),i) \to (\E,w)$ 
which sends the object $f: X_{0'} \to X_0$ to  $F(f)=X_{0'}$ and has duality
compatibility morphism $F(f^*) \to F(f)^*$  
the map $f^*:X_0^* \to X_{0'}^*$.
The form functor induces a map 
$GW_0(\Fun_w(\underline{0},\E),i) \to GW_0(\E,w)$ which equalizes $d_0$ and
$d_1$ in view of the relation \ref{dfn:GW0} \ref{cor:itm1:Kh_0}.
We therefore obtain a map from the co-equalizer of the diagram to $GW_0(\E,w)$.
To construct its inverse, 
consider the map from the free abelian group generated by isomorphism classes 
$[X,\ffi]$ of symmetric spaces in $(\E,w)$ to
$GW_0(\Fun_w(\underline{0},\E),i)$ sending $(X,\ffi)$ to 
$\ffi:X \to X^*$ equipped with the non-singular form $(\eta_X,1)$.
This map is surjective and factors through relations \ref{dfn:GW0}
\ref{cor:itm0:Kh_0} and \ref{cor:itm2:Kh_0}. 
The map induces a surjective map to the co-equalizer which factors through
relation \ref{dfn:GW0} \ref{cor:itm1:Kh_0} and thus induces a surjective map
from $GW_0(\E,w,*,\eta)$ to the co-equalizer.
Since composition with the
map from the co-equalizer to $GW_0(\E,w,*,\eta)$ is the identity, the claim
follows. 
\Qed

\subsection{Remark}
\label{rem:pi0Sdot}
Using the co-equalizer diagram (\ref{eqn:pf:GW0}), one can show directly the
isomorphism 
$\pi_0|(wS^e_{\bullet}\E)_h| \cong W_0(\E,w,*,\eta)$ without the need of the
simplicial resolution lemma. 

\subsection*{\it Proof  of theorem \ref{thm:addtysCx}}
In view of the strictification lemma \ref{lem:strictification} we can assume
the duality on $\E$ to be strong.
By the simplicial resolution lemma \ref{lem:simplResLem} the proof reduces
further to the case of an exact category with duality (in which all weak
equivalences are isomorphisms).
This case was proved in \cite[theorems 7.4, 7.10]{myGWex} in view of
proposition \ref{Qh=Seh}.
\Qed

\subsection{\it Proof of theorem \ref{thm:addtyFun}}
\label{pfofAddtyFun}
Theorem \ref{thm:addtyFun} is a formal consequence of theorem
\ref{thm:addtysCx}.  
Let $(\E,w,*,\eta)$ be an exact category with weak equivalences and duality.
Consider the form functors $\sCx(\E) \to \H\E:E_{\bullet} \mapsto
(E_1,E_{-1}^*)$ and  $\H\E \to \sCx(\E): (E_1,E_{-1}) \mapsto (E_1 \to E_1 \oplus
E_{-1}^* \to E_{-1}^*)$ 
with duality compatibility morphisms
$(1,\eta):(E_{-1}^*,E_1^{**}) \to (E_{-1}^*,E_1)$ and 
$(\eta,\left(\begin{smallmatrix}0&1\\ \eta&0\end{smallmatrix}\right),1): 
(E_{-1} \to E_{-1}\oplus E_1^* \to E_1^*) \to (E_{-1}^{**} \to  E_1^* \oplus
E_{-1}^{**} \to E_1^*)$.
Together with the form functor $H_0:\sCx(\E) \to \E$ 
and the duality preserving functor $\E \to \sCx(\E): E \mapsto (0 \to E
\to 0)$ they define non-singular exact form functors
$$\E \times \H\E \to \sCx(\E) \to \E \times \H\E$$
whose composition is weakly equivalent to the identity functor.
By the additivity theorem for short complexes (theorem \ref{thm:addtysCx}),
the second form functor induces a homotopy equivalence in hermitian
$K$-theory.
It follows that the two form functors induce inverse homotopy equivalences on
Grothendieck-Witt spaces and on hermitian $S_{\bullet}$ constructions.
Therefore, the compositions
$$\A \stackrel{(F_{\bullet},\ffi_{\bullet})}{\longrightarrow} \sCx(\B)
\stackrel{ev_{0}}{\longrightarrow}\B \phantom{123}{\rm and}\phantom{123}  
\A \stackrel{(F_{\bullet},\ffi_{\bullet})}{\longrightarrow} \sCx(\B) \to \B
\times H\B \to \sCx(\B)  
\stackrel{ev_{0}}{\longrightarrow}\B$$ 
 induce homotopic maps in hermitian $K$-theory.
These compositions are (weakly equivalent to) the maps in
(\ref{eqn:AddtyFun}).
\Qed

\subsection{Remark}
Iterated application of theorem \ref{thm:addtyFun} implies homotopy
equivalences 
$$
(wS^e_{\bullet}S^e_n\E)_h\hspace{2ex} \stackrel{\sim}{\longrightarrow}
\hspace{2ex} 
\left\{
\renewcommand\arraystretch{1.8} 
\begin{array}{cl}
(wS^e_{\bullet}\E)_h  \times {\displaystyle \prod_{p=1}^k wS_{\bullet}\E,} &
n=2k+1\\ 
 {\displaystyle \prod_{p=1}^k wS_{\bullet}\E,} & n=2k.
\end{array}
\right.$$
This allows us to identify $(wS^e_{\bullet}S^e_{\bullet}\E)_h$ with the Bar
construction of the $H$-group $wS_{\bullet}\E$ acting on
$(wS^e_{\bullet}\E)_h$ and leads to a homotopy fibration
$$(wS^e_{\bullet}\E)_h \to (wS^e_{\bullet}S^e_{\bullet}\E)_h \to
wS_{\bullet}S_{\bullet}\E$$
in which the first map is 
``inclusion of degree zero simplices'' and the second map is the ``forgetful
map'' $(E,\ffi) \mapsto E$ followed by the canonical homotopy equivalence
$X_{\bullet}^e \to X_{\bullet}$. 
In particular, the iterated hermitian $S_{\bullet}$-construction
$(wS^e_{\bullet}S^e_{\bullet}\E)_h$ is not a delooping of 
$(wS^e_{\bullet}\E)_h$ contrary to the $K$-theory situation, compare
\cite[proposition 1.5.3 and remark thereafter]{wald:spaces}.

\subsection{Remark}
Define the {\em Witt-theory space} $W(\E,w,*,\eta)$ as the colimit of the top
row in the sequence of homotopy fibrations
$$\xymatrix{
GW(\E,w) \ar[r] & (wS^e_{\bullet}\E)_h  \ar[r] \ar[d] &
(wS^e_{\bullet}S^e_{\bullet}\E)_h \ar[r]\ar[d] &
(wS^e_{\bullet}S^e_{\bullet}S^e_{\bullet}\E)_h \ar[r]\ar[d] & \cdots\\
& wS_{\bullet}\E & wS_{\bullet}S_{\bullet}\E &
wS_{\bullet}S_{\bullet}S_{\bullet}\E 
}$$
Since the spaces in the second row get higher and higher
connected, we see that $\pi_0W(\E,w,*,\eta)=W_0(\E,w,*,\eta)$ and
$\pi_1W(\E,i,*,\eta)=W_{form}(\E,*,\eta)$, where the group
$W_{form}(\E,*,\eta)$ is the Witt-group of formations in $(\E,*,\eta)$,
that is, the cokernel of the hyperbolic map 
$K_0(\E) \to GW_{form}(\E):[X] \mapsto [\H X,X,X^*]$ to the Grothendieck-Witt
group of formations $GW_{form}(\E)$ defined in \cite[4.11]{myGWex}.

If $\E$ is a $\Z[\frac{1}{2}]$-linear category and ``complicial'', we show in
\cite{myGWder} that the Grothendieck-Witt space
$W(\E,w,*,\eta)$ is the infinite loop space associated with
(the $(-1)$-connected cover of) Ranicki's $\sptL$-theory
spectrum, and its homotopy groups 
$$\pi_iW(\E,w,*,\eta) = W^{-i}(w^{-1}\E,*,\eta),\phantom{q3}i\geq 0,$$
are Balmer's Witt groups $W^{-i}(w^{-1}\E,*,\eta)$  of the triangulated
category with duality $(w^{-1}\E,*,\eta)$.
At this point, I don't know how to calculate 
$\pi_nW(\E,i,*,\eta)$, $n \geq 2$,
for (complicial) $(\E,w,*,\eta)$ when $\E$ is not
$\Z[\frac{1}{2}]$-linear.

\section{Change of weak equivalences and Cofinality}

In this section we prove in theorems \ref{thm:ChgOfWkEq} and \ref{thm:cofinal}
the higher Grothendieck-Witt theory analogs of Waldhausen's fibration theorem
\cite[1.6.4]{wald:spaces} and of 
Thomason's cofinality theorem \cite[1.10.1]{TT}.

Waldhausen's $K$-theory version of theorem \ref{thm:ChgOfWkEq} below needs a
``cylinder functor''.
The purpose of the next definition is to define the higher Grothendieck-Witt
theory analog.
We first fix some notation.
For an exact category with weak equivalences $(\E,w)$, write $\E^w \subset \E$
for the full subcategory of $w$-acyclic objects, that is, those objects $E$ of
$\E$ 
for which the unique map $0 \to E$ is a weak equivalence.
The category $\E^w$ is closed under extensions in $\E$, and thus inherits an
exact structure from $\E$ such that the inclusion $\E^w \subset \E$ is fully
exact.

\subsection{Definition}
\label{dfn:symCone}
Let $(\E,w,*,\eta)$ be an exact category with weak equivalences and duality.
A {\em symmetric cone} on $(\E,w,*,\eta)$ is given by the following data:
\begin{enumerate}
\item
exact functors $P:\E \to \E^w$, and $C:\E \to \E^w$,
\item
a natural deflation $p_E: PE \twoheadrightarrow E$ and a natural inflation
$i_E:E \rightarrowtail CE$,
\item
a natural map $\gamma_E: P(E^*) \to (CE)^*$ such that
$i^*_E\gamma_E=p_{E^*}$ for all objects $E$ of $\E$.
\end{enumerate}
It is convenient to define $\bar{\gamma}_E:C(E^*) \to (PE)^*$ by 
$\bar{\gamma}_E = P(\eta_E)^*\circ \gamma_{E^*}^*\circ \eta_{C(E^*)}$.
Then $p_E^*=\bar{\gamma}_E\circ i_{E^*}$, and $\gamma_E=C(\eta_E)^*\circ
\bar{\gamma}_{E^*}^* \circ \eta_{P(E^*)}$.
In other words, the sequence $P \to id \to C$ defines an exact form functor
from $\E$ to the category of sequences $E_{-1} \twoheadrightarrow E_0
\rightarrowtail E_1$ in $\E$ with duality compatibility map
$(\gamma,1,\bar{\gamma})$.

For examples of symmetric cones, see \ref{subsec:canSymCone}.
\vspace{2ex}

The proof of the next theorem will occupy most of this section.

\subsection{Theorem \rm (Change of weak equivalences)}
\label{thm:ChgOfWkEq}
{\it 
Let $(\E,w,*,\eta)$ be an exact category with weak equivalences and
duality which has a symmetric cone.
Let $v$ be another set of weak equivalences in $\E$ containing $w$ and which
is closed under the duality.
Then the duality $(*,\eta)$ on $\E$ makes 
$(\E^v,w)$, $(\E^v,v)$, $(\E,v)$ into exact categories with weak equivalences
and duality such that the commutative square of duality preserving inclusions
$$\xymatrix{
(\E^v,w) \ar[r] \ar[d] & (\E^v,v) \ar[d]\\
(\E,w)  \ar[r] & (\E,v)}
$$
induces a 
homotopy cartesian square of associated
Grothendieck-Witt spaces.
Moreover, the upper right
corner has contractible Grothendieck-Witt space. 
}

\subsection{Remark}
\label{rem:HgpsHcart}
A square of homotopy commutative $H$-groups (such as Grothen-dieck-Witt
spaces) is homotopy cartesian if and only if the map
between, say, horizontal homotopy fibres is a homotopy equivalence, and the
map of abelian groups between horizontal cokernels of $\pi_0$'s is a
monomorphism. 

\subsection{Remark}
In theorem \ref{thm:ChgOfWkEq}, the map $GW_0(\E,w,*,\eta) \to
GW_0(\E,v,*,\eta)$ is not surjective, in general, contrary to the $K$-theory
situation in \cite[Fibration theorem 1.6.4]{wald:spaces}.
\vspace{2ex}

We will reduce the proof of theorem \ref{thm:ChgOfWkEq} to idempotent complete
exact categories with weak equivalences and duality.
We recall the relevant definitions.

\subsection{Idempotent completion}
\label{subsec:idempCompl}
Recall that the {\it idempotent completion} $\tilde{\E}$ of
an exact category $\E$ has objects pairs $(A,p)$ with $p=p^2:A \to A$ an
idempotent in $\E$. 
A map $(A,p) \to (B,q)$ is a map $f:A \to B$ in $\E$ such that $f=fp=qf$.
Composition is composition of maps in $\E$.
The idempotent completion $\tilde{\E}$ has a canonical structure of an exact
category such that the inclusion $\E \subset \tilde{\E}:A \mapsto (A,1)$ is
fully exact (see \cite[Theorem A.9.1]{TT}, where ``idempotent completion'' is
called ``Karoubianisation''). 
Any duality $(*,\eta)$ on $\E$ extends to a duality $(A,p)^*=(A^*,p^*)$ on
$\tilde{\E}$  with double dual identification $\eta_A\circ p:(A,p) \to
(A,p)^{**}$. 
If $(\E,w,*,\eta)$ is an exact category with duality, 
call a map in the idempotent completion $\tilde{\E}$ {\em weak
  equivalence} if it is a retract of a weak equivalence in $\E$.
Then $(\tilde{\E},w,*,\eta)$ is an exact category with weak equivalences and
duality.
Note that the natural inclusion $\widetilde{\E^w} \subset (\tilde{\E})^w$ is
an equivalence of categories if $(\E,w,*,\eta)$ has a (symmetric)
cone.
This is because for an object $X$ in $(\tilde{\E})^w$, 
the weak equivalence $0 \to X$ in $\tilde{\E}$ is, by
definition, a retract of a weak equivalence $f:Y \to Z$ in $\E$, and, by
functoriality, also a retract of $0 \to C(f)$, where $C(f)$ is the push-out of
$i_Y:Y \to CY$ along $f$.
Since $C(f)$ is in $\E^w$, the object $X$ is (isomorphic to an object) in
$\widetilde{\E^w}$.

\subsection{Lemma}
\label{lem:EwCofinalI}
{\it
Let $(\E,w,*,\eta)$ be an exact category with weak equivalences and strong
duality which has a symmetric cone.
Then the
commutative diagram of exact categories with duality
$$
\xymatrix{(\E^w,i) \ar[r] \ar[d] & (\tilde{\E}^w,i) \ar[d] \\
(\E,i) \ar[r] & (\tilde{\E},i)}
$$
induces a homotopy cartesian square of Grothendieck-Witt spaces.
}

\subsection*{\it Proof}
By the cofinality theorem in \cite[5.5]{myGWex}, the horizontal homotopy
fibres of 
associated Grothendieck-Witt spaces are contractible. 
Therefore, it suffices to show that the map $GW_0(\tilde{\E}^w)/GW_0(\E^w) \to
GW_0(\tilde{\E})/GW_0(\E)$ between the cokernels of horizontal $\pi_0$'s is
injective. 
For an exact category with duality $\A$, the quotient 
$GW_0(\tilde{\A})/GW_0(\A)$ is the abelian monoid of isometry classes of
symmetric spaces in $\tilde{\A}$ modulo the submonoid of symmetric spaces in
$\A$ \cite[5.2]{myGWex}.
In particular, a symmetric space in $\tilde{\A}$ yields the zero
class in $GW_0(\tilde{\A})/GW_0(\A)$ if and only if it is stably in $\A$.

Let $(A,\alpha)$ be a symmetric space in $\tilde{\E}^w$ 
whose class in $GW_0(\tilde{\E})/GW_0(\E)$ is zero.
Then there are symmetric spaces $(X,\ffi)$ and $(Y,\psi)$ in $\E$
and an isometry $(A,\alpha) \perp (Y,\psi) \cong (X,\ffi)$.
In  particular, there is an exact sequence $0 \to A \to X
\stackrel{f}{\to} Y \to 0$ in $\tilde{\E}$. 
The push-out $C(f)$ of $f$ along the $\E$-inflation $X \rightarrowtail CX$ is
in $\E$.
In the exact sequence $0 \to A \to CX \to C(f)\to 0$, we have $A$ and $CX$ in
$\tilde{\E}^w$, so that $C(f)$ is also in
$\tilde{\E}^w$, hence $C(f) \in \tilde{\E}^w\cap \E = \E^w$.
Choose $\bar{A}\in \tilde{\E}^w$ such that $A \oplus \bar{A} \in \E^w$.
Then $\bar{A}\oplus CX$ is in ${\E}^w$ as it is an extension of $A \oplus
\bar{A}$ and $C(f)$, both being in ${\E}^w$.
It follows that the symmetric spaces $(A,\alpha) \oplus \H(\bar{A} \oplus CX)$
and $\H(\bar{A} \oplus CX)$ lie in $\E^w$, 
where $\H E=(E \oplus E^*,\left(\begin{smallmatrix}0&1\\
    \eta&0\end{smallmatrix}\right))$
is the hyperbolic space associated with $E$.
This implies that $(A,\alpha)$ is trivial in $GW_0(\tilde{\E}^w)/GW_0(\E^w)$.
\Qed
\vspace{2ex}

For an exact category with weak equivalences and duality $(\E,w,*,\eta)$, the
category $\Mor\E=\Fun([1],\E)$ of morphisms in $\E$ is an exact category with
weak equivalences and duality such that the fully exact inclusion 
$\E \subset \Mor\E:E \mapsto id_E$, induced by the unique map $[1] \to [0]$, 
is duality preserving.
The inclusion factors through the fully exact subcategory 
$\Mor_w\E=\Fun_w([1],\E) \subset \Mor\E$ of weak equivalences in $\E$ and
defines a duality preserving functor
$$I:\E \to \Mor_w\E.$$
Note that $(\Mor_w\E)^w=\Mor_w(\E^w)=\Mor(\E^w)$.

The following proposition is the key to proving
Theorem \ref{thm:ChgOfWkEq}.

\subsection{Proposition}
\label{prop:E=Morw}
{\it
Let $(\E,w,*,\eta)$ be an exact category with weak equivalences and strong
duality. 
Assume that $(\E,w,*,\eta)$ has a symmetric cone.
Then the commutative diagram 
\begin{equation}
\label{eqn:prop:E=Morw}
\xymatrix{{\E}^w \ar[r] \ar[d] &  {\Mor_w}{\E}^w \ar[d]\\
           {\E} \ar[r]^I & {\Mor}_w{\E}}
\end{equation}
of duality preserving inclusions of exact categories with duality (all weak
equivalences being isomorphisms)
induces a homotopy cartesian square of associated Grothendieck-Witt spaces.
}
\vspace{2ex}

The proof uses the cone category construction of \cite[\S9]{myGWex}.
We recall the relevant definitions and facts.

\subsection{Cone exact categories}
\label{subsect:ConeEx}
Let $\A \subset \U$ be a duality preserving fully exact inclusion of
idempotent complete exact categories with duality $(*,\eta)$.
In \cite[\S9]{myGWex}, we constructed a duality preserving fully exact
inclusion 
$\Gamma: \U \subset \C(\U,\A)$ of exact categories with duality, 
depending functorially on the pair $\A \subset \U$,
such that the duality preserving commutative square
\begin{equation}
\label{eqn:ConeCat}
\xymatrix{{\A} \ar[r]^{\hspace{-2ex}\const} \ar[d] & {\C(\A,\A)} \ar[d]\\
          {\U} \ar[r]^{\hspace{-2ex}\const}        & {\C( \U,\A)}
}
\end{equation}
of exact categories with duality induces a homotopy cartesian square of
associated Grothendieck-Witt spaces, and the upper right corner $\C(\A,\A)$
(also written as $\C(\A)$) 
has contractible Grothendieck-Witt space \cite[theorem 9.9]{myGWex}.

We recall the definition of the {\em cone category} $\C(\U,\A)$, details can be
found in \cite[\S9.1-9.3]{myGWex}. 
One first constructs a category $\C_0(\U,\A)$, a localization of which
is $\C(\U,\A)$.
Objects of $\C_0(\U,\A)$ are commutative diagrams in $\U$
\begin{equation}
\label{eqn:ConeObj}
\xymatrix{
U_0 \xymono[r]^{\sim} \ar[d] & U_1 \xymono[r]^{\sim} \ar[d] & U_2
\xymono[r]^{\sim} \ar[d] & U_3 \xymono[r]^{\sim} \ar[d] & \cdots \\
U^0 & U^1 \xyepi[l]_{\sim} & U^2 \xyepi[l]_{\sim} & U^3 \xyepi[l]_{\sim}  &
\cdots \xyepi[l]_{\sim}
}
\end{equation}
such that the maps $U_i \stackrel{_{\sim}}{\rightarrowtail} U_{i+1}$ and
$U^{i+1} \stackrel{_{\sim}}{\twoheadrightarrow}  U^{i}$, $i\in \N$ are
inflations with cokernel in 
$\A$ and deflations with kernel in $\A$, respectively.
Moreover, there has to be an integer $d$ such that 
for every $i \geq j$, the map
$U_j \to U^{i+d}$ is an inflation with cokernel in $\A$ and the map $U_{i+d}
\to U^j$ is a deflation with kernel in $\A$.
If the maps in diagram (\ref{eqn:ConeObj}) are understood, we may
abbreviate the diagram as $(U_{\bullet} \to U^{\bullet})$.
Maps in $\C_0(\U,\A)$ are natural transformations of diagrams.
A sequence of diagrams in  $\C_0(\U,\A)$
 is exact if at each $U_i$, $U^j$ spot it is
an exact sequence in $\U$.
The dual of the diagram (\ref{eqn:ConeObj}) is obtained by applying the duality to the diagram: 
$(U_{\bullet}\to U^{\bullet})^*=((U^{\bullet})^*\to (U_{\bullet})^*)$.

For each diagram (\ref{eqn:ConeObj}), forgetting the upper left corner $U_0$
gives us a new object $(U_{\bullet}\to U^{\bullet})_{[1]}=(U_{\bullet +1}\to
U^{\bullet})$ and a canonical map 
$(U_{\bullet}\to U^{\bullet}) \to (U_{\bullet +1}\to U^{\bullet})$.
Similarly, forgetting the lower left corner $U^0$ defines a new object
$(U_{\bullet}\to U^{\bullet})^{[1]}= (U_{\bullet}\to U^{\bullet +1})$ and a
canonical map $(U_{\bullet}\to U^{\bullet +1}) \to
(U_{\bullet}\to U^{\bullet})$.
Finally, the category $\C(\U,\A)$ is the localization of $\C_0(\U,\A)$ with
respect to the two types of canonical maps just defined. 
A sequence in $\C(\U,\A)$ is a conflation if and only if it is isomorphic in
$\C(\U,\A)$ to the image under the localization functor $\C_0(\U,\A) \to
\C(\U,\A)$ of a conflation in $\C_0(\U,\A)$. 

There is a fully exact duality preserving inclusion $\const: \U
\subset \C(\U,\A)$ which sends an object $U$ of 
$\U$ to the constant diagram
$$
\xymatrix{
U \xymono[r]^{=} \ar[d]^{id} & U \xymono[r]^{=} \ar[d]^{id} & U
\xymono[r]^{=} \ar[d]^{id} & U \xymono[r]^{=} \ar[d]^{id} & \cdots \\
U & U \xyepi[l]_{=} & U \xyepi[l]_{=} & U \xyepi[l]_{=}  &
\cdots \xyepi[l]_{=}
}
$$
\vspace{2ex}

\subsection*{\it Proof of proposition \ref{prop:E=Morw}}
By Lemma \ref{lem:EwCofinalI}, we can (and will) assume $\E$ to be idempotent
complete. 
Then all categories in diagram (\ref{eqn:prop:E=Morw}) are
idempotent complete.

We will write $\C(\E,w)$ and $\C(\E^w)$ instead of the categories
$\C(\E,\E^w)$ and $\C(\E^w,\E^w)$ of
\ref{subsect:ConeEx}, and  
we will write $(F,\ffi)\sim (G,\psi)$ if the two non-singular exact form
functors $(F,\ffi)$, 
$(G,\psi)$ induce homotopic maps on Grothendieck-Witt spaces.
The strategy of proof is as follows.
We will extend diagram (\ref{eqn:prop:E=Morw}) to a commutative diagram of
exact categories with duality and non-singular exact form functors
\begin{equation}
\label{eqn:pf:E=Morw}
\xymatrix{{\E}^w \ar[r] \ar[d] &  {\Mor_w}{\E}^w \ar[d] \ar[r] & {\C}({\E}^w)
  \ar[d]\ar[r] & {\C}({\Mor_w}{\E}^w) \ar[d]\\
   {\E} \ar[r]^{\hspace{-1.5ex}I} & {\Mor}_w{\E} \ar[r]^{(F,\ffi)} & {\C}({\E},w)
   \ar[r]^{\hspace{-3ex}\C(I)} & {\C}({\Mor}_w{\E},w)} 
\end{equation}
where the right hand square is obtained from (\ref{eqn:prop:E=Morw}) by
functoriality of the cone category construction in \ref{subsect:ConeEx}.
We will show:
\vspace{1ex}

\hspace{-6ex}
\parbox{5in}{
\begin{itemize}
\item[(\dag)] ($GW$ applied to) the compositions
$(F,\ffi)\circ(I,id)$ and $(\C(I),id)\circ(F,\ffi)$ 
are homotopic to the constant diagram inclusions $(\const,id)$ of
\ref{subsect:ConeEx} in such a way that the homotopies restricted to $\E^w$
and $\Mor_w\E^w$ have images in the Grothendieck-Witt spaces of $\C(\E^w)$ and
$\C(\Mor_w\E^w)$, respectively.
\end{itemize}
}
\vspace{1ex}

\noindent
Assuming (\dag), the outer diagram of the left two squares and the outer diagram of the
right two squares induce homotopy cartesian diagrams of Grothendieck-Witt
spaces, by \cite[theorem 9.9]{myGWex}.
By remark \ref{rem:HgpsHcart}, 
this implies the claim of proposition \ref{prop:E=Morw}.

To construct diagram (\ref{eqn:pf:E=Morw}), we will define the
non-singular exact form functor $(F,\ffi): {\Mor}_w{\E}  \to  {\C}({\E},w)$ as
the 
composition of a non-singular exact form functor $(F_0,\ffi): \Mor_w\E \to \C_0(\E,w)$ and the
localization 
functor $\C_0(\E,w) \to \C(\E,w)$ such that its restriction to $\Mor\E^w$
has image in $\C_0(\E^w)$.
Recall that $(\E,w,*,\eta)$ is assumed to have a symmetric cone (see
definition \ref{dfn:symCone}), where $i:id \rightarrowtail C$ and $p:P
\twoheadrightarrow id$ denote the natural inflation and deflation which are
part of the structure of a symmetric cone. 
The functor $(F,\ffi)$ (or rather $(F_0,\ffi)$) sends an object $g:X
\to Y$ of $\Mor_w\E$ to the object $F(g)$ given by the diagram
\begin{equation}
\label{eqn:pf:E=Morw2}
\hspace{6ex}
\xymatrix{
X\xymono[r] \ar[d]^g & X \oplus PY \ar[d]^{\left(\begin{smallmatrix}g&
      p \\i & 0\end{smallmatrix}\right)} \xymono[r] 
& X \oplus PY \oplus C X 
\ar[d]^{\left(\begin{smallmatrix}g&p & 0 \\i & 0 & 1\\
      0&1&0\end{smallmatrix}\right)} \xymono[r] 
&X \oplus PY \oplus C X \oplus PY \ar[d] \xymono[r] & \cdots\\
 Y & Y\oplus CX \xyepi[l]
& Y\oplus CX \oplus PY \xyepi[l] & Y\oplus CX \oplus PY \oplus CX \xyepi[l]
& \cdots \xyepi[l]
}
\end{equation}
of $\C_0(\E,w)$.
In the notation $(U_{\bullet} \to U^{\bullet})$ of
\ref{subsect:ConeEx} corresponding to diagram (\ref{eqn:ConeObj}), the object
$F(g)$ is given by 

\hspace{3ex}
\renewcommand\arraystretch{1.5} 
\begin{tabular}{ll}
$U_n = X \oplus PY \oplus CX \oplus PY
\oplus CX\oplus \cdots$ & ($n+1$ summands),\\
$U^n = Y \oplus CX \oplus  PY \oplus CX \oplus
PY \oplus \cdots$& ($n+1$ summands).
\end{tabular}

\noindent
The maps $U_{n} \to U_{n+1}$ and $U^{n+1} \to U^{n}$ are the 
canonical inclusions into the first $n+1$ summands and the canonical
projections onto the first $n+1$ factors.
The maps $U_n\to U^n$ are given by the matrix 
$$(u_{rs})_{0\leq r,s \leq n}=
\left(\begin{smallmatrix}g& p  &0&0&0&\cdots\\
                            i& 0 & 1&0&&0\\
                           0&1&0&1&0&\phantom{vdots}\\
                           0&0&1&0&1&\\
                           0& &0&1&&\\
                           \vdots&0&&&&\ddots
\end{smallmatrix}\right)$$
with $u_{rs}=0$ unless $r=s=0$ or $|r-s|=1$, and $u_{0,0}=g$, $u_{0,1}=p_Y$,
$u_{1,0}=i_X$, $u_{r,r+1}=1$, $u_{r+1,r}=1$ for $r \geq 1$.
The construction of diagram (\ref{eqn:pf:E=Morw2}) is functorial in $g$,
so that $F:\Mor_w\E \to \C(\E,w)$ is indeed a functor.
The duality compatibility map $\ffi_g:F(g^*) \to (Fg)^*$ for $g:X \to Y$ 
is the identity on $X^*$ and $Y^*$, it is $\gamma_X$ on the summands
$P(X^*)$ and $\bar{\gamma}_Y$ on $C(Y^*)$.
It is clear that $F$ sends the subcategory $\Mor_w\E^w$ to the full
subcategory $\C(\E^w)$ of $\C(\E,w)$.
This defines diagram (\ref{eqn:pf:E=Morw}).

We are left with proving (\dag).
Since $U_0=X$ is the initial object of diagram (\ref{eqn:pf:E=Morw2}),
it defines a map $j=j_g: X = U_0 \to F(g)$, where  $U_0$ (and $X$) is considered
an object 
of $\C(\E,w)$ via the constant diagram embedding $\const: \E \to \C(\E,w)$. 
The map $j:X \to F(g)$ is an inflation in
$\C(\E,w)$ with cokernel in $\C(\E^w)$
because $j:X \to (U_{\bullet} \to U^{\bullet +1})$ is an inflation in
$\C_0(\E,w)$ with cokernel in $\C_0(\E^w)$.
Varying $g$, the map $j_g$ defines a natural transformation $j:\const \to F$.
Similarly, $U^0$ is the final object of diagram (\ref{eqn:pf:E=Morw2}) and
thus defines a (functorial) map
$q=q_g:F(g) \to U^0 = Y$ with kernel in $\C(\E^w)$. 
We have $g=q j$.

Let $\hat{\ffi}:F \to F^{\sharp}=*F*$ be the symmetric form on the functor $F$
associated with the duality compatibility map $\ffi$ (see
\ref{subsec:castWduals}).
The form $\hat{\ffi}:F \to F^{\sharp}$  fits into a commutative
diagram in $\C_0(\E,w)$
\begin{equation}
\label{eqn:pf:E=Morw3}
\xymatrix{
 X \ar[r]^{j_g} \ar[d]^{\eta_X} & F(g) \ar[d]^{\hat{\ffi}_g}
  \ar[r]^{q_g} &  Y \ar[d]^{\eta_Y}\\
 X^{**} \ar[r]^{q^*_{g^*}}  & F(g^*)^*
  \ar[r]^{j^*_{g^*}} &  Y^{**}.
}
\end{equation}

Write $(FI,\ffi_{FI})$ for the composition $(F,\ffi)\circ (I,id)$ of
form functors.
The natural transformation $j$ above makes the canonical
inclusion $(\const,id): \E \to \C(\E,w)$ into an admissible subfunctor
 $j: \const \subset F I$ of $FI$.
Commutativity of diagram (\ref{eqn:pf:E=Morw3}) for $g=id_X$ implies that
$\eta =j^{\sharp}\hat{\ffi}_{FI} j$, that is, $j$ defines a map $(\const,\eta)
\to (FI, \hat{\ffi}_{FI})$ of symmetric spaces associated with the form
functors $(\const, id)$ and $(FI,\ffi_{FI})$.
Since the maps $id$ and  $\ffi_{FI}$ are
isomorphisms, the symmetric space $(FI,\hat{\ffi}_{FI})$
decomposes in $\Fun(\E,\C(\E,w))$ as 
$(\const,\eta) \perp (A,\hat{\ffi}_A)$ with $(A,\hat{\ffi}_A)$ the orthogonal
complement of $(\const,\eta)$ in $(FI,\hat{\ffi}_{FI})$. 
As mentioned above, the cokernel of $j: X \to FI(X)=F(1_X)$ is in 
$\C(\E^w)$.
Therefore, the form functor
$(A,\ffi_A)$ factors through the category $\C(\E^w)$ (whose hermitian
$K$-theory space is contractible \cite[corollary 9.8]{myGWex}), so that
$(A,\ffi_A)\sim 0$. 
 Hence, 
$(FI,\ffi_{FI}) \cong (\const,id) \perp (A,\ffi_A) \sim
 (\const,id).$
By construction, the homotopy restricted to $\E^w$ has image in the
Grothendieck-Witt space of $\C(\E^w)$.
This shows the first half of the claim (\dag).

For the second half,
write $(I F,\ffi_{IF})$ and $(I F_0,\ffi_{IF_0})$ for the compositions of form
functors 
$(\C(I),id) \circ (F,\ffi)$ and $(\C(I),id) \circ
(F_0,\ffi)$, and note that 
$(IF,\ffi_{IF})$ is just the composition of $(IF_0,\ffi_{IF_0})$ with the
localization functor $\C_0(\Mor_w\E,w) \to \C(\Mor_w\E,w)$.
There is an obvious isomorphism of exact categories with duality
$\Mor_w\C_0(\E,w) \cong 
\C_0(\Mor_w\E,w)$ such that the composition $IF_0$ sends the object $(g:X \to
Y) \in \Mor_w\E$ to $id_{F(g)}: F(g) \to F(g)$, and the duality compatibility
morphism $\ffi_{IF_0}$ becomes
$(\ffi_g,\ffi_g): id_{F(g^*)} \to id_{F(g)^*}$.
Consider the functorial bicartesian square
$$
\xymatrix{
( X \stackrel{j}{\to} F(g) ) \ar[r]^{(j,1)} \ar[d]^{(1,q)} &
(F(g) \stackrel{1}{\to} F(g) ) \ar[d]^{(1,q)}\\
( X \stackrel{g}{\to}  Y ) \ar[r]^{(j,1)} &
(F(g) \stackrel{q}{\to}  Y) 
}$$
in $\Mor_w\C_0(\E,w)$.
The total complex of the square (considered as a bicomplex) is a conflation in
$\Mor_w\C_0(\E,w) = \C_0(\Mor_w,w)$.
It is therefore also a conflation in $\C(\Mor_w,w)$, hence the square
is also bicartesian in $\C(\Mor_w,w)$.
In $\C(\Mor_w,w)$, the horizontal maps in the square are inflations 
with cokernel in $\C(\Mor_w\E^w)$
since $(j,1)$ is isomorphic to the $\C_0(\Mor_w,w)$-inflation $(j,1)^{[1]}$
which has cokernel in $\C_0(\Mor_w\E^w)$.
Similarly, the vertical maps in the square are deflations in $\C(\Mor_w,w)$
with kernel in $\C(\Mor_w\E^w)$
since $(1,q)$ is isomorphic to the $\C_0(\Mor_w,w)$-deflation $(1,q)_{[1]}$
with kernel in $\C_0(\Mor_w\E^w)$.

Commutativity of diagram (\ref{eqn:pf:E=Morw3}) implies that the form
$\hat{\ffi}_{IF_0}$ on the upper right corner $IF_0$ of the square
extends to 
a form on the whole bicartesian square such that its restriction to the lower
left corner is the constant diagram inclusion $(\const, \eta): \Mor_w\E \to
\C_0(\Mor_w\E,w)$. 
It follows that $\ker(1,q) \subset IF$ is a totally isotropic subfunctor of
$(IF,\ffi_{IF})$ with induced form on
$(X\stackrel{g}{\to}Y)=\ker(1,q)^{\perp}/\ker(1,q)$  
isometric to the constant diagram inclusion $(\const,
id): \Mor_w\E \to \C(\Mor_w\E,w)$. 
By the additivity theorem \cite[theorem 7.1]{myGWex} (or its generalization
in theorem \ref{thm:addtyFun}),
the form functors $(\const, id) \perp
\H\ker(1,q)$ and $(IF,\ffi_{IF})$ induce homotopic maps
on Grothendieck-Witt spaces.
Since $\H\ker(1,q)$ has image in $\C(\Mor_w\E^w)$ whose
Grothendieck-Witt space is contractible, we have
$(IF,\ffi_{IF}) \sim (\const,id) \perp \H\ker(1,q) \sim
(\const,id)$.
The homotopies restricted to $\Mor_w\E^w$ have image in the Grothendieck-Witt
space of $\C(\Mor_w\E^w)$.
This is clear for the second homotopy, and for the first, it follows from the
proof of additivity in \ref{pfofAddtyFun}.
\Qed
\vspace{2ex}

Next, we prove a variant of the Change-of-weak-equivalence theorem.

\subsection{Proposition}
\label{prop:iCOW}
{\it
Let $(\E,w,*,\eta)$ be an exact category with weak equivalences and strong
duality which has a symmetric cone.
Then the 
following commutative diagram of exact categories with weak equivalences and
duality induces a homotopy 
cartesian square of Grothendieck-Witt spaces with contractible upper right corner
\begin{equation}
\label{eqn:prop:Change1}
\xymatrix{({\E}^w,i) \ar[r] \ar[d] &  ({\E}^w,w) \ar[d]\\
           ({\E},i) \ar[r] & ({\E},w).}
\end{equation}
}

\subsection*{\it Proof}
From lemma \ref{lem:htpyOfFunctors} it is clear that $GW({\E}^w,w)$
is contractible since $0 \to id$ is a natural weak equivalence in
$({\E}^w,w)$. 
Consider the commutative diagram of (simplicial) exact categories with
dualities (all weak equivalences being isomorphisms)
$$\xymatrix{
{\E}^w \ar[r] \ar[d] & \Fun_w(\underline{0},{{\E}}^w) \ar[d] \ar[r] &
\left\lgroup n \mapsto \Fun_w(\underline{n},{{\E}}^w) \right\rgroup \ar[d]\\
{\E} \ar[r]  & \Fun_w(\underline{0},{{\E}})  \ar[r] &
\left\lgroup n \mapsto \Fun_w(\underline{n},{{\E}}) \right\rgroup }$$
in which the left hand square can be identified with the square of proposition
\ref{prop:E=Morw} and induces therefore a homotopy cartesian square of
Grothendieck-Witt spaces.
On Grothendieck-Witt spaces, the right vertical map can be identified with the
map $(\E^w,w) \to (\E,w)$ in view of the simplicial resolution lemma
\ref{lem:simplResLem}. 
The right hand square is the inclusion of degree zero simplices.
The proof of proposition \ref{prop:iCOW} is thus reduced to showing
that the right hand square of the diagram induces a homotopy cartesian square
of Grothendieck-Witt spaces.

Let $\Fun^1_w(\underline{n},\E) \subset \Fun_w(\underline{n},\E)$ be the full
subcategory of those functors $A:\underline{n} \to \E$ for which $A_p \to A_q$
is an inflation, and $A_{q'} \to A_{p'}$ is a deflation, $0 \leq p \leqq q
\leq n$.
It inherits the structure of an exact category with duality from 
$\Fun_w(\underline{n},\E)$.
Further, let $\Fun^0_w(\underline{n},\E)$ be the category which is equivalent
to $\Fun^1_w(\underline{n},\E)$ but where an object is an object $A$ of
$\Fun^1_w(\underline{n},\E)$ together with a choice of subquotients $A_{p,q} =
A_q/A_p = \coker(A_p \stackrel{\sim}{\rightarrowtail} A_q) \in \E^w$ and
induced maps $A_{p,q} \to A_{p,q}$,
and together with a choice of kernels $A_{q',p'} = \ker(A_{q'}
\stackrel{\sim}{\twoheadrightarrow} A_{p'})\in \E^w$ for $0 \leq p \leq q\leq
n$.
The category $\Fun^0_w(\underline{n},\E)$ is an exact category with duality
such that the forgetful functor $\Fun^0_w(\underline{n},\E) \to
\Fun^1_w(\underline{n},\E)$ is an equivalence of exact categories with duality.
We have an exact functor 
$\Fun^0_w(\underline{n},\E) \to S_n\E^w: A \mapsto (A_{p,q})_{0\leq p\leq q
  \leq n}$.
By the additivity theorem \cite[theorem 7.1]{myGWex} (or
\ref{thm:addtyFun}), this functor induces a map which is part of a 
split homotopy fibration
$$GW(\Fun_w(\underline{0},\E)) \to GW(\Fun^0_w(\underline{n},\E)) \to K(S_n\E^w).$$
The same argument applies to $(\E^w,w)$ instead of $(\E,w)$ so that, varying
$n$, we obtain a map of homotopy fibrations after topological realization
$$\xymatrix{
 GW\Fun_w(\underline{0},{{\E}}^w) \ar[d] \ar[r] &
| n \mapsto GW\Fun_w^0(\underline{n},{{\E}}^w)| \ar[d] \ar[r] & |n \mapsto
K(S_n\E^w)| \ar@{=}[d]\\
GW\Fun_w(\underline{0},{{\E}})  \ar[r] &
| n \mapsto GW\Fun_w^0(\underline{n},{{\E}})| \ar[r] & |n \mapsto
K(S_n\E^w)| }$$
which shows that the left square is homotopy cartesian.

Since $\Fun_w^0 \to \Fun_w^1$ is an equivalence of exact categories with
duality, the proposition follows once we show that the inclusion 
$I: \Fun_w^1(\underline{n},\A) \subset \Fun_w(\underline{n},\A)$ induces a
homotopy equivalence on Grothendieck-Witt spaces for $\A=\E, \E^w$.
We illustrate the argument for $\A=\E$ and $n=1$.
The general case is mutatis mutandis the same.

We define two functors $F,G:\Fun_w(\underline{1},\E) \to
\Fun_w^1(\underline{1},\E)$.
The functor $F$ sends $E_{1'} \stackrel{\sim}{\to}E_{0'}
\stackrel{\sim}{\to}E_{0}  \stackrel{\sim}{\to}  E_1$ to 
$E_{1'}\oplus PE_{0'} \stackrel{\sim}{\twoheadrightarrow}E_{0'}
\stackrel{\sim}{\to}E_{0}  \stackrel{\sim}{\rightarrowtail} E_1 \oplus CE_0$,
the functor $G$ sends the same object to 
$ PE_{0'} \stackrel{\sim}{\twoheadrightarrow} 0
\stackrel{\sim}{\to} 0  \stackrel{\sim}{\rightarrowtail} CE_0$.
Both functors $F$ and $G$ are equipped with canonical duality compatibility
morphisms, induced by $\gamma$ and $\bar{\gamma}$ from definition
\ref{dfn:symCone}, such that $F$ and $G$ are non-singular exact form
functors. 
By the additivity theorem, we have $IF \sim id \perp IG$ and $FI \sim id \perp
GI$, so that $GW(F)-GW(G)$ defines an inverse of $GW(I)$, up to homotopy.
\Qed
\vspace{2ex}

\subsection*{\it Proof of theorem \ref{thm:ChgOfWkEq}}
By lemma \ref{lem:htpyOfFunctors}, $GW({\E}^w,w)$ is contractible.
Let $\A=\E^{\str}_w$ be the strictification of $\E$ from lemma
\ref{lem:strictification}, and recall that it has a strong duality.
Consider the commutative diagram of exact categories
with weak equivalences and duality
$$\xymatrix{
(\A^w,i) \ar[r]\ar[d] & (\A^w,w)\ar[d] & \\
(\A^v,i) \ar[r]\ar[d] & (\A^v,w) \ar[r]\ar[d] & (\A^v,v)\ar[d]\\
(\A,i) \ar[r] & (\A,w) \ar[r] & (\A,v).}$$
By the strictification lemma \ref{lem:strictification} and lemma
\ref{lem:htpyOfFunctors}, the square in theorem \ref{thm:ChgOfWkEq} is
equivalent to the lower right square in the diagram.
By proposition \ref{prop:iCOW}, the upper square and the outer diagram of the
left two squares are homotopy cartesian in hermitian $K$-theory.
Since the left vertical maps are surjective on $GW_0$ (because
$\A=\E^{\str}_w$), it follows that the 
lower left square is homotopy cartesian in hermitian $K$-theory, by remark
\ref{rem:HgpsHcart}. 
Again, by proposition \ref{prop:iCOW}, the outer diagram of the two lower
squares induces a homotopy cartesian square of Grothendieck-Witt spaces.
Together with the facts that the lower left square is homotopy cartesian in
hermitian $K$-theory and that the lower left vertical map is surjective on
$GW_0$, this implies that the lower right square induces a homotopy cartesian
square of Grothendieck-Witt spaces.
\Qed

\subsection{Theorem \rm (Cofinality)}
\label{thm:cofinal}
{\it
Let $(\E,w,*,\eta)$ be an exact category with weak equivalences and 
duality which has a symmetric cone.
Let $A \subset K_0(\E,w)$ be a subgroup closed under the duality action on
$K_0(\E,w)$, and let $\E_A \subset \E$ be the full 
subcategory of those objects whose class in $K_0(\E,w)$ belongs to $A$.
Then the category $\E_A$ inherits the structure of an exact category with weak
equivalences and duality from $(\E,w,*,\eta)$, and
the induced map on Grothendieck Witt spaces
$$GW(\E_A,w,*,\eta) \longrightarrow GW(\E,w,*,\eta)$$ is an isomorphism on
$\pi_i$, $i\geq 1$, and a monomorphism 
on $\pi_0$.
}

\subsection*{\it Proof}
Let $\U=\E^{\str}_w$, and consider the 
the diagram of exact categories with weak equivalences and duality
$$\xymatrix{
(\U_A^w,i) \ar@{=}[d] \ar[r] & (\U_A,i) \ar[d] \ar[r] & (\U_A,w) \ar[d]\\
(\U^w,i) \ar[r] & (\U,i) \ar[r] & (\U,w).}$$
On Grothendieck-Witt spaces, the right vertical map can be identified (up to
homotopy) with the map in the theorem, by lemmas \ref{lem:strictification}
and \ref{lem:htpyOfFunctors}. 
The rows are homotopy fibrations, by proposition \ref{prop:iCOW}, and the
right horizontal maps are surjective on $GW_0$ (as $\U=\E^{\str}_w$).
It follows that the right square induces a homotopy cartesian square of
Grothendieck-Witt spaces.
Since $\U_A \subset \U$ is a cofinal inclusion of exact categories with
duality, the cofinality theorem of \cite[corollary 5.5]{myGWex} shows that
the homotopy fibre 
of the Grothendieck-Witt spaces of the middle vertical map is contractible.
As the right square is homotopy cartesian in hermitian $K$-theory, the same is
true for the right vertical map.
\Qed

\section{Approximation, Change of exact structure and Resolution}

In this section we prove in theorems \ref{thm:ApprFunctorialCalcFrac} and
\ref{thm:ApprCalcOfFrac2} variants of Waldhausen's approximation theorem
\cite[Theorem 1.6.7]{wald:spaces} which hold for higher
Grothendieck-Witt groups.
We explain two immediate consequences, one concerning the conditions under
which a change of exact structure 
has no effect on Grothendieck-Witt groups (lemma \ref{lem:ChExStr},
compare \cite[Theorem 1.9.2]{TT}) and the other concerning an 
analog of Quillen's resolution theorem (lemma \ref{lem:resolution},
compare \cite[\S 4 Corollary 1]{quillen:higherI}). 

\subsection{Theorem \rm (Approximation I)}
\label{thm:ApprFunctorialCalcFrac}
{\it
Let $(F,\ffi):\A \to \B$ be a non-singular exact form functor between exact
categories with weak equivalences and duality.
Assume the following.
\begin{itemize}
\item[(a)]
Every map in $\A$ can be written as the composition of an inflation followed
by a weak equivalence.
\item[(b)]
  A map in $\A$ is a weak equivalence iff its image in $\B$ is a weak
  equivalence.
\item[(c)]
For every map $f: FA \to E$ in $\B$, there is a map $a:A \to B$ in $\A$ and 
a weak equivalence $g: FB \stackrel{_{\sim}}{\to} E$ in $\B$ such that
$f=g\circ Fa$: 
$$\xymatrix{FA \ar@{-->}[r]^{\forall f} \ar[d]_{Fa} & E\\
  \exists\phantom{a} FB \ar@{-->}[ur]_{\sim} &}$$
\item[(d)]
The duality compatibility morphism $\ffi_A: F(A^*) \to F(A)^*$ is an
isomorphism for every $A$ in $\A$.
\item[(e)]
For every map $f:FA \to FB$ in $\B$, there is an exact functor $L:\A \to \A$, 
a natural weak equivalence $\lambda: L \stackrel{_{\sim}}{\to} id_{\A}$ and a map
$a: LA \to B$ in $\A$ such that $Fa=f\circ F\lambda_A$:
$$\xymatrix{ \exists\phantom{a} FLA \ar[r]^{\sim} \ar@/_1pc/[rr]_{Fa}
& FA \ar@{-->}[r]^{\hspace{-3ex} \forall f}  & FB.}
$$
\item[(f)]
For every map $a:A \to B$ in $\A$ such that $Fa=0$ in $\B$, there is an
exact functor $L:\A \to \A$ and a natural weak equivalence $\lambda:L
\stackrel{_{\sim}}{\to} id_{\A}$ such that $a\circ \lambda_A=0$ in $\A$: 
$$\xymatrix{ \exists\phantom{a} LA \ar[r]^{\sim} \ar@/_1pc/[rr]_{0}
& A \ar[r]^{\hspace{-1ex} \forall a}  &  B, &Fa=0.}
$$
\end{itemize}
Then $(F,\ffi)$ induces homotopy equivalences
$$(wS^e_{\bullet}\A)_h
\stackrel{\sim}{\longrightarrow}(wS^e_{\bullet}\B)_h\hspace{5ex}{\rm \it and}
\hspace{5ex} 
GW(\A,w,*,\eta) \stackrel{\sim}{\longrightarrow} GW(\B,w,*,\eta) .$$
}

\subsection*{\it Proof}
For the purpose of the proof, we call {\it lattice} a pair $(L,\lambda)$ with
$L: \A \to \A$ an exact functor and $\lambda:L \stackrel{_{\sim}}{\to} id_{\A}$ a natural weak
equivalence.
Lattices form an associative monoid under composition
$$(L_2,\lambda_2) \circ (L_1,\lambda_1) := (L_2L_1, \lambda_2 \circ
L_2(\lambda_1)).$$ 
Since $\lambda_2\circ L_2(\lambda_1)=\lambda_1\circ \lambda_{2,L_1}$
(as $\lambda_2$ is a natural transformation),
lattices behave like a multiplicative set in a commutative ring. 
More precisely, composition of lattices allows us to generalize properties
(e) and (f) to finite families of maps: 
\begin{itemize}
\item[(e')]
for any finite set of maps $f_i:FA_i\to FB_i$ in $\B$, $i=1,...,n$, there
are a lattice $(L,\lambda)$ and maps $a_i:LA_i \to B_i$ such that
$Fa_i=f_i\circ F\lambda_{A_i}$, $i=1,...,n$, and 
\item[(f')]
for any finite set of maps $a_i:A_i \to B_i$ in $\A$ such that $Fa_i=0$ in
$\B$, $i=1,...,n$, there is a lattice $(L,\lambda)$ such that
$a_i\lambda_{A_i}=0$ in $\A$, $i=1,...,n$.
\end{itemize}
We will refer to (e') and
(f') as ``clearing denominators'', in analogy with
the localization of a commutative ring with respect to a multiplicative subset.
``Clearing denominators'' together with (b)
and (d) implies that we can lift non-degenerate
symmetric forms from $\B$ to $\A$ in the following sense:
\begin{itemize}
\item[(\dag)]
For any non-degenerate symmetric form $(FA,\alpha) \in (w\B)_h$ on the image
$FA$ of an 
object $A$ of $\A$, there is a lattice $(L,\lambda)$ and a non-degenerate
symmetric form $(LA,\beta) \in (w\A)_h$ on $LA$ such that $F(\lambda_A)$ is
a map of symmetric spaces $F(\lambda_A):F(LA,\beta) \stackrel{_{\sim}}{\to}
(FA,\alpha)$. 
\end{itemize}
The proof is the same as the classical proof which shows that a non-degenerate
symmetric form over the fraction field of a Dedekind domain can be lifted to a
(usual) lattice in the ring. 
In detail, the map $\ffi_A^{-1}\alpha:FA \to F(A^*)$ lifts to a map $a_1:L_1A
\to A^*$ such that $\ffi_A Fa_1=\alpha F\lambda_{1,A}$ for some lattice
$(L_1,\lambda_1)$.
The map $a:\lambda_{1,A}^*a_1:L_1A \to (L_1A)^*$ is a weak equivalence but not
necessarily symmetric.
However, the difference $\delta = a - a^*\eta_{L_1A}$ satisfies $F\delta=0$, so that
there is a second lattice $(L_2,\lambda_2)$ such that $\delta\lambda_{2,L_1A}=0$.
Then $\beta=\lambda_{2,L_1A}^*a\lambda_{2,L_1A}:L_2L_1A \to (L_2L_1A)^*$ is a
non-degenerate symmetric form on $LA$ with
$(L,\lambda)=(L_2,\lambda_2)\circ(L_1,\lambda_1)$, and 
$F(\lambda_A)$ defines a map of symmetric spaces
$F(LA,\beta)=(FLA,\ffi_{LA}F\beta) \stackrel{_{\sim}}{\to} (FA,\alpha)$.

Apart from ``clearing denominators'', the proof of the first homotopy
equivalence in the theorem proceeds now as
the proof of \cite[theorem 10]{myMZ} which was based on the proof of 
\cite[theorem 1.6.7]{wald:spaces}.
We first note that under the assumptions of the theorem, the non-singular
exact form functors
$S_n(F,\ffi):S_n\A \to S_n\B$ also satisfy
(a) - (f), $n\in \N$.
For (a) - (c), this is
in \cite[lemma 1.6.6]{wald:spaces} (using the fact that in the presence of
(a), the map $a:A \to B$ in
(c) can be replaced by an inflation), 
(d) extends by
functoriality, and the extension of (e), (f) to $S_n$ easily follows by
induction on $n$ by successively clearing denominators.

In order to show that $(wS^e_{\bullet}\A)_h \to (wS^e_{\bullet} \B)_h$ is a
homotopy equivalence, 
it suffices to prove that $(wS^e_n\A)_h \to (wS^e_n\B)_h$ is a homotopy
equivalence for every $n\in \N$, which, by the argument of the previous
paragraph, only needs 
to be checked for $n=0$, that is, it is sufficient to
prove that 
$$F:(w\A)_h \to (w\B)_h$$
is a homotopy equivalence.
The last claim will follow from Quillen's theorem A once we show that for
every object 
$X=(X,\psi)$ of $(w\B)_h$, the comma category $(F\downarrow X)$ is non-empty and
contractible.

By (a) with $A=0$, there is an object $B$ of $\A$
and a weak equivalence $FB\stackrel{_{\sim}}{\to} X$, hence a map
$(FB,\psi_{|FB}) \to (X,\psi)$ in $(w\B)_h$. 
By (\dag) above, there is a symmetric space $(C,\gamma)$ in $\A$ and a map
$F(C,\gamma) \to (FB,\psi_{|FB})$ in $(w\B)_h$.
Hence, the category $(F\downarrow X)$ is non-empty.

In order to show that $(F\downarrow X)$ is contractible it suffices to show
that every functor $\cP \to (F\downarrow X)$ from a finite poset $\cP$ to the
comma category $(F\downarrow X)$ is
homotopic to a constant map (see for instance \cite[Lemma 14]{myMZ}).
Such a functor is given by a triple $(\underline{A}, \alpha, f)$ where
$\underline{A}$ is a functor
$\underline{A}: \cP \to w\A:i\mapsto A_i, (i\leq
j) \mapsto a_{j,i}$ together with a collection 
$\alpha$ of non-degenerate symmetric forms $\alpha_i:A_i \to
A_i^*$ in $\A$, $i \in \cP$, such that $\alpha_i=\alpha_{j|A_i}$ whenever
$i\leq j$ in $\cP$, and $f$ is
a collection of compatible maps of symmetric spaces
$f_i:F(A_i,\alpha_i) \to (X,\psi)$ in $(w\B)_h$ such that $f_i=f_jFa_{j,i}$
whenever $i\leq j \in \cP$. 
It is convenient to consider
 $f$ as a map $F(\underline{A},\alpha) \to (X,\psi)$ of functors $\cP \to
 (w\B)_h$, where objects in $(w\B)_h$ (or in $\A$, $\B$, $(w\A)_h$) such
 as $(X,\psi)$ are interpreted as constant $\cP$-diagrams. 
By \cite[Lemma 13]{myMZ}, there is a map $b: \underline{B} \to \underline{A}$
of functors $\cP \to w\A$ such that $b_i:B_i \to A_i$ is a weak equivalence,
$i\in \cP$, and the $\cP$-diagram $\underline{B}:\cP \to \A$ has a colimit
in $\A$ such that $F(\underline{B}) \to F(\colim_{\cP}\underline{B})$
represents the 
colimit of $F(\underline{B})$ in $\B$ (the diagram $\underline{B}$ is a
cofibrant replacement of $\underline{A}$ in a suitable cofibration structure
on the category of functors $\cP \to \A$, see \cite[Appendix A.2]{myMZ},
cofibrant objects have colimits, and $F$, being an exact functor, preserves
cofibrant objects and their colimits). 
By (c), the natural map
$F(\colim_{\cP}\underline{B}) = \colim_{\cP}F(\underline{B}) \to X$ induced
by $b\circ f$ factors as $F(\colim_{\cP}\underline{B}) \to F(C)
\stackrel{c}{\to} X$ where 
the first map is in the image of $F$ and the second map is a
weak equivalence in $\B$.
Let $g:\underline{B} \to C$ be the composition $\underline{B} \to
\colim_{\cP}\underline{B} \to C$ which, by (b), is
a weak equivalence since $b$, $f$, and $c$ are.
The null-homotopy to be constructed can be read off the following diagram
\begin{equation}
\label{eqn:approxDiag}
\xymatrix{
L'L{\underline{B}},\beta
\ar@{}[d]|-{\begin{smallmatrix}\beta= \alpha_{|L'LB} = \gamma_{|L'LB}
  \end{smallmatrix}} 
\ar[r]^{\lambda'} &  
L{\underline{B}} \ar[r]^{\lambda} \ar[d]^{Lg} & {\underline{B}} \ar[r]^{b}
\ar[d]^g & 
{\underline{A}, \alpha} \ar@{-->}[d]^{f}\\
& LC,\gamma \ar[r]^{\lambda} & C \ar@{-->}[r]^{c} & X,{\psi}
}
\end{equation}
of which we have constructed the right hand square, so far.
In the diagram, a dashed arrow $A \dashrightarrow X$ stands for an arrow $FA
\to X$ in $\B$, and solid arrows are arrows in $\A$.
By (\dag), there is a lattice $(L,\lambda)$ and a non-degenerate symmetric
form 
$(LC,\gamma)\in (w\A)_h$  such that $F\lambda_C:F(LC,\gamma) \to
(FC,\psi_{|FC})$ defines a map in $(w\B)_h$.
The restrictions of $\gamma$ and
$\alpha_i$ to $LB_i$ may not coincide, but their images under $F$ coincide
since in $\B$, both are (up to composition with the isomorphism
$\ffi_{B_i}$) the restriction of $\psi$ to $FLB_i$. 
Clearing denominators, we can find a lattice $(L',\lambda')$ such that 
$\gamma_{|L'LB_i}=\alpha_{i|L'LB_i}=:\beta_i$ for all $i\in \cP$ simultaneously.
The outer part of the diagram involving $(L'L\underline{B},\beta)$,
$(\underline{A},\alpha)$, $(LC,\gamma)$, $(X,\psi)$ and the maps between
them, defines 
a homotopy between $(\underline{A},\alpha,f)$ and the constant functor
$(LC,\gamma,cF\lambda):\cP \to (F\downarrow X)$ via
the functor $(L'L\underline{B},\beta,fF(b\lambda\lambda')):\cP \to (F\downarrow
X)$. 
This finishes the proof of the first homotopy equivalence in the theorem.

The same proof (forgetting forms), or \cite[1.6.7]{wald:spaces}, implies that
$wS_{\bullet}\A \to wS_{\bullet}\B$ is a homotopy equivalence.
Therefore, the map $GW(\A,w,*,\eta) \to GW(\B,w,*,\eta)$ is also a homotopy
equivalence. 
\Qed

\subsection{Change of exact structure}
\label{subsec:ChExStr}
Let $(\A,w,*,\eta)$ be an additive category with weak equivalences and duality,
and assume that $\A$ can be equipped with two exact structures $\textgoth{E}_1$
and $\textgoth{E}_2$ one smaller than the other, $\textgoth{E}_1 \subset
\textgoth{E}_2$, so that the identity functor $(\A,\textgoth{E}_1) \to
(\A,\textgoth{E}_2)$ is exact.
Assume furthermore that the duality functor $*:\A^{op} \to \A$ is exact for
both exact structures $\textgoth{E}_2$ and $\textgoth{E}_2$, so that the
identity defines a duality preserving exact functor 
\begin{equation}
\label{eqn:rmkChgExStr}
(\A,\textgoth{E}_1,w,*,\eta) \to (\A,\textgoth{E}_2,w,*,\eta).
\end{equation}
The following is an immediate consequence of theorem
\ref{thm:ApprFunctorialCalcFrac}. 

\subsection{Lemma}
\label{lem:ChExStr}
{\it
In the situation of \ref{subsec:ChExStr}, 
if every map in $\A$ can be written as the composition of an inflation in
$\textgoth{E}_1$ followed by a weak equivalence, then the map
(\ref{eqn:rmkChgExStr}) induces an equivalence of hermitian
$S_{\bullet}$-constructions and of associated Grothendieck-Witt spaces.
}
\vspace{-2.5ex}

\Qed
\vspace{2ex}

The purpose of the next lemma (lemma \ref{lem:C'replacesC} below) is to
simplify some of the hypothesis of 
theorem \ref{thm:ApprFunctorialCalcFrac} provided the exact categories with
weak equivalences in the theorem are ``categories of complexes''.

\subsection{Definition}
\label{dfn:CatOfCxs}
We call an exact category with weak equivalences $(\C,w)$ a {\it category of
complexes} (with underlying additive category $\C_0$)
if 
$\C \subset \Ch\C_0$ is a full additive subcategory of the category
$\Ch\C_0$ of chain
complexes $\C_0$ such that
\begin{enumerate}
\item
$\C$ is closed under degree-wise split extensions in $\Ch\C_0$,
\item
degree-wise split exact sequences are exact in $\C$,
\item
with a complex $A$ in $\C$ its usual cone $CA$ (that is, the cone on the
identity map of $A$) and all its shifts $A[i]$, $i\in \Z$ are in $\C$, and
\item
the set of weak equivalences $w$ contains at least all usual homotopy
equivalences between complexes in $\C$. 
\end{enumerate}

\subsection{Lemma}
\label{lem:C'replacesC}
{\it
Let $(\A,w)$ and $(\B,w)$ 
be categories of complexes with associated additive categories  
$\A_0$ and $\B_0$, and let $F:(\A,w) \to (\B,w)$ be an exact functor which is
the induced functor on 
chain complexes of an additive functor $\A_0 \to \B_0$.
Then condition \ref{thm:ApprFunctorialCalcFrac}(a)
holds and condition \ref{thm:ApprFunctorialCalcFrac}  (c)
is implied by conditions \ref{thm:ApprFunctorialCalcFrac}
(b), (e)
and the following condition. 
\begin{itemize}
\item[(c')]
For every object $E$ of $\B$ there is an object $A$ of $\A$ and a weak
equivalence
$$FA \stackrel{\sim}{\longrightarrow}  E.$$
\end{itemize}
}

\subsection*{\it Proof}
Every map $f:A \to B$ in $\A$ can be written as the composition
$$\xymatrix{A \xymono[r]^{\hspace{-2ex}\left(\begin{smallmatrix}f\\
      i\end{smallmatrix}\right)} &  B \oplus CA  
\ar[r]^{\hspace{2ex}\left(\begin{smallmatrix}1&0
      \end{smallmatrix}\right)}_{\hspace{2ex}\sim} & B, }$$
where $i:A \rightarrowtail CA$ is the canonical inclusion of $A$ into its cone
$CA$. 
This shows that condition \ref{thm:ApprFunctorialCalcFrac}
(a) holds.

To prove condition \ref{thm:ApprFunctorialCalcFrac}  (c)
assuming conditions \ref{thm:ApprFunctorialCalcFrac}  (b),
(e) 
and (c') are satisfied, let $f:FA \to E$ be a map in $\B$ with $A$ in $\A$.
By (c'), there is a weak equivalence $s:F(B') \to E$ in $\B$ with $B'$ in $\A$.
Let $M$ be the pull-back 
of $f:FA \to E$ along the degree-wise split surjection $(s,p): F(B')\oplus PE
\to E$ where $p: PE\to E$ is the canonical degree-wise split surjection from
the 
contractible complex $PE=CE[-1]$ to $E$, so that we have a homotopy
commutative diagram 
$$\xymatrix{M \ar[r] \ar[d] 
& F(A) \ar[d]^f\\
            F(B') \ar[r]_{\hspace{1ex}\sim}^s & E.}$$ 
Since $F(B')\oplus PE \to E$ is a deflation and a weak equivalence, its
pull-back, the map
$M \to FA$, is a deflation and a weak equivalence, too.
By condition (c'), there is a weak equivalence $F(A') \to M$ with $A'$ in
$\A$, and by condition 
\ref{thm:ApprFunctorialCalcFrac}  (e)
we can assume the compositions $F(A') \to FA$ and $F(A')
\to F(B')$ to be the images $F\alpha$ and $F\beta$ of maps $\alpha:A' \to A$
and $\beta:A' \to B$ in $\A$. 
The resulting square involving 
$F(A')$, $F(B')$, $FA$ and $E$ homotopy commutes, and the map $\alpha:A' \to
A$ is 
a weak equivalence, by \ref{thm:ApprFunctorialCalcFrac} (b).
Replacing $B'$ with $B'\oplus CA'$, we obtain a commutative diagram
$$\xymatrix{
F(A') \xymono[r]^{F\alpha}_{\sim}
\ar[d]_{\left(\begin{smallmatrix}F\beta\\ Fi\end{smallmatrix}\right)} & F(A)
\ar[d]^f \\  
 F(B') \oplus F(CA')
 \ar[r]^{\hspace{7ex}\left(\begin{smallmatrix}s&g\end{smallmatrix}\right)} &
 E,}$$ 
where $i:A' \rightarrowtail CA'$ is the canonical inclusion of $A'$ into its
cone, $g: F(CA')=CF(A') \to E$ a map such that $g\circ Fi$ is the
null-homotopic map $f\circ F\alpha - s\circ F\beta: F(A') \to E$.
Let $B$ be the push-out of $\alpha: A'\to A$ along the degree-wise split
inclusion 
$A' \rightarrowtail B' \oplus CA'$ and call $a:A \to B$ and $b:
B'\oplus CA' \to B$ the induced maps.
Since $\alpha:A' \to A$ is a weak equivalence, so is $b$. 
The functor $F$ preserves push-outs along inflations so that we obtain an
induced map $t:F(B) \to E$ which is a weak equivalence since $F(b)$ and
$(s,g)$ are.
By construction, we have $t\circ Fa = f$.
\Qed
\vspace{2ex}

Next, we prove an analog of Quillen's resolution theorem
\cite[\S4]{quillen:higherI}.

\subsection{Lemma \rm (Resolution)}
\label{lem:resolution}
{\it
Let $(\B,w,*,\eta)$ be an exact category with weak equivalences and duality
such that $(\B,w)$ is a category of complexes.
Let $\A\subset \B$ be a full subcategory closed under the duality, degree-wise
split extensions and under taking cones and shifts in $\B$.
Restricting $(w,*,\eta)$ to $\A$ makes $\A$ into an exact category with weak
equivalences and duality.
Assume that the following resolution condition holds:
\vspace{1ex}

For every object
$E$ of $\B$, there is an object $A$ of $\A$ and a weak equivalence
$$A \stackrel{\sim}{\longrightarrow} E.$$

\noindent
Then the duality preserving inclusion $\A \subset \B$ induces homotopy
equivalences
$$(wS^e_{\bullet}\A)_h
\stackrel{\sim}{\longrightarrow}(wS^e_{\bullet}\B)_h\hspace{5ex}{\rm \it and}
\hspace{5ex} 
GW(\A,w,*,\eta) \stackrel{\sim}{\longrightarrow} GW(\B,w,*,\eta) .$$
}

\subsection{\it Proof}
This follows from theorem \ref{thm:ApprFunctorialCalcFrac} in view of lemma
\ref{lem:C'replacesC} and the fact that conditions \ref{thm:ApprFunctorialCalcFrac}
(e) and (g) trivially hold
since $\A \subset \B$ is fully faithful, and condition
(d) holds since $\A \to \B$ is duality preserving.
\Qed
\vspace{2ex}

We finish the section with theorem \ref{thm:ApprCalcOfFrac2} below.
It is a variant of theorem
\ref{thm:ApprFunctorialCalcFrac} 
when conditions \ref{thm:ApprFunctorialCalcFrac} (e) and
(f) can not be achieved in a functorial way but the
form functor is a localization by a calculus of
right fractions.

\subsection{Definition}
\label{dfn:locnCalcFrac}
Call a functor $F: \A \to \B$ between small categories a {\em localization
by a calculus of right fractions} if the following three conditions
hold (compare theorem \ref{thm:ApprFunctorialCalcFrac} (e), (f)).
\begin{enumerate}
\item
\label{enum:dfn:locnCalcFrac1}
The functor $F:\A \to \B$ is essentially surjective.
\item
\label{enum:dfn:locnCalcFrac2}
For every map $f:F(A) \to F(B)$ between the images of objects
$A$, $B$ of $\A$, there are maps $s:A' \to A$ and $g:A' \to B$ in $\A$ with
$F(s)$ an isomorphism in $\B$ and $f\circ F(s) = F(g)$.
\item
\label{enum:dfn:locnCalcFrac3}
For any two maps $a,b:A \to B$ in $\A$ such that $F(a)=F(b)$ there is a map
$s:A' \to A$ such that $F(s)$ is an isomorphism in $\B$ and $as=bs$.
\end{enumerate}

\subsection{Remark}
A functor $F: \A \to \B$ between small categories is a localization by
  a calculus of right fractions if and only if
the set $\Sigma$ of maps $f$ in $\A$ such that $F(f)$ is an
isomorphism in $\B$ satisfies a calculus of right fractions
(the dual of \cite[\S I.2.2]{gabrielZisman}) and the induced functor
$\A[\Sigma^{-1}] \to \B$ is an equivalence of categories.

\subsection{Context for theorem \ref{thm:ApprCalcOfFrac2}}
\label{context:thmAppx2}
Consider an additive functor $\A \to \B$ between additive categories.
It induces an exact functor $F:\Ch^b\A \to \Ch^b\B$ between the associated exact
categories of bounded chain complexes
where we call a sequence of chain complexes exact if it is degree-wise split
exact. 
We assume that $F$ is part of an exact form functor
$$(F,\ffi):(\Ch^b\A,w,*,\eta) \to (\Ch^b\B,w,*,\eta)$$
between exact categories 
with weak equivalences
and duality such that the duality compatibility map $F* \to *F$ is a natural
isomorphism.

\subsection{Theorem \rm (Approximation II)}
\label{thm:ApprCalcOfFrac2}
{\it
If in the situation \ref{context:thmAppx2} above, a map in $\Ch^b\A$ is a weak
equivalence iff its image in $\Ch^b\B$ is a weak 
equivalence, and 
if the functor $\A \to \B$ is a localization by a calculus of right
fractions, then $(F,\ffi)$ induces homotopy equivalences
\begin{equation}
\renewcommand\arraystretch{1.5} 
\begin{array}{cccl}
(wS^e_{\bullet}\Ch^b\A)_h &
\stackrel{\sim}{\longrightarrow} &(wS^e_{\bullet}\Ch^b\B)_h & {\rm \it and}\\
GW(\Ch^b\A,w,*,\eta) & \stackrel{\sim}{\longrightarrow} &
GW(\Ch^b\B,w,*,\eta).&
\end{array}
\end{equation}
}
\vspace{2ex}

The proof of theorem \ref{thm:ApprCalcOfFrac2} is a consequence of the
following two lemmas. 

\subsection{Lemma}
\label{lem:extendingCalcFrac}
{\it
Let $F:\A \to \B$ be a localization by a calculus of right fractions.
Then the following holds.
\begin{enumerate}
\item
\label{enum:lem:extendingCalcFrac1}
For every integer $n \geq 0$, the induced functor 
$\Fun([n],\A) \to \Fun([n],\B)$
on diagram categories is a localization by a calculus of right fractions.
\item
\label{enum:lem:extendingCalcFrac2}
If $(F,\ffi):(\A,*,\eta) \to (\B,*,\eta)$ is a form functor between categories
with duality such that the duality compatibility map $\ffi:F* \to *F$ is a
natural isomorphism, then 
the induced functor $\A_h \to \B_h$
on associated categories of symmetric forms is a localization by a
calculus of right fractions. 
\item
\label{enum:lem:extendingCalcFrac3}
If $F$ is an additive functor between additive categories, then the induced
functors $\Ch^b\A \to \Ch^b\B$ and $S_n\A \to S_n\B$ are localizations by a
calculus of right fractions.
\end{enumerate}
}

\subsection*{\it Proof}
The proof is an exercise in clearing denominators, and we omit the
details. 
\Qed

\subsection{Lemma}
\label{lem:calcOfFrac=HtpyEq}
{\it
Let $F:\A \to \B$ be a functor between small categories $\A$ and $\B$ which
is a localization by a calculus of right fractions.
Then $F$ induces a homotopy equivalence on classifying spaces
$$|\A| \stackrel{\sim}{\longrightarrow} |\B|.$$
}

\subsection*{\it Proof}
For a category $\C$, and a subcategory $w\C$, recall that the category 
$\Fun_w([n],\C)$  is the full subcategory of the category $\Fun([n],\C)$ of
functors $[n] \to \C$ which have image in $w\C$.
Maps are natural transformations of functors $[n] \to \C$.
There are
homotopy equivalences of topological realizations of simplicial categories (a
variant of which already appeared in the proof of lemma \ref{lem:simplResLem}) 
\begin{equation}
\label{eqn:lem:calcOfFrac=HtpyEq}
|\C| \stackrel{\sim}{\to} |n \mapsto \Fun_w([n],\C)| \simeq |n \mapsto
w\Fun([n],\C)|
\end{equation}
which is functorial in the pair $(\C,w\C)$.
In the first map, the category $\C$ is considered as a constant simplicial
category, and the functor $\C \to \Fun_w([n],\C)$ sends an object $C \in \C$
to the string consisting of only identity maps on $C$.
This functor is a homotopy equivalence with inverse the functor 
$\Fun_w([n],\C) \to \C$ which sends a string of maps $C_0 \to \cdots
\to C_n$ to $C_0$.
The composition of the two functors is the identity in one direction, and in
the other, it is homotopic to the identity, where the homotopy is given by the
natural transformation from $C_0 \stackrel{1}{\to} C_0 \stackrel{1}{\to} \cdots
C_0$ to $C_0 \to C_1 \to \cdots C_n$ induced by the structure maps of the last
string. 
Therefore, the first map in (\ref{eqn:lem:calcOfFrac=HtpyEq}) is a homotopy
equivalence.
The second map in (\ref{eqn:lem:calcOfFrac=HtpyEq}) is in fact a homeomorphism
as it is the realization in two different orders of the same bisimplicial set.

In order to prove the lemma, let $\sigma\A \subset \A$ be the subcategory of
$\A$ whose maps are the maps which are sent to isomorphisms in $\B$.
By the natural homotopy equivalences in (\ref{eqn:lem:calcOfFrac=HtpyEq}), the
map $|\A| \to |\B|$ in the lemma 
is equivalent to the map 
$|n\mapsto \sigma\Fun([n],\A)| \to |n\mapsto i\Fun([n],\B)|$ 
which is the realization of a map of simplicial categories, so
that it suffices to show that for each integer $n \geq 0$, the functor 
$\sigma\Fun([n],\A) \to i\Fun([n],\B)$ is a homotopy equivalence.
By
lemma \ref{lem:extendingCalcFrac} \ref{enum:lem:extendingCalcFrac1},
this functor is a localization by a calculus of fractions, so that we are
reduced to proving that
$G:\sigma\A \to i\B$ is a homotopy equivalence whenever $\A \to \B$ is a
localization by a calculus of right fractions.
In this case, the comma category $(G\downarrow B)$ is left filtering for every
object $B$ of $\B$, hence contractible.
By theorem A of Quillen, the functor $\sigma\A \to i\B$ is a homotopy
equivalence. 
\Qed

\subsection{\it Proof of theorem \ref{thm:ApprCalcOfFrac2}}
We only prove the first homotopy equivalence in the theorem, the second
homotopy equivalence follows from this and the homotopy
equivalence  
$wS_{\bullet}\Ch^b\A \to wS_{\bullet}\Ch^b\B$ which is proved in the same way
(forgetting forms). 

A map of simplicial categories which is degree-wise a homotopy equivalence
induces a homotopy equivalence after topological realization.
Therefore, it suffices to show that for every $n \geq 0$, the form functor
$(F,\ffi)$ induces a homotopy equivalence
\begin{equation}
\label{eqn:pf:thm:ApprCalcOfFrac2}
(wS^e_{n}\Ch^b\A)_h
\longrightarrow (wS^e_{n}\Ch^b\B)_h.
\end{equation}
By lemma \ref{lem:extendingCalcFrac} \ref{enum:lem:extendingCalcFrac3} 
and in view of the isomorphism $S_n\Ch^b\E =\Ch^bS_n\E$ of exact categories
applied to $\E=\A,\B$, we are reduced to showing that the map 
(\ref{eqn:pf:thm:ApprCalcOfFrac2})
is a homotopy equivalence for $n=0$.
The functor $\Ch^b\A \to \Ch^b\B$ is a localization by a
calculus of right fractions, by \ref{lem:extendingCalcFrac}
\ref{enum:lem:extendingCalcFrac3}.
The assumption that $\Ch^b\A \to
\Ch^b\B$ preserves and detects weak equivalences, implies 
that, on subcategories of weak equivalences, the functor $w\Ch^b\A \to
w\Ch^b\B$ is also a localization by a calculus of right fractions.
By \ref{lem:extendingCalcFrac} \ref{enum:lem:extendingCalcFrac2}, the induced
functor on categories of symmetric forms
$(\Ch^b\A)_h \to (\Ch^b\B)_h$, which is the map
(\ref{eqn:pf:thm:ApprCalcOfFrac2}) in degree $n=0$, 
is a localization by a calculus of right fractions and therefore a homotopy
equivalence, by lemma \ref{lem:calcOfFrac=HtpyEq}.
\Qed

\section{From exact categories to chain complexes}

The purpose of this section is to prove proposition \ref{prop:gilletWald}
which allows us the replace the Grothendieck-Witt space of an exact category
with duality by the Grothendieck-Witt space of the associated category of
bounded chain complexes.

\subsection{Chain complexes and dualities}
\label{subsec:chainCx}
Let $(\E,*,\eta)$ be an exact category with duality, 
that is, an exact category with weak equivalences and duality where all weak
equivalences are isomorphisms.
Let $\Ch^b(\E)$ be the category of bounded chain complexes 
$$(E,d):\hspace{4ex} \cdots \to E_{n-1}
\stackrel{d_{n-1}}{\to} E_n \stackrel{d_{n}}{\to} E_{n+1} \to \cdots
\hspace{4ex}, d_nd_{n-1}=0,$$
in $\E$.
A sequence of chain complexes 
$(E',d) \to (E,d) \to
(E'',d)$ is {\em exact} if it is degree-wise exact in
$\E$, that is, if the sequence $E'_n \rightarrowtail E_n \twoheadrightarrow
E''_n$ is exact for all $n$.
Call a chain complex $(E,d)$ in $\E$ {\em strictly
  acyclic} if every 
differential $d_n$ is the composition $E_n \twoheadrightarrow \im d_n
\rightarrowtail E_{n+1}$  of a deflation followed by an inflation, and the
sequences $\im d_{n-1} \rightarrowtail E_n \twoheadrightarrow \im d_n$  are
exact in $\E$.
A chain complex is called {\em acyclic} if it is homotopy equivalent to a
strictly acyclic chain complex.
A map of chain complexes is a {\em quasi-isomorphism}  if its cone is acyclic.
Write $\quis$ for the set of quasi-isomorphisms, then the triple
$$(\Ch^b(\E),\quis)$$
is an exact category with weak equivalences.

For $n \in \Z$, the duality $(*,\eta)$ induces a (naive) duality
$(*^n,\eta^n)$ on $\Ch^b\E$ which 
on objects $(E,d)$ and on chain maps $f:(E,d) \to (E',d)$ is given by the
formulas
\renewcommand\arraystretch{1.5} 
$$\begin{array}{lcllcl}
(E^{*^n})_i &=&(E_{-i-n})^*,&
 (f^{*^n})_i&=&(f_{-i-n})^*,\\
 (d^{*^n})_i&=& (d_{-i-1-n})^*,&
 (\eta^n_{E})_i&=& (-1)^{\frac{n(n-1)}{2}}\eta_{E_i}.
\end{array}
$$
With these definitions, we have an exact category with weak equivalences and
duality
$$(\Ch^b(\E),\quis,*^n,\eta^n).$$
If $n=0$ we may simply write $(*,\eta)$ for $(*^0,\eta^0)$.

\subsection{Remark}
\label{rmk:Chbmod4}
The functor $T:\Ch^b\E \to \Ch^b\E$ given by the formula
$$(TE)_i=E_{i+1},\hspace{2ex} T(f)_i=f_{i+1},\hspace{2ex}(d_{TE})_i=d_{i+1}$$
defines a duality preserving isomorphism of exact categories with duality
$$T: (\Ch^b\E,*^n,\eta^n) \cong (\Ch^b\E,*^{n+2},-\eta^{n+2}).$$

\subsection{Remark}
\label{rmk:ChoiceOfSign}
There is another (more natural) sign choice for defining induced dualities
$(\sharp^n,\can^n)$ on $\Ch^b\E$ coming from the internal hom of chain
complexes, compare \ref{subsec:DGmodDual}.
They are given by the formulas
\renewcommand\arraystretch{1.5} 
$$\begin{array}{lcllcl}
(E^{\sharp^n})_i &=&(E_{-i-n})^*,&
 (f^{\sharp^n})_i&=&(f_{-i-n})^*,\\
 (d^{\sharp^n})_i&=& (-1)^{i+1}(d_{-i-1-n})^*,&
 (\can^n_{E})_i&=& (-1)^{i(i+n)}\eta_{E_i}.
\end{array}
$$
The identity functor on $\Ch^b\E$ together with the duality compatibility
isomorphism 
$\eps^n_E:E^{*^n} \to E^{\sharp^n}$ which in degree $i$ is 
$$(\eps^n_E)_i=(-1)^{\frac{i(i+1)}{2}}\phantom{a} id_{E^*_{-i-n}}: (E^{*^n})_i \to
(E^{\sharp^n})_i$$ 
defines an isomorphism of exact categories with duality
$$(id,\eps^n):(\Ch^b\E,*^n,\eta^n) \stackrel{\cong}{\longrightarrow}
(\Ch^b\E,\sharp^n,\can^n).$$
For the purpose of proving proposition \ref{prop:gilletWald} below, the duality
$(*,\eta)$ is convenient.
For most other purposes, the duality $(\sharp,\can)$ is more natural.
It is the latter duality, which we
will use from \S \ref{sec:HigherGWschemes} on.
In any case, both give rise to isomorphic exact categories with duality and
thus have isomorphic Grothendieck-Witt spaces.

\subsection{Exercise}
Show that the exact categories with weak equivalences and duality 
$(\Ch^b\E,*^n,\eta^n)$ and $(\Ch^b\E,\sharp^n,\can^n)$ have symmetric cones in
the sense of definition \ref{dfn:symCone} (hint: see \ref{subsec:canSymCone}).
\vspace{2ex}

For an exact category with duality $(\E,*,\eta)$, 
inclusion as complexes concentrated in degree $0$, defines a duality
preserving exact functor
\begin{equation}
\label{eqn:GilletWald}
(\E,i,*,\eta) \to (\Ch^b\E,\quis,*,\eta).
\end{equation}
The following proposition generalizes
\cite[Theorem 1.11.7]{TT}, see also remark \ref{rmk:exStrOnChb}.

\subsection{Proposition}
\label{prop:gilletWald}
{\it 
For an exact category with duality $(\E,*,\eta)$, 
the functor (\ref{eqn:GilletWald})
induces a homotopy equivalence of Grothendieck-Witt spaces
$$GW(\E,*,\eta) \stackrel{\sim}{\longrightarrow}
GW(\Ch^b\E,\quis,*,\eta).$$ 
}
\vspace{2ex}

We will reduce the proof of the proposition to ``semi-idempotent complete''
exact categories. 
This has the advantage that for such categories, every acyclic complex is
strictly acyclic.
Here are the relevant definitions and facts.

\subsection{Semi-idempotent completions}
Call an exact category $\E$ {\em semi-idempotent complete} if any map $p: A
\to B$ which has a section $s:B \to A$, $pi=1$, is a deflation in $\E$.
A semi-idempotent complete exact category has the following property:
any map $B \to C$ for which there is a map $A \to B$ such that the
composition $A \to C$ is a deflation, is itself a deflation.
This is because a semi-idempotent complete exact category satisfies Thomason's
axiom \cite[A.5.1]{TT} so that the standard embedding of $\E$ into the
category of left exact functors $\E^{op} \to \langle ab \rangle$ into the
category of abelian groups is closed
under kernels of surjections \cite[Theorem A.7.1 and Proposition A.7.16
(b)]{TT}. 
For a semi-idempotent complete exact category, every acyclic complex is
strictly acyclic.

The semi-idempotent completion of an exact category $\E$ is the full
subcategory $\widetilde{\E}^{semi} \subset \widetilde{\E}$ of the idempotent
completion $\widetilde{\E}$ of $\E$ of those objects which are stably in $\E$.
Clearly, $\widetilde{\E}^{semi}$ is semi-idempotent complete, and the map 
$K_0(\E) \to K_0(\widetilde{\E}^{semi})$ is an isomorphism.
If $(\E,*,\eta)$ is an exact category with duality, then 
$(\widetilde{\E}^{semi},*,\eta)$ is an exact category with duality such that
the fully exact inclusion $\E \subset \widetilde{\E}^{semi}$ is duality
preserving.
If $(\A,w,*,\eta)$ is an exact category with weak equivalences and duality,
then  $(\widetilde{\A}^{semi},w,*,\eta)$ inherits the structure of an exact
category 
with weak equivalence and duality from $(\widetilde{\A},w,*,\eta)$, see
\ref{subsec:idempCompl}, so that the inclusion $\A \subset
\widetilde{\A}^{semi}$ is duality preserving.

\subsection{Lemma}
\label{lem:strongCof}
{\it
Let $(\A,w,*,\eta)$ be an exact category with weak equivalences and strong
duality, then the inclusion $\A \subset \widetilde{\A}^{semi}$ induces a
homotopy equivalence of Grothendieck-Witt spaces
$$GW(\A,w,*,\eta) \to GW(\widetilde{\A}^{semi},w,*,\eta).$$
}

\subsection*{\it Proof}
For an exact category with duality $(\E,*,\eta)$,
the map $GW_i(\E) \to GW_i(\widetilde{\E}^{semi})$ is an isomorphism for $i
\geq 1$ and it is injective for $i=0$, by cofinality proved in
\cite{myGWex}.
The map $GW_0(\E) \to GW_0(\widetilde{\E}^{semi})$
is also surjective, hence an isomorphism, since for every symmetric space
$(X,\ffi)$ 
in $\widetilde{\E}^{semi}$, we can find an $A$ in $\E$ such that $X \oplus A$
is in $\E$, and thus, $[X,\ffi]=[(X,\ffi)\perp \H A]-[\H A]$ is in the image of
the map.
Therefore, lemma \ref{lem:strongCof} holds when $w$ is
the set of isomorphisms.
The lemma now follows from this case and the simplicial resolution lemma
\ref{lem:simplResLem} since $\Fun_w(\underline{n},\widetilde{\A}^{semi})$ is
the semi-idempotent completion of $\Fun_w(\underline{n},\A)$ (for a string
$X_{n'} \to \cdots \to X_n$ of weak equivalences in $\widetilde{\A}^{semi}$,
there is an object $A$ of $\A$ such that $X_i \oplus A$ is in $\A$ for all $i
\in \underline{n}$, so that $(X_{n'} \to \cdots \to X_n) \oplus (A
\stackrel{1}{\to}  \cdots \stackrel{1}{\to} A)$ is a string of weak
equivalences in $\A$).
\Qed

\subsection*{\it Proof of proposition \ref{prop:gilletWald}}
In view of lemma \ref{lem:strongCof}, we can assume
$\E$ to be semi-idempotent complete, so that acyclic complexes are strictly
acyclic.
Moreover, by the strictification lemma \ref{lem:strictification}, we can
assume $\E$ to have a strict duality, since an exact category with duality is
equivalent to its strictification. 
Let $\Ac^b\E\subset \Ch^b\E$ be the full subcategory of acyclic chain
complexes.
It inherits the structure of an exact category with weak equivalences and
strict duality from $\Ch^b\E$.
Consider the commutative diagram of exact categories with weak equivalences and
strict duality 
$$\xymatrix{
    0 \ar[r]\ar[d] &  (\Ac^b\E,i) \ar[d] \ar[r]& (\Ac^b\E,\quis) \ar[d]\\
   (\E,i) \ar[r]  & (\Ch^b\E,i) \ar[r] & (\Ch^b\E,\quis).}$$
We will show that 
\begin{enumerate}
\item
\label{pf:gilletWald1}
the left square induces a homotopy cartesian square of
Grothendieck-Witt spaces, and that 
\item
\label{pf:gilletWald2}
the map $GW_0(\E,i) \to GW_0(\Ch^b\E,\quis)$ is surjective.
\end{enumerate}
By proposition \ref{prop:iCOW}, the right hand square induces a homotopy
cartesian square of Grothendieck-Witt spaces, so that 
by \ref{pf:gilletWald1} the same is true for the outer square.
Together with \ref{pf:gilletWald2}, this implies the proposition.

We prove \ref{pf:gilletWald1}.
For $n \geq 0$, let $Ch^b_{[-n,n]}\E \subset Ch^b\E$ and $Ac^b_{[-n,n]}\E
\subset Ac^b\E$ be the 
  full subcategories of those chain complexes which are concentrated in
  degrees $[-n,n]$. 
They inherit a structure of exact categories with duality.
Note that the inclusion $Ac^b_{[-n,n]}\E \subset Ch^b_{[-n,n]}\E$ is $0
\subset \E$ for $n=0$.
The natural inclusions induce a commutative diagram of exact categories with
duality
\begin{equation}
\label{pf:diag1}
\xymatrix{
0 \ar[r] \ar[d] & Ac^b_{[-n,n]}{\E} \ar[r] \ar[d] & Ac^b_{[-n-1,n+1]}{\E}
\ar[d] \\
{\E} \ar[r] &  Ch^b_{[-n,n]}{\E}  \ar[r] & Ch^b_{[-n-1,n+1]}{\E}.}
\end{equation}
We will show that the right-hand square induces a homotopy cartesian square of
Grothendieck-Witt spaces.
Then, by induction, the outer square induces a homotopy cartesian
square, too.
Taking the colimit over $n$ of the Grothendieck-Witt spaces of the outer
squares yields the desired homotopy cartesian square, since the
Grothendieck-Witt space functor $GW$ commutes with filtered colimits.

Consider the following form functors between exact categories with duality
(equipped with the obvious duality compatibility maps):

\hspace{-4ex}
\renewcommand\arraystretch{2} 
\begin{tabular}{lcl}
$\H\E \times Ac^b_{[-n,n]}{\E} $&$\stackrel{s}{\to}$&$ Ac^b_{[-n-1,n+1]}{\E}:$\\
\phantom{12}$(A,B),C$  &$\mapsto$&$ A \stackrel{(1\ 0)}{\longrightarrow} A \oplus C_{-n} \to C_{-n+1} \to \cdots \to C_{n-1} \to C_n \oplus B^*  \stackrel{(0\  1)}{\longrightarrow} B^*$\\
$Ch^b_{[-n-1,n+1]}{\E} $&$\stackrel{\rho}{\to}$&$ \H\E \times Ch^b_{[-n,n]}{\E}:$\\
\phantom{1223}$C$&$ \mapsto$&$ (C_{-n-1},(C_{n+1})^*), C_{-n} \to
  \cdots \to C_n$\\
$\H\E \times Ac^b_{[-n,n]}{\E} $&$ \to $&$\H\E \times Ch^b_{[-n,n]}{\E}:$\\
\phantom{12}$(A,B),C$&$ \mapsto $&$(A,B),\ \   A \oplus C_{-n} \to C_{-n+1} \to \cdots \to
C_{n-1} \to C_n   \oplus B^*.$
\end{tabular}
These functors, together with the natural inclusions, fit into a commutative diagram of exact
  categories with duality 
$$
\xymatrix{
Ac^b_{[-n,n]}{\E} \ar[rrr] \ar[dr]^{\left(\begin{smallmatrix}0\\1\end{smallmatrix}\right)} \ar[dd]&& & Ch^b_{[-n,n]}{\E} \ar[ld]_{\left(\begin{smallmatrix}0\\1\end{smallmatrix}\right)} \ar[dd]\\
& {\H}\E \times Ac^b_{[-n,n]}{\E} \ar[r] \ar[dl]^{s} &{\H}\E \times Ch^b_{[-n,n]}{\E}&\\
Ac^b_{[-n-1,n+1]}{\E} \ar[rrr] &&& Ch^b_{[-n-1,n+1]}{\E}. \ar[lu]^{\rho}
}
$$
The upper square induces a homotopy cartesian square of Grothendieck-Witt
spaces.
By additivity (theorem \ref{thm:addtyFun} or \cite[7.1]{myGWex}), the diagonal maps
$s$ and $\rho$ in the lower square induce 
homotopy equivalences of Grothendieck-Witt spaces (details below).
Therefore, the outer square induces a  homotopy cartesian square as well.
To see that the functors $s$ and $\rho$ induce homotopy equivalences, consider
the following functors of categories with duality 

\hspace{-4ex}
\renewcommand\arraystretch{2} 
\begin{tabular}{lcl}
$ Ac^b_{[-n-1,n+1]}{\E}$&$\stackrel{r}{\to}$&$\H\E \times Ac^b_{[-n,n]}{\E}: C
\mapsto (C_{-n-1},(C_{n+1})^*),$\\
&&$  C_{-n}/C_{-n-1} \to
  C_{-n+1} \to \cdots \to C_{n-1} \to {\ker}(C_n \to C_{n+1})$\\
$ \H\E \times
Ch^b_{[-n,n]}{\E}$&$\stackrel{\sigma}{\to}$&$Ch^b_{[-n-1,n+1]}{\E}: (A,B),C
\mapsto $\\
&&$A \stackrel{0}{\to} C_{-n} \to
  C_{-n+1} \to \cdots \to C_{n-1} \to C_n \stackrel{0}{\to} B^*.$
\end{tabular}
We have  $\rho\sigma=id$.
The identity functor $id$ on ${Ch^b_{[-n-1,n+1]}}$ has a totally isotropic
subfunctor $G\subset id$ given by 
$$\xymatrix{
G \xymono[d] & 0 \ar[r] \ar[d] & \cdots \ar[r] & 0 \ar[r] \ar[d] & C_{n+1}
\ar[d]^1\\
id & C_{-n-1} \ar[r] & \cdots \ar[r] & C_n \ar[r] & C_{n+1}.
}
$$
By additivity \cite[7.1]{myGWex}, the identity functor $id$ and $\sigma\rho=
G^{\perp}/G \oplus {\H}G$ 
induce homotopic maps on Grothendieck-Witt spaces, thus $\rho$ induces a
homotopy equivalence.
Similarly, we have $rs=id$,  and the identity functor $id$ on
$Ac^b_{[-n-1,n+1]}{\E}$ has a totally 
isotropic subfunctor $F \subset id$ given by
$$\xymatrix{
F \xymono[d] & C_{-n-1} \ar[r]^1 \ar[d]^1 & C_{-n-1} \xymono[d] \ar[r] & 0
\ar[d]\ar[r] & \cdots \ar[r] & 0 \ar[d]\\
id & C_{-n-1} \xymono[r]  & C_{-n} \ar[r] & C_{-n+1} \ar[r] & \cdots \xyepi[r] &
C_{n+1}. 
}
$$
Again, by additivity \cite{myGWex}, the identity functor $id$ on
$Ac^b_{[-n-1,n+1]}{\E}$ and $sr= F^{\perp}/F \oplus {\H}F$
induce homotopic maps on Grothendieck-Witt spaces.
It follows that $s$ induces a homotopy equivalence.
This finishes the prove of \ref{pf:gilletWald1}.

We are left with proving \ref{pf:gilletWald2}.
We will show that 
\begin{itemize}
\item[(c)]
a symmetric space
$(A,\alpha)$, where $A$ is supported in $[-n,n]$, $n\geq 1$, equals 
$[A,\alpha]=[B,\beta]+[\H(C)]$
in the Grothendieck-Witt group $GW_0(\Ch^b\E,\quis)$,
where $B$ is supported in $[-n+1,n-1]$.
\end{itemize}
By induction, $[A,\alpha]$ is then a sum of hyperbolic objects plus a
symmetric space supported in degree $0$.
Since the latter two kinds of symmetric spaces are obviously in the image of
$GW_0(\E,i) \to GW_0(\Ch^b\E,\quis)$, this proves $(b)$.
To show (c), let $n\geq 1$ and let $(A,\alpha)$ be a symmetric space supported
in $[-n,n]$.
Since the cone of $\alpha$ is acyclic (hence strictly acyclic, by
semi-idempotent completeness of $\E$), the map 
$\left(\begin{smallmatrix}d_{-n}\\
    \alpha_{-n}\end{smallmatrix}\right):A_{-n} \to A_{-n+1}\oplus A_n^*$ 
is an inflation.
Define a complex $\tilde{A}$, also supported in $[-n,n]$, by
$$A_{-n} \stackrel{\left(\begin{smallmatrix}d_{-n}\\
    \alpha_{-n}\end{smallmatrix}\right)}{\rightarrowtail} 
A_{-n+1} \oplus A_n^* \stackrel{\left(\begin{smallmatrix}d &
    0 \end{smallmatrix}\right)}{\longrightarrow} A_{-n+2}
\stackrel{d}{\to} \cdots
\stackrel{d}{\to}  A_{n-2} 
\stackrel{\left(\begin{smallmatrix}d\\
    0\end{smallmatrix}\right)}{\longrightarrow} A_{n-1} \oplus A_n 
\stackrel{\left(\begin{smallmatrix}d_{n-1}&
    1\end{smallmatrix}\right)}{\twoheadrightarrow} A_n.
$$
The complex $\tilde{A}$ is equipped with a non-singular symmetric form
$\tilde{\alpha}$ which is $\alpha_i$ in degree $i$ except in degrees $i=-n+1,
n-1$ where it is $\alpha_i \oplus 1$.
We have a symmetric space in the category of admissible short complexes in
$\Ch^b\E$ 
$$A^*_n[n-1] \rightarrowtail  \tilde{A} \twoheadrightarrow  A_n[-n+1]$$
with form $(1,\tilde{\alpha},\eta)$, where for an object $E$ of $\E$, we
denote by $E[i]$ the complex which is $E$ in degree $-i$ and $0$ elsewhere.
The maps $A^*_n \rightarrowtail  \tilde{A}_{-n+1} = A_{-n+1} \oplus A_n^*$ and 
$A_{n-1} \oplus A_n = \tilde{A}_{n-1} \to A_n$ are the canonical inclusions
and projections, respectively.
Since $(A,\alpha)$ is the zero-th homology of this admissible short
complex equipped with its form, we have $[\tilde{A},\tilde{\alpha}]=
[A,\alpha]+[\H (A_n[-n+1])]$ in $GW_0(\Ch^b\E,\quis)$.
There is another symmetric space in 
the category of admissible short complexes in $\Ch^b\E$ 
$$CA_{-n}[n-1] \rightarrowtail  \tilde{A} \twoheadrightarrow  CA_n[n]$$
with non-singular from $(\alpha_{-n},\tilde{\alpha},\alpha_n)$, where
for an object $E$ of $\E$, we write $CE[i]$ for the complex $E
\stackrel{1}{\to} E$ placed in degrees $i$ and $i-1$.
The maps $CA_{-n}[n-1] \rightarrowtail  \tilde{A}$ and $\tilde{A}
\twoheadrightarrow  CA_n[n]$ 
are the unique maps which are the identity in degree $-n$ and $n$,
respectively.
Since $CA_{-n}[n-1]$ and $CA_n[n]$ are acyclic, the form on the admissible short
complex is non-singular and its zero-th homology, which is concentrated
in degrees $[-n+1,n-1]$, has the same class in $GW_0(\Ch^b\E,\quis)$ as
$(\tilde{A},\tilde{\alpha})$.
\Qed

\subsection{Remark}
\label{rmk:exStrOnChb}
We can equip $(\Ch^b\E,\quis,*,\eta)$ with two exact structures, the
one defined in 
\ref{subsec:chainCx}, and the degree-wise split exact structure.
By lemma \ref{lem:ChExStr}, the two yield homotopy equivalent
Grothendieck-Witt spaces.

\section{DG-Algebras on ringed spaces and dualities}
\label{sec:DGalg}

In this section we recall basic definitions and facts about differential
graded algebras and modules over them.
Besides fixing terminology, the main point here is the construction of the
canonical symmetric cone in \ref{subsec:canSymCone}, and the interpretation of
certain form functors as symmetric forms in dg bimodule categories, see
\ref{subsec:FormFunAsTensProd}. 

\subsection{DG $\kappa$-modules}
\label{subsec:DGk}
Let $\kappa$ be a commutative ring.
Unless otherwise indicated, modules will always mean left module, tensor
product $\otimes$ will be tensor product $\otimes_{\kappa}$ over $\kappa$, and
homomorphism sets $Hom(M,N)$ between $\kappa$-modules means set of $\kappa$-linear
homomorphisms, and is itself a $\kappa$-module.
Recall that a differential graded $\kappa$-module $M$ is a graded $\kappa$-module
$\bigoplus_{n\in \Z}M_n$ together with a $\kappa$-linear map $d:M_n \to
M_{n+1}$, $n \in \Z$, called differential of $M$, satisfying $d\circ d=0$.
In other words, $M$ is a chain complex of $\kappa$-modules.
A map of dg $\kappa$-modules is a map of graded $\kappa$-modules commuting with the
differentials. 
For two dg $\kappa$-modules $M$, $N$, the tensor product dg $\kappa$-module $M\otimes N$
and the homomorphism dg $\kappa$-module $[M,N]$ are defined
by the usual formulas
$$
\renewcommand\arraystretch{1.5} 
\begin{array}{ll}
(M\otimes N)_n = \bigoplus_{i+j}M_i\otimes N_j & d(x\otimes
y)=dx\otimes y + (-1)^{|x|}x\otimes dy\\
{[M,N]}_n
= \prod_{j-i=n} Hom(M_i,N_j) &  d(f) = d_N \circ f - (-1)^{|f|}f \circ
d_M.
\end{array}$$

\subsection{DGAs and modules over them}
A differential graded $\kappa$-algebra (dga) is a
dg $\kappa$-module $A$ equipped with dg $\kappa$-module maps $\cdot: A \otimes A \to A$
and $\kappa \to A$, called multiplication and unit, making the usual
associativity 
and unit diagrams commute \cite[diagrams (1), (2), p. 166]{MacLane}. 
In other words, $A$ is an associative graded $\kappa$-algebra with
multiplication 
satisfying $d(a\cdot b)=(da)\cdot b+(-1)^{|a|}a\cdot (db)$.
For dg algebras $A$ and $B$, we denote by 
$A\LRmod B$ the category
of left $A$ and right $B$-modules. 
Its objects are the dg $\kappa$-modules $M$ 
equipped with dg $\kappa$-maps $A\otimes M \to M$ and $M \otimes B
\to M$, called multiplication, both of which are associative and unital, and
furthermore, $(am)b=a(mb)$ for all $a\in A$, $m\in M$ and $b\in B$.
We also denote by $A\Lmod = A\LRmod \kappa$, $\Rmod A=
\kappa\LRmod A$, $A\Bimod = A\LRmod A$ 
the categories of dg left
$A$-modules, right $A$-modules, and of dg $A$-bimodules.

Let $A$, $B$, $C$ be dg algebras.
Recall that for a right $B$ module $M$ and left $B$-module $N$,
tensor product $M\otimes_BN$ over of $M$ and $N$ over $B$
is the dg $\kappa$-module which is the co-equalizer
$$\xymatrix{
M\otimes B \otimes N \ar@<.7ex>[r]^{\hspace{2ex}1\otimes \mu} 
\ar@<-.7ex>[r]_{\hspace{2ex}\mu \otimes 1}  & M\otimes N \ar[r] & M \otimes_B
N
}
$$
in the category of dg $\kappa$-modules, where $\mu$ stands for the
multiplications $M\otimes B \to M$ and $B \otimes N \to N$.
For two dg right $C$-modules $M$ and $N$, the dg $\kappa$-module 
$[M,N]_C$ of right $C$-module morphisms is the equalizer in the category of dg $\kappa$-modules
$$\xymatrix{
[M,N]_C \ar[r] & [M,N] \ar@<.7ex>[rr]^{\hspace{-2ex}[\mu,1]} 
\ar@<-.7ex>[rr]_{\hspace{-2ex}[1,\mu]\circ (?\otimes 1_C)} && [M\otimes C,N],
}$$
where $(?\otimes 1_C): [M,N] \to [M\otimes C, N \otimes C]$ is the dg
$\kappa$-module map $f \mapsto f\otimes 1_C$ defined by 
$(f\otimes 1_C)(x\otimes c) = f(x)\otimes c$ for $x \in M$ and $c \in C$.  
Similarly, one can define the dg $\kappa$-module of left $C$-module morphisms
$_C[M,N]$ for two dg left $C$-modules.

Tensor product $\otimes_B$ and right $C$-module morphisms
$[\phantom{M},\phantom{N}]_C$  define functors
$$
\renewcommand\arraystretch{1.5} 
\begin{array}{ccc}
\otimes_B : &  A\LRmod B \times B\LRmod C \longrightarrow A\LRmod C: &
M,N \mapsto M\otimes_BN\\
{[\phantom{A},\phantom{B}]}_C: 
& (B\LRmod C)^{op} \times A\LRmod C
\longrightarrow A\LRmod B:& M, N \mapsto [M,N]_C.
\end{array}$$

\subsection{DGAs with involution}
For a dg $\kappa$-module $M$, we denote by $M^{op}$ the opposite dg $\kappa$-module
which, as a dg $\kappa$-module, is simply $M$ itself.
To avoid confusion we may sometimes write $x^{op}, y^{op}$, etc, for 
the elements in $M^{op}$ corresponding to $x,y \in M$, so that
$d(x^{op})=(dx)^{op}$, for instance, denotes an equation in $M^{op}$ as
opposed to in $M$. 
For a dg $\kappa$-algebra $A$, its opposite dga $A^{op}$ has underlying dg
$\kappa$-module the opposite module of $A$ and multiplication 
$x^{op}y^{op}=(-1)^{|x||y|}(yx)^{op}$.
A {\em dg algebra with involution} is a dg $\kappa$-algebra $A$ together with an
isomorphism  
$A \to A^{op}: a \mapsto \bar{a}$ of dgas satisfying $\bar{\bar{a}}=a$ for
all $a\in A$.
The base commutative ring $\kappa$ is always considered as a dga with trivial
involution $\kappa \to \kappa^{op}:x \mapsto x$.

Let $A,B$ be dgas with involution.
For a left $A$, right $B$-module $M \in A\LRmod B$ its opposite module
$M^{op}$ is a left $B^{op}$ and right $A^{op}$-module which we consider as a
left $B$ and right $A$-module via the isomorphisms $A \to A^{op}$, $B \to
B^{op}$, that 
is, $m^{op}\cdot a = (-1)^{|a||m|}(\bar{a}\cdot m)^{op}$, $b\cdot m^{op} =
(-1)^{|b||m|}(m\cdot \bar{b})^{op}$ for $a\in A$,
$b\in B$, $m \in M$.
For dgas with involution $A$, $B$, $C$ and $M\in A\LRmod
B$, $N\in B\LRmod C$, the commutativity
isomorphism of the 
tensor product defines an $A\LRmod C$-isomorphism
$$c:M\otimes_B N \stackrel{\cong}{\longrightarrow} (N^{op} \otimes_B
M^{op})^{op}:\phantom{a} x\otimes y \mapsto (-1)^{|x||y|}(y^{op}\otimes
x^{op})^{op}.$$ 
This can be iterated to obtain for rings with involution $A$, $B$, $C$, $D$
and $M\in A\LRmod B$, $N\in B\LRmod C$, $P \in C\LRmod D$ 
an isomorphism $c_3$ in $A\LRmod D$ defined by 
$$
\renewcommand\arraystretch{1.5} 
\begin{array}{ccl}
(M\otimes_BN\otimes_CP) & \stackrel{c_3}{\longrightarrow} &
(P^{op}\otimes_CN^{op}\otimes_BM^{op})^{op}\\
 x\otimes y \otimes z & \mapsto &
(-1)^{|x||y|+|x||z|+|y||z|} (z^{op}\otimes y^{op}\otimes x^{op})^{op}.
\end{array}
$$

\subsection{DG-modules and dualities}
\label{subsec:DGmodDual}
Let $A$ be a dga with involution, and 
let $I$ be an $A$-bimodule equipped with an $A$ bimodule isomorphism 
$i:I \to I^{op}$ such that $i^{op}\circ i=id$, for instance $A$ itself with
$i(x)=\bar{x}$.
We call the pair $(I,i)$ a {\em duality coefficient} for the category $A\Lmod$
of dg $A$-modules, as it
defines a duality $\sharp_{(I,i)}:(A\Lmod)^{op} \to A\Lmod$ by
$$M^{\sharp_{(I,i)}}=[M^{op},I]_A.$$
The canonical double dual identification $\can_{(I,i),M}:M \to
M^{\sharp_{(I,i)}\sharp_{(I,i)}}$ 
is the left $A$-module map given by the formula 
$$\can_{(I,i),M}(x)(f^{op})=(-1)^{|x||f|}i(f(x^{op}))$$
for $f\in [M^{op},I]_A$ and $x\in M$.
It is a straight forward to check the identity
$\can_{(I,i),M}^{\sharp}\can_{(I,i),M^{\sharp}} = 1_{M^{\sharp}}$, so that the triple
$(A\Lmod,\sharp_{(I,i)}, \can_{(I,i)})$ 
is a category with duality.
In this paper, for the duality $\sharp_{(I,i)}$, the double dual identification
  will {\em always} be the natural map $\can_{(I,i)}$, 
so that we will write 
$$(A\Lmod,\sharp_{(I,i)})$$ for the category with duality
$(A\Lmod,\sharp_{(I,i)}, \can{(I,i)})$, the double dual identification being 
$\can_{(I,i)}$.
If $i:I \to I^{op}$ is understood, we may write $\sharp_I$ instead of
$\sharp_{(I,i)}$. 

To give a symmetric form $\ffi:M \to [M^{op},I]_A$ in
$(A\Lmod,\sharp_I)$ is the same as to give an
$A$-bimodule map $\hat{\ffi}:M \otimes M^{op} \to I$ such that the diagram
$$\xymatrix{M \otimes M^{op} \ar[r]^{\hspace{2ex}\hat{\ffi}} \ar[d]_c & I
  \ar[d]^i\\ 
            (M \otimes M^{op})^{op} \ar[r]^{\hspace{6ex}\hat{\ffi}^{op}} &
            I^{op} }$$
commutes.
The bijection is given by the identity $\hat{\ffi}(x\otimes
y^{op})=\ffi(x)(y^{op})$.

\subsection{The canonical symmetric cone}
\label{subsec:canSymCone}
Let $C=\kappa\cdot 1_C \oplus \kappa\cdot \epsilon$ be the dg
$\kappa$-module whose underlying $\kappa$-module is free of rank 
$2$ with basis $1_C$ and $\epsilon$ in degrees $0$ and $-1$, respectively, and
differential $d\epsilon= 1_C$.
In fact, $C$ is a commutative dga with unit $1_C$ and unique multiplication.
Let $A$ be a dga and  $M$ a (left) dg $A$-module.
We write $C:A\Lmod \to A\Lmod$ for the
functor $M \mapsto M\otimes C$, and 
$i_M$ for the natural inclusion $M \to CM=M\otimes C: x \mapsto x\otimes 1_C$.
Similarly, we write $P:A\Lmod \to A\Lmod$ for the functor
$M\otimes [C^{op},\kappa]$ and $p_M$ for the natural surjection
$PM = M\otimes [C^{op},\kappa] \to M:m\otimes g \mapsto m\cdot g(1_C^{op})$.
If $A$ is a dga with involution, and $(I,i)$ a duality coefficient for
$A\Lmod$, we define a natural transformation 
$$\gamma_M:[M^{op},I]_A\otimes [C^{op},\kappa] \longrightarrow
[(M\otimes C)^{op},I]_A$$
by the formula $\gamma_M(f\otimes g)((x\otimes
a)^{op})=(-1)^{|a||x|}f(x^{op})\cdot g(a^{op})$.
One checks the equality 
$i_M^{\sharp_I}\circ \gamma_M = p_{M^{\sharp_I}}$.

Therefore, an exact category with weak equivalences and duality
$(\E,w,*,\eta)$ 
which admits a fully faithful and duality preserving embedding into 
$(A\Lmod,\sharp_I)$ has a symmetric cone in the sense
of \ref{dfn:symCone} provided
the functors $C$ and $P$ restrict to exact endofunctors of $\E$, 
the natural maps $i_M$ and $p_M$ are inflation and deflation, and the objects
$CM$ and $PM$ are $w$-acyclic for all $M \in
\E$.

\subsection{Symmetric forms in bimodule categories and their tensor product}
\label{subsec:symmAmodB}
Let $A$ and $B$ be dgas with involution, and let $(I,i)$, $(J,j)$ be duality
coefficients for $A\Lmod$ and $B\Lmod$, respectively.
A {\em symmetric form in $A\LRmod B$, with respect to the duality coefficients
$(I,i)$ and $(J,j)$}, is a pair $(M,\ffi)$ where $M\in A\LRmod B$ is 
a left $A$ and right $B$-module, and 
$\ffi:M \otimes_B J \otimes_B M^{op} \to I$ is an $A$-bimodule map
making the diagram of $A$-bimodule maps
\begin{equation}
\label{eqn:AModBsymmForm}
\xymatrix{M \otimes_B J \otimes_B M^{op} \ar[r]^{\hspace{6ex}\ffi} \ar[d]_{c_3\circ
    (1\otimes j \otimes 1)} &  I \ar[d]^i \\
(M \otimes_B J \otimes_B M^{op})^{op} \ar[r]_{\hspace{8ex}\ffi^{op}} & I^{op}}
\end{equation}
commute.
Isometries and orthogonal sums of symmetric forms in $A\LRmod B$ are defined
in the obvious way.

Tensor product of symmetric forms is defined as follows.
Let $A$, $B$, $C$ be dg algebras with involution,
and let 
$(I,i)$, $(J,j)$, $(K,k)$ be duality coefficients for $A\Lmod$, $B\Lmod$ and
$C\Lmod$, respectively.
Further, let $M\in A\LRmod B$ and $N \in B\LRmod C$ be equipped with 
symmetric forms given by the
$A$- and $B$-bimodule maps $\ffi:M \otimes_B J \otimes_B M^{op} \to I$ and
$\psi:N \otimes_C K \otimes_C N^{op} \to J$ (making diagram
(\ref{eqn:AModBsymmForm}) and its analog for $\psi$ commute).
The tensor product 
$$(M,\ffi)\otimes_B(N,\psi)$$
of the symmetric forms
$(M,\ffi)$ and $(N,\psi)$ has
$M\otimes_BN$ as underlying 
left $A$ and right $C$-module, and is equipped with the symmetric form which
is the $A$-bimodule map
$$
\renewcommand\arraystretch{1.5} 
\hspace{-1ex}
\begin{array}{ccccc}
(M\otimes_BN)\otimes_CK\otimes_C(M\otimes_BN)^{op} 
&\stackrel{c}{\to}&
M\otimes_BN\otimes_CK\otimes_CN^{op}\otimes_BM^{op}&&\\
&\stackrel{\psi}{\to}& M\otimes_BJ\otimes_BM^{op}& \stackrel{\ffi}{\to}& I.
\end{array}$$

\subsection{Form functors as tensor product with symmetric forms}
\label{subsec:FormFunAsTensProd}
Let $A$ and $B$ be dgas with involution, and let $(I,i)$, $(J,j)$ be duality
coefficients for $A\Lmod$ and $B\Lmod$, respectively.
We want to think of (certain) form functors $(B\Lmod,\sharp_J) \to
(A\Lmod,\sharp_I)$ as tensor product with symmetric forms in
$A\LRmod B$.  
For that, let $(M,\ffi)$ be a symmetric form in $A\LRmod B$.
It defines a form functor
$$(M,\ffi)\otimes_B ? : (B\Lmod,\sharp_J)
\stackrel{(F,\Phi)}{\longrightarrow} (A\Lmod,\sharp_I),$$
where $F(P)=M\otimes_BP$ and the duality
compatibility map is the left $A$-module homomorphism
$$M \otimes_B [P^{op},J]_B \stackrel{\Phi_P}{\longrightarrow}
[(M\otimes_B P)^{op},I]_A$$ 
defined by
$$\Phi(x\otimes f)((y\otimes t)^{op}) = (-1)^{|y||t|}\ffi(x\otimes
f(t^{op})\otimes y^{op})$$
for $x, y\in M$, $f\in [P^{op},J]_B$, and $t\in P$.

\subsection{Basic properties of $(M,\ffi)\otimes_B?$}
\label{subsec:basicMffiOtimes}
Let $A$, $B$, $C$ be dg algebras with involution,
and let 
$(I,i)$, $(J,j)$, $(K,k)$ be duality coefficients for $A\Lmod$, $B\Lmod$ and
$C\Lmod$, respectively.
Further, let $(M,\ffi)$, $(M',\ffi')$  be symmetric forms in $A\LRmod B$ and
$(N,\psi)$ a symmetric form in  $B\LRmod C$.
Form functors induced by tensor
product with symmetric forms have the following elementary
properties.
\begin{enumerate}
\label{subsec:basicMffiOtimesA}
\item
Tensor product $(A,\bmu_I)\otimes_A?$ with the symmetric form 
$$\bmu_I: A \otimes_A I \otimes_A A^{op} \to
I: a \otimes t \otimes b \mapsto a\cdot t\cdot \bar{b}$$ 
on the $A$-bimodule $A$ induces the identity form functor on $(A,\sharp_I)$.
\item
An isometry $(M,\ffi) \cong (M',\ffi')$ between symmetric forms in $A\LRmod B$
defines an isometry of associated form functors 
$$(M,\ffi)\otimes_B ? \cong
(M',\ffi')\otimes_B ?: (B\Lmod,\sharp_J) \to (A\Lmod,\sharp_I).$$
\item
Orthogonal sum $(M,\ffi) \perp (M',\ffi')$ of symmetric forms in $A\LRmod B$
corresponds to orthogonal sum of associated form functors:
$$(M,\ffi)\otimes_B ?\hspace{1ex} \perp \hspace{1ex} (M',\ffi')\otimes_B ?
\hspace{1ex}\cong\hspace{1ex} [(M,\ffi)\perp (M',\ffi')]\otimes_B ?$$
\item
Tensor product of symmetric forms corresponds to composition of
associated form functors:
$$[(M,\ffi)\otimes_B (N,\psi)]\otimes_C ? \cong [(M,\ffi)\otimes_B ?] \circ [
(N,\psi) \otimes_C ?]$$
\end{enumerate}
These properties follow directly from the definitions, and we omit the details.

\subsection{Tensor product of dgas with involution}
\label{subsec:AotimesB}
For two dgas $A$, $V$, the tensor product dg $\kappa$-module
$AV=A\otimes_{\kappa}V$ is a dga with multiplication $(a\otimes v)\cdot
(b\otimes w) = (-1)^{|b||v|}(a\cdot b) \otimes (v\cdot w)$.
If $A$, $V$ are dgas with involution, then the tensor product dga $AV$ is a
dga with involution $AV \to (AV)^{op}:a\otimes v \mapsto (\bar{a}\otimes \bar{v})^{op}$.
Furthermore, if $(I,i)$ and $(U,u)$ are duality coefficients for $A\Lmod$ and
$V\Lmod$, then $IU=I\otimes U$ is an $AV$-bimodule with left multiplication 
$(a\otimes v)\cdot (t\otimes x)= (-1)^{|v||t|}a\cdot t \otimes v\cdot x$ and
right multiplication $(t\otimes x)\cdot (a\otimes v) = (-1)^{|a||x|}t\cdot a
\otimes x\cdot v$, for $a\in A$, $v \in V$, $t \in I$, $x \in U$, and the
$AV$-bimodule map $iu=i\otimes u:
I\otimes U \to (I\otimes U)^{op}: t \otimes x \mapsto (i(t)\otimes u(x))^{op}$
makes the pair $(IU, iu)$ into a duality coefficient for $AV\Lmod$.

If $B$ is another dga with involution, and $(J,j)$ is a duality coefficient
for $B\Lmod$, then, with the same formulas as in \ref{subsec:symmAmodB},
\ref{subsec:FormFunAsTensProd}, any symmetric from $(M,\ffi)$ in $A\LRmod B$
with respect to the duality coefficients $(I,i)$ and $(J,j)$
defines a form functor 
$$(M,\ffi) \otimes_B ?:( BV\Lmod, \sharp_{(JU,ju)}) \to ( AV\Lmod,
\sharp_{(IU,iu)})$$
satisfying the properties in \ref{subsec:basicMffiOtimes}.

\subsection{Extension to ringed spaces}
\label{subsec:ExtRngSp}
Let $(X,O_X)$ be a ringed space with $O_X$ a sheaf of commutative rings on
a topological space $X$.
Replacing in \ref{subsec:DGk} - \ref{subsec:AotimesB} the ground ring
$\kappa$ with the sheaf of commutative rings $O_X$, all definitions and
properties from \ref{subsec:DGk} - \ref{subsec:AotimesB} extend to
modules over differential graded sheaves of $O_X$-algebras.
Definitions are extended by applying the definitions of \ref{subsec:DGk} -
\ref{subsec:AotimesB} to sections over open subsets of $X$.
For instance, let $\A$ be a sheaf of dg $O_X$-algebras with involution, 
$(I,i)$ a duality coefficient for $\A\Lmod$, and
$P$ a sheaf of left dg $\A$-modules.
The canonical double dual identification $\can: P \to
[[P^{op},I]_{\A}^{op},I]_{\A}$ is defined by sending a section $x \in P(U)$, $U
\subset X$, to the map of sheaves of dg modules
$\can(x): ([P^{op},I]_{\A}^{op})_{|U} \to I_{|U}$ defined on $V \subset U$ by
$\can(x)(f^{op})=(-1)^{|x||f|}i(f(x^{op}_{|V}))$ for $f\in
[P^{op},I]_{\A}(V)$.

\section{Higher Grothendieck-Witt groups of schemes}
\label{sec:HigherGWschemes}

Let $X$ be a scheme, $\A_X$ be a quasi-coherent sheaf of $O_X$-algebras with
involution, $L$ a line bundle on $X$, $Z\subset X$ a
closed subscheme and $n \in \Z$ an integer.
The purpose of this section is to introduce the
Grothendieck-Witt space 
$GW^n(\A_X\phantom{i}on\phantom{i}Z,\phantom{i}L)$
of symmetric spaces over $\A_X$ with coefficients in the $n$-th shifted
line bundle $L[n]$ and support in $Z$. 
We work in this generality 
in order to be able to extent the localization and excision theorems of
\S\ref{sec:LocnZar} to negative degrees.  

Recall that, unless otherwise indicated, ``module'' will always mean ``
left module''. 
In what follows, we will denote by $\otimes$ the tensor product
$\otimes_{O_X}$ of $O_X$-modules.

\subsection{Vector bundles and strictly perfect complexes} 
Let $\A_X$ be a quasi-coherent sheaf of $O_X$-algebras with involution.
The category of quasi-coherent left
$\A_X$-modules (dg-modules concentrated in degree $0$) is a fully exact
abelian subcategory of the abelian category of left $\A_X$-modules. 
We denote by 
$$\Vect(\A_X)$$
the full subcategory of $\A_X$ vector bundles, that is, of those
quasi-coherent left $\A_X$-modules $F$ for which $F(U)$ is a finitely generated
projective $\A_X(U)$-module for every affine $U \subset X$. 
As usual, the last condition only needs to be checked for those $U$ running
through a choice of an affine open cover of $X$.
The category of $\A_X$ vector bundles inherits the notion of exact sequences
from the category of all (quasi-coherent) $\A_X$-modules.
Note that $\A_X$ vector bundles need not be locally free, since $\A_{X,x}$ may
not be commutative nor a local ring for $x \in X$.
In case $\A_X=O_X$, the category $\Vect(X)$ is the usual exact category of
vector-bundles on $X$.

A {\em strictly perfect} complex of $\A_X$-modules is a
dg left $\A_X$-module $M$ such that  $M_n=0$ for
all but finitely many $n\in \Z$ and $M_n$ is an $\A_X$ vector bundle for all
$n \in \Z$. 
Denote by $\sPerf(\A_X)$ the category of strictly perfect complexes of
$\A_X$-modules, in oder words, the category of bounded chain complexes of
$\A_X$ vector bundles.

Let $\L$ be a line bundle on $X$.
Then $\A_XL[n]= \A_X\otimes L[n]$ is a dg $\A_X$-bimodule via the
multiplication defined on sections by
$a(x\otimes l)b= (-1)^{|b||l|}axb\otimes l$ for $a,b,x \in \A_X$ and $l\in L[n]$.
We equip $\A_X L[n]$ with a dg $\A_X$-bimodule isomorphism
$i:\A_X L[n] \to (\A_X L[n])^{op}: a\otimes l \mapsto
\bar{a}\otimes l$ satisfying $i^{op}\circ i=1$, so that $(\A_XL[n],i)$
is a duality coefficient for $\A_X\Lmod$.
If $\eps\in \{+1,-1\}$, then $(\A_XL[n],\eps i)$ is also a duality coefficient
for $\A_X\Lmod$.
In the notation of \ref{subsec:AotimesB}, the duality coefficient
$(\A_XL[n],\eps i)$ is the tensor 
product of the duality coefficient $(\A_X,\mu)$ for $\A_X\Lmod$ and
the duality coefficient $(L[n],\eps)$ for $O_X\Lmod$.

For a strictly perfect complex of $\A_X$-modules $M$, the left dg
$\A_X$-module 
$$M^{\sharp_{\eps L}^n}=[\, M^{op},\A_XL[n]\, ]_{\A_X}$$
is also strictly perfect,
the functor $M \mapsto M^{\sharp_{\eps L}^n}$ is exact and preserves
quasi-isomor-phisms.
Moreover, the double dual identification $\can_{(\A_XL[n],\eps i)}$ 
defined in \ref{subsec:DGmodDual} is an isomorphism.
Therefore, the triple
$$(\sPerf(\A_X),\quis,\sharp_{\eps L}^n)$$
defines an exact category with weak equivalences and duality, the double dual
identification being understood as $\can_{(\A_XL[n],\eps i)}$.
If $n=0$ (or $\eps = 1$, or $L=O_X$), we may omit the label corresponding to
$n$ (or $\eps$, or $L$, respectively), so that $(A_X\Lmod, \sharp^n)$ means
$(A_X\Lmod, \sharp_{1,O_X}^n)$, for instance.
By restriction of structure, we have an exact category with duality
$$(\Vect(\A_X),\sharp_{\eps L}).$$

Let $Z \subset X$ be a closed subscheme with open complement $U=X-Z$.
A strictly perfect complex $M$ of $\A_X$-modules has {\em cohomological
  support in $Z$} if the 
restriction $M_{|U}$ of $M$ to $U$ is acyclic.
We write $\sPerf(\A_X\phantom{i}on\phantom{i}Z)$ for the category
of strictly perfect complexes of $\A_X$-modules which have cohomological
support in $Z$. 
By restriction of structure, we have exact categories with weak equivalences
and duality 
\begin{equation}
\label{eqn:sPerfWD}
(\, \sPerf(\A_X\phantom{i}on\phantom{i}Z),\, \quis,\, \sharp_{\eps L}^n,\, ).
\end{equation}

\subsection{Definition}
\label{dfn:GWnXonZ}
Let $X$ be a scheme, $\A_X$ be a quasi-coherent sheaf of $O_X$-algebras with
involution, $L$ a line bundle on $X$, $Z\subset X$ a closed subscheme, $n
\in \Z$ an integer, and $\eps \in \{+1,-1\}$.
The Grothendieck-Witt space 
$${_{\eps}GW}^n(\A_X\phantom{i}on\phantom{i}Z,\phantom{i}L)$$
of $\eps$-symmetric spaces over $\A_X$ with coefficients in the $n$-th shifted
line bundle $L[n]$ and (cohomological) support in $Z$ is
the Grothendieck-Witt space of the
exact category with weak equivalences and duality (\ref{eqn:sPerfWD}).
If $Z=X$ (or $\eps =1$, or $L=O_X$, or $n=0$), we may omit the label
corresponding to $Z$ ($\eps$, $L$, $n$, respectively).
For instance, the space $GW(\A_X,L)$ denotes the Grothendieck-Witt space 
${_{1}GW}^0(\A_X\phantom{i}on\phantom{i}X,\phantom{i}L)$.

\subsection{Remark}
By \ref{subsec:canSymCone}, the exact category with weak equivalences and
duality (\ref{eqn:sPerfWD}) has a symmetric cone in the sense of
\ref{dfn:symCone}. 
\vspace{2ex}

In the following proposition, we write $\mu:O_X\otimes O_X \to O_X$ for the
multiplication in $O_X$. 

\subsection{Proposition}
\label{prop:Periodicity}
{\it
Tensor product with the $(-1)$-symmetric space $(O_X[1],\mu)$ induces an
equivalence of exact categories with weak equivalences and duality 
$$(\, \sPerf(\A_X\phantom{i}on\phantom{i}Z),\, \quis,\, \sharp_{\eps L}^n,\,
) \cong (\, \sPerf(\A_X\phantom{i}on\phantom{i}Z),\, \quis,\,
\sharp_{-\eps L}^{n+2},\, ).$$
In particular, we have homotopy equivalences of Grothendieck-Witt spaces
$$
\renewcommand\arraystretch{1.5} 
\begin{array}{llcl}
(O_X[1],\mu)\otimes ?: &
{_{\eps}GW}^n(\A_X\phantom{i}on\phantom{i}Z,\phantom{i}L) & \simeq &
{_{-\eps}GW}^{n+2}(\A_X\phantom{i}on\phantom{i}Z,\phantom{i}L)\hspace{4ex}{\rm
  and}\\
(O_X[2],\mu)\otimes ?:&
{_{\eps}GW}^n(\A_X\phantom{i}on\phantom{i}Z,\phantom{i}L) & \simeq &
{_{\eps}GW}^{n+4}(\A_X\phantom{i}on\phantom{i}Z,\phantom{i}L).
\end{array}$$
}

\subsection*{\it Proof}
The pair $(O_X[1],\mu)$ defines a symmetric space in $O_X\Lmod$ with respect
to the duality coefficient $(O_X[2],-1)$.
Tensor product $(O_X[1],\mu)\otimes_{O_X} ?$ defines a
form functor 
$(\A_X\Lmod,\sharp_{\eps L}^n) \to
  (\A_X\Lmod, \sharp_{-\eps L}^{n+2})$  
as explained in \ref{subsec:FormFunAsTensProd} -- \ref{subsec:ExtRngSp}.
Since $(O_X[1],\mu)\otimes_{O_X}(O_X[-1],-\mu)$ and
$(O_X[-1],-\mu)\otimes_{O_X}(O_X[1],\mu)$ are isometric to $(O_X,\mu)$
which induces the identity form functor, the equivalence of categories with
duality and the first homotopy equivalence follow.
The second map of spaces is a homotopy equivalence 
with inverse given by the tensor product with the symmetric space
$(O_X[-2],\mu)$.
\Qed

\subsection{Corollary}
\label{cor:sPerf=exVect}
{\it 
For $n\in \Z$,
there are functorial homotopy equivalences
$$
\renewcommand\arraystretch{1.5} 
\begin{array}{ccl}
GW^{4n}(\A_X,L) & \simeq &
GW(\Vect(\A_X),\sharp_L,\can_L),\phantom{a}{\rm \it and}\\
GW^{4n+2}(\A_X,L) & \simeq & GW(\Vect(\A_X),\sharp_L,-\can_L).
\end{array}
$$
where the Grothendieck-Witt spaces on the right hand side are the 
ones associated with the exact categories with duality
$(\Vect(\A_X),\sharp_L,\pm\can_L)$
as defined in \cite[4.6]{myGWex}.}

\subsection*{\it Proof}
The homotopy equivalences follow from proposition \ref{prop:gilletWald},
remark \ref{rmk:ChoiceOfSign}, and proposition \ref{prop:Periodicity}.
\Qed

\section{Localization and Zariski-excision in positive degrees}
\label{sec:LocnZar}

\subsection{Schemes with an ample family of line bundles}
A scheme $X$ has
an {\em ample family of line bundles}  if there is a finite set
$L_1,...,L_n$ of line bundles with global sections 
$s_i\in \Gamma(X,L_i)$ such that the non-vanishing loci $X_{s_i}=\{x\in
  X|\phantom{a}s_i(x)\neq 0\}$ form an open affine cover of $X$, see
\cite[Definition 2.1]{TT}, \cite[II 2.2.4]{SGA6}.

Recall that if  $f\in \Gamma(X,L)$ is a global section of a line bundle $L$ on
a scheme 
$X$, then the open inclusion $X_f \subset X$ is an affine map (as can be seen by
choosing an open affine cover of $X$ trivializing the line bundle $L$).
As a special case, $X_f$ is affine whenever $X$ is affine.
Thus, for the affine cover above $X = \bigcup X_{s_i}$, all finite
intersections of the $X_{s_i}$'s are affine.
In particular, a scheme with an ample family of line bundles is quasi-compact
(as a finite union of affine, hence 
quasi-compact, subschemes) and it is quasi-separated.
Recall that the latter means that the intersection of
any two quasi-compact open subsets is quasi-compact 
(a condition which only needs to be checked for the pair-wise intersections
$U_i\cap U_j$ 
of a cover of $X = \bigcup_iU_i$ by quasi-compact open subsets $U_i$; in
our case, we can take $U_i=X_{s_i}$). 

For a scheme $X$ with an ample family of line bundles, there is a set $L_i$,
$i \in I$, of line
bundles on $X$ with global sections $s_i \in \Gamma(X,L_i)$ such that the open
subsets $X_{s_i}$, $i \in I$, form an open affine basis for the topology of
$X$ \cite[2.1.1 (b)]{TT}.
\vspace{2ex}

For examples of schemes with an ample family of line bundles, see
\cite[2.1.2]{TT}.
Any quasi-compact open or closed subscheme of a scheme with an
ample family of line bundles has itself an ample family of line bundles.
Any scheme 
quasi-projective over an affine scheme, and any separated regular noetherian
scheme has an ample family of line bundles.
\vspace{2ex}

The main purpose of this section is to prove the following two theorems.

\subsection{Theorem \rm (Localization)}
\label{thm:Locn1}
{\it
Let $X$ be a scheme with an ample family of line-bundles, let $Z \subset X$
be a closed subscheme with quasi-compact open complement
 $j:U \subset X$, and let $L$ a line bundle on $X$.
Let $\A_X$ be a quasi-coherent sheaf of $O_X$-algebras with involution.
Then for every $n\in \Z$ there is a homotopy fibration of Grothendieck-Witt
spaces 
$$GW^n(\A_X\phantom{i}on\phantom{i}Z,\phantom{i}L) \longrightarrow
GW^n(\A_X,\phantom{i}L) \longrightarrow GW^n(\A_U,\phantom{i}j^*L).$$ 
}

\subsection{Theorem \rm (Zariski-excision)}
\label{thm:ZarExc1}
{\it
Let $X$ be a scheme with an ample family of line-bundles, let $Z \subset X$
be a closed subscheme with quasi-compact open complement,
let $L$ be a line bundle on $X$ and
let $\A_X$ be a quasi-coherent sheaf of $O_X$-algebras with involution.
Then for every $n\in \Z$ and every quasi-compact open subscheme $j:V \subset X$
containing $Z$,
restriction of vector-bundles induces a homotopy equivalence 
$$GW^n(\A_X\phantom{i}on\phantom{i}Z,\phantom{i}L)
\stackrel{\sim}{\longrightarrow}
GW^n(\A_V\phantom{i}on\phantom{i}Z,\phantom{i}j^*L).$$
}
\vspace{2ex}

The proofs of theorems \ref{thm:Locn1} and \ref{thm:ZarExc1} will occupy the
rest of this section.
First, we introduce some terminology.
For an open subscheme $U \subset X$, call a map $M \to N$ of left dg
$\A_X$-modules a $U$-isomorphism ($U$-quasi-isomorphism) if its restriction 
$M_{|U} \to N_{|U}$ to $U$ is an isomorphism (quasi-isomorphism).
A left dg $\A_X$-module $M$ is called $U$-acyclic if its restriction 
$M_{|U}$ to $U$ is acyclic.

\subsection{Notation}
\label{subsec:notation1}
In \ref{subsec:truncKosz}, \ref{lem:CalcFracByKoszul}, \ref{lem:UquisAppr}
below, we consider the following objects:
\begin{enumerate}
\item
a scheme $X$ which is quasi-compact and quasi-separated,
\item
a finite set of line bundles $L_i$, $i=1,...,n$ together with global sections
$s_i \in \Gamma(X,L_i)$,
\item
the union $U=\bigcup_{i=1}^nX_{s_i}$ of the non-vanishing loci
$X_{s_i}$ of the $s_i$'s, denoting $j:U \subset X$ the corresponding open
immersion, and
\item
a quasi-coherent sheaf of $O_X$-algebras $\A_X$.
\end{enumerate}

\subsection{Truncated Koszul complexes}
\label{subsec:truncKosz}
In the situation of \ref{subsec:notation1}, 
the global sections $s_i$ define maps $s_i:O_X \to L_i$ of line-bundles whose
$O_X$-duals are denoted by $s_i^{-1}:L_i^{-1} \to O_X$.
We consider the maps $s_i^{-1}$ as
(cohomologically graded) chain-complexes with $O_X$ placed in degree $0$. 
For an $l$-tuple $n=(n_1,...,n_l)$ of negative integers, the Koszul complex 
\begin{equation}
\label{eqn:Koszul}
\bigotimes_{i=1}^l(L_i^{n_i} \stackrel{s^{n_i}}{\to} O_X)
\end{equation}
is acyclic over $U$, since the map 
$s^{n_i}=(s_i^{-1})^{\otimes  |n_i|}:L^{n_i}_i \to O_X$ 
is an $X_{s_i}$-isomorphism, so that
the Koszul complex (\ref{eqn:Koszul}) is acyclic (even contractible) over each
$X_{s_i}$.
Let $K(s^n)$ denote the bounded complex whose non-zero part (which we place in
degrees $-l+1,...,0$) is the Koszul complex (\ref{eqn:Koszul}) in degrees
$\leq -1$. 
The last differential $d_{-1}$ of the Koszul complex defines a map
$$K(s^n)=\left[\bigotimes_{i=1}^l(L_i^{n_i} \stackrel{s^{n_i}}{\to}
  O_X)\right]_{\leq -1}[-1] \stackrel{\eps}{\longrightarrow} O_X$$
of strictly perfect complexes of $O_X$-modules which is a $U$-quasi-isomorphism
since its cone, the Koszul complex, is $U$-acyclic.
For a left dg $\A_X$-module $M$, we write
$\eps_M$ for the tensor product map $\eps_M=1_M \otimes \eps: M \otimes
K(s^n) \to M\otimes O_X \cong M$.
\vspace{2ex}

The following lemma is a consequence of the well-known techniques of
extending a 
section of a quasi-coherent sheaf from an open subset cut out by
a divisor to the scheme itself \cite[Th\'eor\`eme 9.3.1]{EGAI}.
It is implicit in the proof of \cite[Proposition 5.4.2]{TT}.

\subsection{Lemma}
\label{lem:CalcFracByKoszul}
{\it
In the situation \ref{subsec:notation1},
let $M$ be a complex of quasi-coherent left $\A_X$-modules and let $A$ be a
strictly perfect complex of $\A_X$-modules. Then the following holds.
\begin{itemize}
\item[(a)]
For every map $f:j^*A \to j^*M$ of left dg $\A_U$-modules between the
restrictions of $A$ and $M$ to 
$U$, there
is an $l$-tuple of negative integers $n=(n_1,...,n_l)$, and
a map $\tilde{f}: A \otimes K(s^{n}) \to M$ of left dg $\A_X$-modules  such
that $f\circ j^*(\eps_A) = j^*(\tilde{f})$. 
\item[(b)]
For every map $f:A \to M$ of left dg $\A_X$-modules such that $j^*(f)=0$,
there is an $l$-tuple of 
negative integers $n=(n_1,...,n_l)$ such that $f\circ \eps_A=0$.
\end{itemize}
}

\subsection*{\it Proof}
For any complex of $O_X$-modules $K$, to give a map $A\otimes K \to M$ of
chain complexes of $\A_X$-modules 
is the same as to give a map $K \to {_{\A_X}[A,M]}$ of chain complexes of
$O_X$-modules, by adjointness of the tensor product $A \otimes ? :O_X\Lmod \to
\A_X\Lmod$ and the left $\A_X$-module map functor 
${_{\A_X}[A,\phantom{M}]}: \A_X\Lmod \to O_X\Lmod$.
Note that if $A$ is 
strictly perfect, the complex ${_{\A_X}[A,M]}$ is a complex of quasi-coherent
$O_X$-modules and the 
natural map ${_{\A_X}[A,j_*j^*M]} \to j_*j^*\ {_{\A_X}[A,M]}$ is an isomorphism.
The adjunction 
allows us to reduce the proof of the lemma to $\A_X=O_X$ and $A=O_X$
(concentrated in degree $0$).

Every map $O_X \to M$ of chain complexes of $O_X$-modules 
factors through the subcomplex $Z_0M \subset M$ of $M$ which is the complex
$\ker(d_0:M_0 \to M_1)$ concentrated in degree $0$.
By adjunction, every map $O_X \to j_*j^*M$ factors as $O_X \to
j_*Z_0j^*M=j_*j^*Z_0M \to j_*j^*M$.
This allows us to further reduce the proof to $M$ a complex with $M_i=0$,
$i\neq 0$.
In this case, the proof for $l=1$ is classical, see \cite[Th\'eor\`eme
9.3.1]{EGAI},  \cite[Lemma 5.4.1]{TT}, so that the lemma is proved in case
$l=1$.

Before we treat the case $l>1$ (and $\A_X=O_X$; $A=O_X$, $M$ concentrated in
degree $0$), we prove the following.
\begin{itemize}
\item[(\dag)]
For every map $A \to G$ of complexes of quasi-coherent $O_X$-modules with $A$
strictly perfect and $j^*G$ contractible, there is an $l$-tuple of negative
integers $(n_1,...,n_l)$ and a commutative diagram of complexes of
$O_X$-modules
$$\xymatrix{A{\phantom{A}} \ar[r]^{\hspace{-12ex}1_A \otimes \iota} \ar[dr] &
 {\phantom{A}} A \otimes  \left[\, \bigotimes_i^l\,
   (L_i^{n_i}\stackrel{s^{n_i}}{\to} O_X)\, \right]
  \ar[d]\\ 
& G,}$$
where $\iota$ is the canonical embedding of $O_X$ (concentrated in degree $0$)
into the Koszul complex.
\end{itemize}
It suffices to prove (\dag) for $l=1$, since the general case is a repeated
application of the case $l=1$.
For the proof of (\dag) with $l=1$, we can assume $A=O_X$, as above.
The composition $O_X \to G \to j_*j^*G$ has contractible target and therefore
factors through the cone $(O_X \to O_X)$ of $O_X$.
By the case $l=1$ of the lemma (proved above), there is an $n<0$ such that
the composition $(L^{n} \to L^{n}) \stackrel{s^{n}}{\to} (O_X \to O_X)
\to j_*j^*G$ lifts to $G$. 
The two maps $(0 \to L^{n}) \stackrel{s^{n}}{\to}  (0 \to O_X) \to G$
and  $(0 \to L^{n}) \to (L^{n} \to L^{n}) \to G$ may not be the same,
but their compositions with $G \to j_*j^*G$ are, so that, again by case $l=1$
of the lemma, precomposing both maps with $L^{n+t} \stackrel{s^{t}}{\to}
L^{n}$ makes the two maps with target $G$ equal.
Replacing $n$ with $n+t$, we can assume that the two maps $(0 \to L^{n})
\to G$ above coincide.
Then we obtain a map from the push-out $(L^{n} \stackrel{s^{n}}{\to} O_X)$
of $(L^n \to L^n) \leftarrow (0 \to L^n) \stackrel{s^{n}}{\to} (0 \to O_X)$ to
$G$.
This proves (\dag) for $l=1$, hence for all $l$.

For part (a) of the general case of lemma \ref{lem:CalcFracByKoszul} (and
$\A_X=O_X$; $A=O_X$, $M$ 
concentrated in degree $0$), we apply (\dag) to 
the map of chain-complexes of $O_X$-modules
$(0 \to O_X) \to (M \to j_*j^*M)$, and obtain a factorization of that map
through the Koszul complex $\bigotimes_i(L_i^{n_i}\stackrel{s^{n_i}}{\to}
O_X)$ for some $l$-tuple of negative integers $(n_1,...,n_l)$.
The canonical map from the stupid truncation in degrees $\leq -1$ (shifted by
$1$ degree) to its degree $0$ part 
$$\xymatrix{
K(s^n) = \left[\bigotimes_{i=1}^l(L_i^{n_i} \stackrel{s^{n_i}}{\to}
  O_X)\right]_{\leq -1}[-1] \ar[d] \ar[r]^{\hspace{16ex}d_{-1}} & O_X \ar[d]\\
M=[M \to j_*j^*M]_{\leq -1}[-1] \ar[r]^{\hspace{8ex}d_{-1}} & j_*j^*M}
$$
yields (a).
The general case of part (b) is a repeated application of the case $l=1$.
\Qed
\vspace{2ex}

For an exact category with weak equivalences $(\C,w)$, we write 
$D(\C,w)$ for its derived category, that is, the category $\C[w^{-1}]$
obtained from $\C$ by formally inverting the arrows in $w$.
If $(\C,w)$ is a category of complexes in the sense of definition
\ref{dfn:CatOfCxs}, its derived category $D(\C,w)$ is a triangulated category.
In this case, it can also be obtained as the localization by a calculus of
fractions of the homotopy category $\scK(\C)$ of $\C$ which is the factor
category of $\C$ modulo the ideal of maps which are homotopic to zero.

\subsection{Lemma}
\label{lem:UquisAppr}
{\it
In the situation \ref{subsec:notation1},
let $\A \subset \sPerf(\A_X)$ be a full subcategory of the
category of strictly perfect $\A_X$-modules such that the inclusion
$\A \subset \sPerf(\A_X)$ is closed under degree-wise split
extensions, usual shifts and cones.
Assume furthermore that for all $A \in \A$, $k \leq 0$
and $i=1,...n$, we have
$A \otimes L_i^{k} \in \A$.

Then for every $U$-quasi-isomorphism
$M \to A$  of complexes of quasi-coherent
$\A_X$-modules with $A \in \A$,
there is a $U$-quasi-isomorphism $B \to M$
of complexes of $\A_X$-modules with $B \in \A$:
$$\xymatrix{B \ar[r]^{\exists}_{U-quis}& M \ar[r]^{\forall}_{U-quis} & A.}$$
\vspace{1ex}

\noindent
In particular, the inclusion $\A \subset \sPerf(\A_X)$
 induces a fully faithful triangle functor
$$\D(\, \A,\, U\Quis) \subset \D(\, \sPerf(\A_X),\, U\Quis).$$
}

\subsection*{\it Proof}
We first prove the following statement.
\vspace{2ex}

\hspace{-7ex}
\parbox{5in}{
\begin{itemize}
\item[(\dag)]
Let $s \in \Gamma(X,L)$ be a global section of a line-bundle $L$ such that $X_s$
is affine.
Then for every $X_s$-quasi-isomorphism $N \to E$ of complexes of quasi-coherent
$\A_X$-modules with $E$ strictly perfect on $X$, there is an
$X_s$-quasi-isomorphism $E\otimes L^{-k} \to N$ for some integer $k>0$.
\end{itemize}
}
\vspace{2ex}

\noindent
Write $j: X_s \subset X$ for the open inclusion.
Since $X_s$ is affine, we have an equivalence of categories between
quasi-coherent $\A_{X_s}$-modules 
and $\A_X(X_s)$-modules under which the map $j^*N \to j^*E$ becomes a
quasi-isomorphism of 
complexes of $\A_X(X_s)$-modules with $j^*E$ a bounded complex of projectives. 
Such a map always has a section up to homotopy $f: j^*E \to j^*N$ which  
is then a quasi-isomorphism.
By \ref{lem:CalcFracByKoszul} with $l=1$, there is a
map of complexes $\tilde{f}:E\otimes L^{k} \to N$ such that
$j^*\tilde{f}=f\cdot s^k$, for some $k<0$.
In particular, $\tilde{f}$ is a $U$-quasi-isomorphism.

Now we prove the lemma by induction on $n$.
For $n=1$, this is (\dag).
Let $U_0=\bigcup_{i=1}^{n-1}X_{s_i}$.
By our induction hypothesis, there is a $U_0$-quasi-isomorphism $B_0 \to M$
with $B_0\in \A$.
Let $M_0$ and $A_0$ be the cones of the maps $B_0 \to M$ and $B_0 \to A$, that
is, the push-out of these maps along the canonical (degree-wise split)
injection of $B_0$ into its cone $CB_0$.
We obtain a commutative (in fact bicartesian) diagram involving
$M$, $M_0$, $A$ and $A_0$ with $M \to M_0$ and $A \to A_0$ degree-wise 
split injective.
Factor the map $A \to A_0$ as in the diagram
$$
\xymatrix{M \xymono[r] \ar[d] & M_1 \xyepi[r] \ar[d] & M_0 \ar[d]\\
A \xymono[r] &  A\oplus  PA_0 \xyepi[r] & A_0}
$$
with $A \oplus PA_0 \to A_0$ degree-wise split surjective and $PA_0 = CA_0[-1]
\in \Ch^b\A$ contractible.
Then $M \to M_0$ factors through the pull-back $M_1$ of $M_0 \to A_0$ along the
surjection $A \oplus PA_0 \to A$.
The map $M \to M_1$ is degree-wise split injective (as $M \to M_0$ is), and
has cokernel 
the contractible complex $PA_0$.
It follows that $M \to M_1$ is a homotopy equivalence, and we can choose a
homotopy inverse $M_1 \to M$.
By construction, $A_0$ and $M_0$ are acyclic over $U_0$ and
$A_0 \in \A$.
Moreover, the map $M_0 \to A_0$ is an $X_{s_n}$-quasi-isomorphism.
By (\dag), there is an $X_{s_n}$-quasi-isomorphism $A_0\otimes L^{-k} \to M_0$
for some $k>0$.
The complex $A_0\otimes L^{-k}$ is $U_0$-acyclic since $A_0$ is, so that 
the map $A_0\otimes L^{-k} \to M_0$ is in fact a $U$-quasi-isomorphism.
Let $B$ be the pull-back of $A_0\otimes L^{-k} \to M_0$ along the surjection
$M_1 \to M_0$.
The resulting map $B \to M_1$ is a $U$-quasi-isomorphism.
Moreover $B$ is an object of  $\A$
since $B$ is also the pull-back of $A_0\otimes
L^{-k} \to A_0$ along the (degree-wise split) surjection $A \oplus PA_0 \to
A_0$. 
Composing the map $B \to M_1$ with the homotopy inverse $M_1 \to M$ of $M \to
M_1$ yields the desired $U$-quasi-isomorphism $B \to M$.
\Qed

\subsection{The derived category of quasi-coherent $\A_X$-modules} 
Recall that a Grothendieck abelian category is an abelian category $\A$ in
which all set-indexed direct sums exist, filtered colimits are exact, and $\A$
has a set of generators.
We remind the reader that a set $I$ of objects of $\A$ generates the abelian
category $\A$ if for every object $E\in
\A$, there is a surjection $\bigoplus A_i \twoheadrightarrow E$ from a set
indexed direct sum of (possibly repeated) objects $A_i \in I$  to $E$.

If $X$ is a scheme with an ample family of line-bundles, and $\A_X$ a
quasi-coherent $O_X$-algebra, then the category $\Qcoh(\A_X)$ of
quasi-coherent $\A_X$-modules is an abelian category with generating set the
set $\A_X\otimes L_i^k$, $i=1,...,n$, $k\leq 0$, where $L_1,...,L_n$ is a set
of line bundles on $X$ with global sections $s_i \in \Gamma(X,L_i)$ such that
the non-vanishing loci $X_{s_i}$ are affine and cover $X$.

For a Grothendieck abelian category $\A$, write
$D(\A)$ for the unbounded derived category of $\A$, that is, the
triangulated category $D(\Ch\A,\quis)$.
This category has small homomorphism sets, by \cite[remark 10.4.5]{weibel:hom}.
Coproducts of complexes are also coproducts in $D\A$, so that the triangulated
category $D\A$ has all set-indexed coproducts.
This applies in particular to the unbounded derived category $D\Qcoh(\A_X)$ of
quasi-coherent $\A_X$ modules.
If $Z \subset X$ is a closed subset with quasi-compact open complement
$U=X-Z$, we 
write $D_Z\Qcoh(\A_X) \subset D\Qcoh(\A_X)$ for the full triangulated
subcategory of those complexes $E$ of quasi-coherent $\A_X$-modules whose
restriction $E_{|U}$ to $U$ are acyclic.

\subsection{Lemma}
\label{lem:DAgeneration}
{\it
Let $\A$ be a Grothendieck abelian category with generating set of objects
$I$.
Then, an object $E$ of the triangulated category $D\A$ is zero iff
every map $A[j] \to E$ in $D\A$ is the zero map for $A\in I$ and $j\in \Z$.
}

\subsection*{\it Proof}
Let $E$ be an object of the derived category $D\A$ of $\A$ such that every map
$A[j] \to E$ in $D\A$ is the zero map for $A\in I$ and $j\in \Z$.
We can choose a surjection $\bigoplus_J A_j \twoheadrightarrow \ker(d_0)$
in $\A$ with $A_j$ objects in the generating set $I$.
The inclusion of complexes $\ker(d_0) \to E$ yields 
a map of complexes $\bigoplus_J A_j  \to \ker(d_k) \to E$ which
induces a surjective map $\bigoplus_J A_j \twoheadrightarrow \ker(d_0)
\twoheadrightarrow H^0E$ on cohomology.
Since every map $\bigoplus_J A_j \to E$ is zero in $D\A$, the induced
surjective map $\bigoplus_J A_j \twoheadrightarrow H^0E$ is the zero map,
hence $H^0E=0$. 
The same argument applied to $E[k]$ instead of to $E$ shows that 
$H^kE=0$ for all $k\in \Z$, so that $E$ is
quasi-isomorphic to the zero complex.
\Qed
\vspace{2ex}

Next, we recall the concept of compactly generated triangulated categories
due to Neeman \cite{neeman:locn} in the form of
\cite{neeman:grothendieck}.

\subsection{Compactly generated triangulated categories}
Let $\T$ be a triangulated category in which (all set-indexed) coproducts
exists. 
An object $A$ of $\T$ is called {\em compact} 
\cite[Definition 1.6]{neeman:grothendieck} if the natural map
$\bigoplus_{j\in J}Hom(A,M_j) \to Hom(A,\bigoplus_{j\in J}M_j)$ is an
isomorphism for any set $M_j$, $j\in J$, of objects in $\T$.
The full subcategory $\T^c \subset \T$ of compact objects is closed under
shifts and cones and thus is a triangulated subcategory. 

A triangulated
category $\T$ is {\em compactly generated} 
\cite[Definition 1.7]{neeman:grothendieck} if $\T$ contains all set-indexed
direct sums, and if there is a set $I$ of compact objects in $\T$ such that an
object $M$ of $\T$ is the zero object iff all maps $A \to M$ are the zero map
for $A \in I$.

A set $I$ of compact objects in a compactly generated triangulated category
$\T$ is called a {\it generating set} \
\cite[Definition 1.7]{neeman:grothendieck} if $I$ is closed under shifts and
if  an
object $M$ of $\T$ is the zero object iff all maps $A \to M$ are the zero map
for $A \in I$.

The following theorem is due to Neeman 
\cite[Theorem 2.1]{neeman:grothendieck}.

\subsection{Theorem \rm (Neeman)}
\label{thm:neeman}
{\it 
\begin{enumerate}
\item
\label{thm:neemanA}
Let $\T$ be a compactly generated triangulated with generating set of
objects $I$.
Then the full triangulated subcategory $\T^c$ of compact objects in $\T$ is
the smallest idempotent complete triangulated subcategory of $\T$ containing
$I$.
\item
\label{thm:neemanB}
Let $\scR$ be a compactly generated triangulated category, $S \subset
\scR^c$ a set of compact objects closed under taking shifts.
Let $\scS \subset \scR$ be the smallest full triangulated subcategory closed
under formation of coproducts in $\scR$ which contains $S$.
Then $\scS$ and $\scR/\scS$ are compactly generated triangulated categories
with generating set $S$ and the image of $\scR^c$ in $\scR/\scS$.
Moreover, the functor $\scR^c/\scS^c \to \scR/\scS$ induces an equivalence
between the idempotent completion of $\scR^c/\scS^c$ and the category of
compact objects in $\scR/\scS$.
\item
\label{thm:neemanC}
Let $\scS \to \scR$ be a triangle functor between compactly generated
triangulated categories which preserves coproduct and compact objects.
Then $\scS \to \scR$ is an equivalence iff the functor $\scS^c \to \scR^c$
on compact objects is an equivalence.
\end{enumerate}
}
\vspace{2ex}

The following two propositions are essentially due to Thomason \cite{TT}.
We include the proofs here because only the commutative situation is considered
in \cite{TT}, and we need the explicite versions below.

For an exact category $\E$, write $D^b(\E)$ for the triangulated category
$D(\Ch^b\E,\quis)$. 
Recall that a fully faithful functor $\A \to 
\B$ of additive categories is called {\it cofinal} if every object of $B$ is a
direct factor of an object of $\A$.

\subsection{Proposition}
\label{prop:DXgenByLineBdls}
{\it
Let $X$ be a quasi-compact and quasi-separated scheme which is the union
$X=\bigcup_{i=1}^{n}X_{s_i}$ of open affine non-vanishing loci $X_{s_i}$ of global
sections $s_i \in \Gamma(X,L_i)$ of line-bundles $L_i$, $i=1,...,n$.
Let $\A_X$ be a quasi-coherent $O_X$-algebra.
Then the triangulated category $D\Qcoh(\A_X)$ is compactly generated by the
set of objects $\A_X\otimes L^k_i[j]$ for $k\leq 0$, $i=1,...,n$ and $j\in
\Z$.

Moreover, the inclusion $\Vect(\A_X) \subset \Qcoh(\A_X)$ induces a fully
faithful triangle functor $D^b\Vect(\A_X) \subset D\Qcoh(\A_X)$ which
identifies, up to equivalence, the category $D^b\Vect(\A_X)$ with 
the full triangulated subcategory $D^c\Qcoh(\A_X)$ of compact objects in $D\Qcoh(\A_X)$.
}

\subsection*{\it Proof}
In the triangulated category $D\Qcoh(\A_X)$, every strictly perfect complex of
$\A_X$ modules is a compact object. 
To see this, note that for an $\A_X$ vector bundle $A$ and a set $M_j$,
$j\in J$, of quasi-coherent $\A_X$-modules the canonical map of sheaves of
homomorphisms 
$\bigoplus_j\underline{Hom}_{\A_X}(A,M_j) \to
\underline{Hom}_{\A_X}(A,\bigoplus_j M_j)$
is an isomorphism since this can be
checked on an affine open cover of $X$ where the statement is clear.
Taking global sections, we obtain an isomorphism
$\bigoplus_jHom_{\A_X}(A,M_j) \stackrel{\cong}{\to}
Hom_{\A_X}(A,\bigoplus M_j).$
This isomorphism extends to an isomorphism of homomorphism sets of complexes
of $\A_X$-modules for $A$ a strictly perfect complex and $M$ an arbitrary
complex of quasi-coherent $\A_X$-modules.
For such complexes, the isomorphism induces an isomorphism
\begin{equation}
\label{eqn:pf:prop:DXgenByLineBdls}
\bigoplus_jHom_{\scK\Qcoh(\A_X)}(A,M_j) \stackrel{\cong}{\to}
Hom_{\scK\Qcoh(\A_X)}(A,\bigoplus_j M_j)
\end{equation}
of homomorphism sets in the homotopy category $\scK\Qcoh(\A_X)$ of chain
complexes of quasi-coherent $\A_X$-modules.
It follows from lemma \ref{lem:UquisAppr} with $U=X$ and
$\A=\sPerf(\A_X)$ that maps in $D\Qcoh(\A_X)$ from a strictly perfect complex
$A$ to an arbitrary complex $M$ of quasi-coherent $\A_X$-modules can be
computed as the filtered colimit
$$\colim_{B\stackrel{\sim}{\to}A}Hom_{\scK\Qcoh(\A_X)}(B,M)
\stackrel{\cong}{\longrightarrow}Hom_{D\Qcoh(\A_X)}(A,M)$$
of homomorphism sets in $\scK\Qcoh(\A_X)$ where the indexing category is the
left filtering category of homotopy classes of quasi-isomorphisms 
$B \stackrel{\sim}{\to} A$ of 
strictly perfect complexes with target $A$.
Taking the colimit over this filtering category of the isomorphism
(\ref{eqn:pf:prop:DXgenByLineBdls})
yields the isomorphism which proves that $A$ is compact in $\D\Qcoh(\A_X)$.

Since the set $\A_X\otimes L_i^k$, $i=1,...,n$, $k\leq 0$ is a set of
generators for the Grothendieck abelian category $\Qcoh(\A_X)$ all of which
are compact in the derived category $D\Qcoh(\A_X)$, lemma
\ref{lem:DAgeneration} shows that $D\Qcoh(\A_X)$ is a compactly generated
triangulated category with generating set $\A_X\otimes L_i^k[j]$, $i=1,...,n$,
$k\leq 0$, $j\in \Z$.

The inclusion $\Vect(\A_X) \subset \Qcoh(\A_X)$ of vector bundles into all
quasi-coherent $\A_X$-modules induces a triangle functor
$D^b\Vect(\A_X) \to D\Qcoh(\A_X)$ which is fully faithful, by the
existence of an ample family of line bundles and the criterion in
\cite[12.1]{keller:uses}. 
Since the exact category 
$\Vect(\A_X)$ is idempotent complete, its bounded
 derived category $D^b\Vect(\A_X)$ is also idempotent complete
\cite[Theorem 2.8]{balmerme}.
By Neeman's theorem \ref{thm:neeman} \ref{thm:neemanA}, the inclusion 
$D^b\Vect(\A_X) \to D^c\Qcoh(\A_X)$ is an equivalence.
\Qed

\subsection{Reminder on $Rj_*$}
Let $X$ be a scheme with an ample family of line bundles, and let $j:U
\hookrightarrow X$ be an open immersion from a quasi-compact open subset $U$
to $X$.
We recall one possible construction of the right-derived functor
$Rj_*:D\Qcoh(\U) \to D\Qcoh(X)$ of $j_*: \Qcoh(\U) \to \Qcoh(X)$.
For that, choose a finite cover $\mathscr{U}=\{U_0,...,U_n\}$ of $U$ such that 
the inclusion of all finite intersections $U_{i_0}\cap ... \cap U_{i_k}
\subset X$ into $X$ are affine maps, $i_0,...,i_k \in \{0,...,n\}$.
For instance, we can take as $\mathscr{U}$ an open cover of $U$ by a finite
number of non-vanishing loci
$X_{s_i}$ associated with a set of line bundles $L_i$ on $X$ and global
sections $s_i\in \Gamma(X,L_i)$, $i=0,...,n$.
For a $k+1$-tuple $\underline{i}=(i_0,...,i_k)$, set 
$U_{\underline{i}}=U_{i_0}\cap ... \cap U_{i_k}$ and write
write $j_{\underline{i}}:U_{\underline{i}} \subset U$ for the
  corresponding open immersion.

For a quasi-coherent $\A_U$ module $F$, 
consider the sheafified \v{C}ech complex $\check{C}(\mathscr{U},F)$ associated
with this covering.
In degree $k$ it is the quasi-coherent $\A_X$-module 
$$\check{C}(\mathscr{U},F)_k=\bigoplus_{\underline{i}}j_{\underline{i},*}j_{\underline{i}}^*F$$ 
where the indexing set is taken over all $k+1$-tuples
$\underline{i}=(i_0,...,i_k)$ with $0\leq i_0 < \cdots < i_k \leq n$.
The differential $d_k: \check{C}(\mathscr{U},F)_k \to
\check{C}(\mathscr{U},F)_{k+1}$ for the 
component $\underline{i}=(i_0,...,i_{k+1})$ is given by the
formula
$$(d_k(x))_{\underline{i}}=\sum_{l=0}^{k+1}(-1)^l\,
j_{\underline{i},*}j_{\underline{i}}^*\ 
x_{(i_0,...,\hat{i}_l,...,i_{k+1})}.$$
Note that the complex $\check{C}(\mathscr{U},F)$ is concentrated in degrees $0,...,n$.

The units of adjunction $F \to j_{i*}j_i^*F$ define a map $F \to
\check{C}(\mathscr{U},F)_0 = \bigoplus_{i=0}^n j_{i*}j_i^*F$ into the degree
zero part of the \v{C}ech complex with $d_0(F)=0$, and thus a map
of complexes of quasi-coherent $\A_X$-modules
$F \to \check{C}(\mathscr{U},F).$
This map is a quasi-isomorphism for any quasi-coherent $\A_X$-module $F$ as
can be checked by restricting the map to the open subsets $U_{i}$ of the
cover of $U$. 
Since, by assumption, for every $k+1$-tuple, $\underline{i} = (i_0,...,i_k)$,
the open inclusion  $j\circ j_{\underline{i}}: U_{\underline{i}}
\subset X$ is an affine map, the functor
$$j_*\check{C}(\mathscr{U}):\Qcoh(\A_X) \to \Ch\Qcoh(\A_X): F \mapsto
j_*\check{C}(\mathscr{U},F)$$ is exact.
Taking total complexes, this functor extends to a functor on all complexes
$$j_*\Tot\check{C}(\mathscr{U}):\Ch\Qcoh(\A_X) \to \Ch\Qcoh(\A_X): F \mapsto
j_* \Tot \check{C}(\mathscr{U},F) = \Tot j_* \check{C}(\mathscr{U},F) $$ 
which preserves quasi-isomorphisms as it is exact and sends acyclics to
acyclics. 
This functor is equipped with a natural quasi-isomorphism 
\begin{equation}
\label{eqn:Rj*NatQuis}
F \stackrel{\sim}{\longrightarrow} \Tot \check{C}(\mathscr{U},F).
\end{equation}
Therefore, it represents the right
derived functor $Rj_*$ of $j_*$, that is,
$$Rj_* = j_*\Tot\check{C}(\mathscr{U}):D\Qcoh(\A_U) \to D\Qcoh(\A_X).$$

\subsection{Lemma}
\label{lem:DZX=DZU}
{\it
Let $X$ be a scheme with an ample family of line bundles, $j:U \subset X$ a
quasi-compact open subscheme,  $Z \subset X$ a closed subset with
quasi-compact open complement $X-Z$ such that $Z \subset U$, then
we have an equivalence of triangulated categories
$$j^*: D_Z\Qcoh(\A_X) \stackrel{\simeq}{\longrightarrow} D_Z\Qcoh(\A_U)$$ 
with inverse the functor $Rj_*$.
}

\subsection*{\it Proof}
We first check that $Rj_*$ preserves cohomological support.
Denote by $j_U:U-Z \subset U$, $j_X:X-Z \subset X$ and $j_Z:U-Z \subset
X-Z$ the corresponding open immersions, and note that the canonical map 
$j_X^*j_*M \to j_{Z*}j_U^*M$ is an isomorphism for every quasi-coherent
$\A_U$-module $M$.
By the existence of an ample family of line bundles on $X$, we can choose a
finite open cover $\mathscr{U}=\{U_0,...,U_n\}$ of $U$ such that 
all inclusions $U_{i_0}\cap ... \cap U_{i_k} \subset X$ are affine maps.
For a complex $F$ of quasi-coherent $\A_U$-modules, we have
$$j_X^*Rj_*F=j_X^*j_*\Tot\check{C}(\mathscr{U},F) =
j_{Z*}j_U^*\Tot\check{C}(\mathscr{U},F) =
j_{Z*}\Tot\check{C}(\mathscr{U}-Z,j_U^*F)$$ 
where $\mathscr{U}-Z$ is the cover $\{U_0-Z,...,U_n-Z\}$ of $U-Z$.
As pull-backs of affine maps, all inclusions $(U_{i_0}-Z)\cap ... \cap
(U_{i_k}-Z)  \subset X-Z$ are also affine maps, so that the functor
$j_{Z*}\Tot\check{C}(\mathscr{U}-Z)$ represents $Rj_{Z*}$, and we obtain a
natural isomorphism of functors
$$j_X^*\circ Rj_* \stackrel{\cong}{\longrightarrow} Rj_{Z*}\circ j_U^*.$$
In particular, $Rj_*$ sends $D_Z\Qcoh(\A_U)$ into $D_Z\Qcoh(\A_X)$.

For the proof of the lemma, note that the composition
$j^*Rj_* = j^*j_*\Tot\check{C}(\mathscr{U}) = \Tot\check{C}(\mathscr{U})$ is
naturally quasi-isomorphic to the identity functor via the map
(\ref{eqn:Rj*NatQuis}).
Furthermore, the unit of adjunction $F \to Rj_*\circ j^*(F)$, which is 
adjoint to the map (\ref{eqn:Rj*NatQuis}) applied to $j^*F$, is a
quasi-isomorphism for $F \in D_Z\Qcoh(\A_X)$ since both complexes have
cohomological support 
in $Z$, so that we can check this property by restricting the map to $U$,
where it is the quasi-isomorphism (\ref{eqn:Rj*NatQuis}).
\Qed

\subsection{Proposition}
\label{prop:DXDUcofinal2}
{\it
Let $X$ be a scheme with an ample family of line-bundles.
Let $Z \subset X$ be a closed subscheme with quasi-compact open complement
$X-Z$.
Let $j:U \subset X$ be a quasi-compact open subscheme.
Let $\A_X$ be a quasi-coherent sheaf of $O_X$-algebras.
Then the following hold.
\begin{enumerate}
\item
\label{prop:DXDUcofinal2A}
Restriction of vector bundles induces a fully faithful triangle functor 
$$j^* : \D\sPerf(\A_X\phantom{i}on\phantom{i}Z,U\Quis) \hookrightarrow 
\D\sPerf(\A_U\phantom{i}on\phantom{i}Z\cap U,U\Quis).$$
\item
\label{prop:DXDUcofinal2B}
The triangulated category $\D_Z\Qcoh(\A_X)$ is compactly generated and the
triangle functor 
$\D\sPerf(\A_X\phantom{i}on\phantom{i}Z,\quis) \to \D_Z\Qcoh(\A_X)$ induces an
equivalence of $\D\sPerf(\A_X\phantom{i}on\phantom{i}Z,\quis)$ with 
the full triangulated subcategory of compact objects in $\D_Z\Qcoh(\A_X)$.
\item
\label{prop:DXDUcofinal2C}
If $Z \subset U$, then restriction of vector bundles induces an equivalence of
triangulated categories 
$$j^* : \D\sPerf(\A_X\phantom{i}on\phantom{i}Z,\quis) \stackrel{\simeq}{\to} 
\D\sPerf(\A_U\phantom{i}on\phantom{i}Z,\quis).$$
\item
\label{prop:DXDUcofinal2D}
The triangle functor in \ref{prop:DXDUcofinal2A} is cofinal.
\end{enumerate}
}

\subsection*{\it Proof}
The functor in \ref{prop:DXDUcofinal2A} is clearly conservative.
It is full by the following argument.
Let $A$ and $B$ be strictly perfect complexes of $\A_X$-modules with support
in $Z$, and let $j^*A \stackrel{_{\sim}}{\leftarrow} E \to j^*B$ be a diagram in
$\sPerf(\A_U\phantom{i}on\phantom{i}U\cap Z)$ representing a map $f:j^*A \to
j^*B$ in $\D^c(\A_U\phantom{i}on\phantom{i}Z\cap U,U\Quis)$ where $E
\stackrel{_{\sim}}{\to} j^*A$ is a $U$-quasi-isomorphism.
Let $M$ be the pull-back of $j_*E \to j_*j^*A$ and the $U$-isomorphism $A \to
j_*j^*A$. 
The induced maps $M \to j_*E$ and $M \to A$ are $U$-isomorphism and
$U$-quasi-isomorphism, respectively.
By lemma \ref{lem:UquisAppr} with $\A = \sPerf(\A_X\phantom{i}on\phantom{i}
Z)$, there is a $A_0\in \A$ and a $U$-quasi-isomorphism $A_0 \to M$.
By lemma \ref{lem:CalcFracByKoszul}, there is a map $A_0\otimes K(s^n) \to B$
such that the two maps $A_0\otimes K(s^n) \to A_0 \to M \to j_*E \to j_*j^*B$
and $A_0\otimes K(s^n) \to B \to j_*j^*B$ coincide.
If follows that the map $f:j^*A \to j^*B$ is the image of the map in 
$\D^c(\A_X\phantom{i}on\phantom{i}Z,U\Quis)$ which is represented by the diagram
$A \stackrel{\sim}{\leftarrow} A_0\otimes K(s^n) \to B$.
Therefore, the functor in \ref{prop:DXDUcofinal2A} is full.
Any conservative and full triangle functor is faithful, hence the triangle
functor in \ref{prop:DXDUcofinal2A} is fully faithful.

It follows from proposition \ref{prop:DXgenByLineBdls} 
and Neeman's theorem \ref{thm:neeman} \ref{thm:neemanA}
that 
the functor in \ref{prop:DXDUcofinal2A} is cofinal for $Z=X$ since 
in this case both categories 
contain as a cofinal subcategory the triangulated category generated by
$\A_X\otimes L$ where $L$ runs through the line bundles on $X$.
This shows part \ref{prop:DXDUcofinal2D} when $Z=X$.

In order to prove \ref{prop:DXDUcofinal2B}, write $\scR$ for the compactly
generated 
triangulated category $D\Qcoh(\A_X)$ with category of compact objects 
$\scR^c=D\sPerf(\A_X,\quis)$, see \ref{prop:DXgenByLineBdls}.
Let $\scS\subset \scR$ be the full triangulated subcategory closed under the
formation of coproducts in $\scR$ which is generated by the set
$S=\sPerf(\A_X\phantom{i}on\phantom{i}Z,\quis) \subset \scR^c$ of compact
objects.
By part \ref{prop:DXDUcofinal2D} for $Z=X$ proved above and proposition
\ref{prop:DXgenByLineBdls}, we have a cofinal inclusion  
$\scR^c/\scS^c \to D^c\Qcoh(\A_U)$.
By Neeman's theorem \ref{thm:neeman} \ref{thm:neemanB} and \ref{thm:neemanC},
this implies that the functor $\scR/\scS \to D\Qcoh(\A_U)$ is an equivalence.
In particular $\scS$ is the kernel category of the functor 
$D\Qcoh(\A_X) \to D\Qcoh(\A_U)$, so that $\scS = D_Z\Qcoh(\A_U)$ is compactly
generated by $D\sPerf(\A_X\phantom{i}on\phantom{i}Z,\quis)$.
Since $D\sPerf(\A_X,\quis) = D^c\Qcoh(\A_X)$ is idempotent complete, its epaisse
subcategory $D\sPerf(\A_X\phantom{i}on\phantom{i}Z,\quis)$ is also idempotent
complete, so that we have the identification
$D\sPerf(\A_X\phantom{i}on\phantom{i}Z,\quis) = 
D^c_Z\Qcoh(\A_X)$, by Neeman's theorem \ref{thm:neeman} \ref{thm:neemanA}.

In view of \ref{prop:DXDUcofinal2B}, the map in \ref{prop:DXDUcofinal2C} is
the restriction to compact objects of the equivalence of lemma
\ref{lem:DZX=DZU}.
It is therefore also an equivalence.

For the proof of \ref{prop:DXDUcofinal2D}, we simplify notation by writing 
$D^c(\A_X\phantom{i}on\phantom{i}Z,w)$ for the triangulated category
$D\sPerf(\A_X\phantom{i}on\phantom{i}Z,w)$.
Let $V=U\cup (X-Z)$ and consider the commutative diagram
of triangulated categories
$$\xymatrix{
D^c(\A_X\phantom{i}on\phantom{i}Z,V\Quis) \ar[r] \ar[d] & 
D^c(\A_X\phantom{i},V\Quis) \ar[r] \ar[d] &
D^c(\A_X\phantom{i},X-Z\Quis)\ar[d]\\
D^c(\A_V\phantom{i}on\phantom{i}Z\cap U,V\Quis) \ar[r] & 
D^c(\A_V\phantom{i},V\Quis) \ar[r] &
D^c(\A_V\phantom{i},X-Z\Quis)
}$$
in which all vertical functors are fully faithful, by \ref{prop:DXDUcofinal2A},
and the middle and the right
vertical functors are cofinal since all four categories have as cofinal
subcategory the triangulated category generated by $\A_X\otimes L$, where $L$
runs through the line-bundles on $X$, by proposition
\ref{prop:DXgenByLineBdls}.
Therefore, the right two vertical functors are equivalences after idempotent
completion. 
Since the two left triangulated categories are the ``kernel categories'' of
the two right horizontal functors, and this property is preserved under
idempotent completion, the left vertical functor is also an equivalence after
idempotent completion.
Thus, the left vertical functor is cofinal.

The functor 
$\D^c(\A_X\phantom{i}on\phantom{i}Z,U\Quis) \hookrightarrow 
\D^c(\A_U\phantom{i}on\phantom{i}Z\cap U,U\Quis)$
in \ref{prop:DXDUcofinal2A} can be identified with the left vertical functor
in the diagram since  
$U$-quasi-isomorphisms are $V$-quasi-isomorphisms for complexes of
$A_X$-modules cohomologically supported in
$Z$, and since the functor 
$\D^c(\A_V\phantom{i}on\phantom{i}Z\cap U,V\Quis) \to
\D^c(\A_U\phantom{i}on\phantom{i}Z\cap U,U\Quis)$
is an equivalence, by \ref{prop:DXDUcofinal2B}.
\Qed

\subsection{Corollary}
\label{cor:K0Lift}
{\it
Let $X$ be a scheme which has an ample family of line-bundles, let $Z
\subset X$ be a closed subset with quasi-compact open complement $X-Z$, and
let $j:U \subset X$ be a quasi-compact open subscheme. 
Let $M$ be a quasi-coherent $\A_X$-module 
such that $j^*M$ is strictly perfect on $U$ and has cohomological support in
$Z \cap U$.
If the class $[j^*M] \in
K_0(\A_U\phantom{,}on\phantom{,}Z\cap U)$ is in the image of the map
$K_0(\A_X\phantom{,}on\phantom{,}Z) \to K_0(\A_U\phantom{,}on\phantom{,}Z\cap
U)$, then there is a $U$-quasi-isomorphism 
$$\xymatrix{A\hspace{2ex} \ar[r]^{\exists}_{U-quis}& \hspace{2ex}M}$$
with $A$  a strictly perfect complex of $\A_X$-modules which has cohomological
support in $Z$.
}

\subsection*{\it Proof}
We start with a standard fact about $K_0$ of triangulated categories.
Let $\T_0 \subset \T_1$ be a (fully faithful and) cofinal functor between
triangulated categories.
Then an object $T$ of $\T_1$ is isomorphic to an object of $\T_0$ if and only
if its class $[T] \in K_0(\T_1)$ is in the image of $K_0(\T_0) \to K_0(\T_1)$.
This is because the cokernel of $K_0(\T_0) \to K_0(\T_1)$ can be identified
with the quotient monoid of the abelian monoid of isomorphism classes of
objects in $\T_1$ under direct sum modulo the submonoid of isomorphism classes
of objects in $\T_0$, so that
an object of $\T_1$ defines the zero class in the cokernel if
and only if it is stably in $\T_0$.
But for triangulated categories, an object is stably in $\T_0$ iff it is
isomorphic to an object in $\T_0$.

For the proof of the corollary, we apply this argument to the 
inclusion in proposition \ref{prop:DXDUcofinal2} \ref{prop:DXDUcofinal2A}
which is cofinal, by \ref{prop:DXDUcofinal2} \ref{prop:DXDUcofinal2D}. 
We see that $j^*M$ is isomorphic 
in $\D^c(\A_U\phantom{i}on\phantom{i}Z\cap U,U\Quis)$ to
an object $j^*B$, where $B$ is a perfect complex of $\A_X$-modules with
cohomological support in $Z$. 
It follows that there is a zig-zag of $U$-quasi-isomorphisms $j^*M
\stackrel{\sim}{\leftarrow} F \stackrel{\sim}{\to} j^*B$.
Let $P$ be the pull-back of $j_*F\to j_*j^*B$ along the $U$-isomorphism $B \to
j_*j^*B$.
Then $P \to j_*F$ is a $U$-isomorphism, and it follows that $P \to B$ is a
$U$-quasi-isomorphism.
By lemma \ref{lem:UquisAppr} with $\A$ the subcategory of those strictly
perfect complexes which are cohomologically supported in $Z$, there is a
$U$-quasi-isomorphism 
$B' \to P$ with $B'$ strictly perfect and cohomologically supported in $Z$.
Since $X$ has an ample family of line-bundles and $U$ is quasi-compact, we can
choose line-bundles $L_i$ and global
sections $s_i \in \Gamma(X,L_i)$, $i=1,...,l$ such that
the set of non-vanishing loci $X_{s_i}$, $i=1,...,l$ is an affine open cover
of $U$.
By Lemma \ref{lem:CalcFracByKoszul}, we can find an $l$-tuple of negative
integers $n$ such that the composition of $U$-quasi-isomorphisms
$A=B' \otimes K(s) \to B' \to j_*j^*M$ lifts to $M$.
\Qed

\subsection{Proposition}
\label{prop:MainLcn}
{\it
Let $X$ be a scheme which has an ample family of line bundles, let $Z
\subset X$ be a closed subscheme with quasi-compact open complement, and let
$j:U \subset X$ be a quasi-compact open subscheme.
Let $\A_X$ be a quasi-coherent $O_X$-algebra with involution.
Then for any line-bundle $L$ on $X$ and any integer $n \in \Z$, restriction of
$\A_X$ vector bundles to $U$ defines non-singular exact form functors
$$
\renewcommand\arraystretch{2} 
\begin{array}{l}
\hspace{-1ex}
(\, \sPerf(\A_X\phantom{i}on\phantom{i}Z),\, U\Quis,\,
\sharp^n_L\, )\\
\hspace{25ex} \longrightarrow \hspace{1ex}
(\, \sPerf(\A_U\phantom{i}on\phantom{i}U\cap Z),\, U\Quis,\,
\sharp^n_{j^*L}\, )
\end{array}
$$
which induce isomorphisms on higher Grothendieck-Witt groups $GW_i$ for $i
\geq 1$ and a monomorphism for $GW_0$. 

If, moreover, we have $Z \subset U$, then the form functors induce 
isomorphisms for all higher Grothendieck-Witt groups $GW_i$, $i \geq 0$.
}

\subsection*{\it Proof}
Let $\sPerf_{K_0}(\A_U\phantom{,}on\phantom{,}\U\cap Z) \subset 
\sPerf(\A_U\phantom{,}on\phantom{,}\U\cap Z)$
be the full subcategory of those strictly perfect complexes of $\A_U$-modules
with cohomological support in $U\cap Z$ which have class in the image of
$K_0(\A_X\phantom{i}on\phantom{i}Z) \to K_0(\A_U\phantom{i}on\phantom{i}Z\cap
U)$.  
By cofinality \ref{thm:cofinal}, the duality preserving inclusion 
\begin{equation}
\label{eqnK0Eq}
\sPerf_{K_0}(\A_U\phantom{,}on\phantom{,}\U\cap Z) \to 
\sPerf(\A_U\phantom{,}on\phantom{,}\U\cap Z)
\end{equation}
of exact categories with weak equivalences the $U$-quasi-isomorphisms and
duality $\sharp_{j^*L}^n$ 
induces maps on higher Grothendieck-Witt groups $GW_i$ which are isomorphisms
for $i \geq 1$ and a monomorphism for $i=0$.
Restriction of vector-bundles defines a non-singular exact form functor
\begin{equation}
\label{eqnK0Eq2}
(\, \sPerf(\A_X\phantom{,}on\phantom{,}Z),\, U\Quis\, ) \to 
(\, \sPerf_{K_0}(\A_U\phantom{,}on\phantom{,}\U\cap Z),\, U\Quis\, )
\end{equation}
which induces a homotopy equivalence of Grothendieck-Witt spaces by theorem
\ref{thm:ApprFunctorialCalcFrac}, 
where \ref{thm:ApprFunctorialCalcFrac} (c) follows from
corollary \ref{cor:K0Lift} and lemma \ref{lem:C'replacesC};  
\ref{thm:ApprFunctorialCalcFrac} (e) and (f)
are proved in lemma \ref{lem:CalcFracByKoszul}; the remaining hypothesis
of theorem \ref{thm:ApprFunctorialCalcFrac} being trivially satisfied.

If $Z \subset U$, then $K_0(\A_X\phantom{i}on\phantom{i}Z) =
K_0(\A_U\phantom{i}on\phantom{i}Z\cap U)$, by proposition
\ref{prop:DXDUcofinal2} \ref{prop:DXDUcofinal2C},
so that (\ref{eqnK0Eq}) is the
identity inclusion, and the map (\ref{eqnK0Eq2}), which induces a homotopy
equivalence of Grothendieck-Witt spaces, is the map in the proposition.
\Qed

\subsection*{\it Proof of theorem \ref{thm:Locn1}}
By the ``change-of-weak-equivalence 
theorem'' \ref{thm:ChgOfWkEq}, the sequence
$$
(\, \sPerf(\A_X\phantom{i}on\phantom{i}Z),\, \quis,\,
\sharp^n_L\, ) \to 
(\, \sPerf(\A_X),\, \quis,\,
\sharp^n_L\, ) \to 
(\, \sPerf(\A_X),\, U\Quis,\,
\sharp^n_L\, )
$$
induces a homotopy fibration of Grothendieck-Witt spaces.
By proposition \ref{prop:MainLcn}, the form functor
$$
(\, \sPerf(\A_X),\, U\Quis,\, \sharp^n_L\, ) \to
(\, \sPerf(\A_U),\, \quis,\, \sharp^n_{j^*L}\, )
$$
induces isomorphisms on $GW_i$ for $i\geq 1$ and a monomorphism for $i=0$.
\Qed

\subsection*{\it Proof of theorem \ref{thm:ZarExc1}}
The theorem is a special case of proposition \ref{prop:MainLcn}.
\Qed

\section{Extension to negative Grothendieck-Witt groups}
\label{sec:NegGW}

For an open subscheme $U \subset X$, the restriction map $GW_0(X) \to GW_0(U)$
is not surjective, in general, not even if $X$ is regular.
The purpose of this section is to extend the long exact sequence associated
with the homotopy fibration of theorem 
\ref{thm:Locn1} to negative degrees.
Theorems \ref{thm:Locn1} and \ref{thm:ZarExc1} will be extended to a fibration
and a weak equivalence of non-connective spectra.

\subsection{Cone and suspension of $\A_X$}
\label{subsec:ConeAndSusp}
The {\em cone ring} is the ring $C$ of infinite matrices
$(a_{i,j})_{i,j\in\N}$ with coefficients $a_{i,j}$ in $\Z$ for which 
each row and each column has only finitely many non-zero entries.
Transposition of matrices ${^t(a_{i,j})} = (a_{j,i})$ makes $C$ into a ring
with involution.
As a $\Z$-module $C$ is torsion free, hence flat.

The {\em suspension ring} $S$ is the factor ring of $C$ by the two sided ideal
$M_{\infty} \subset C$ of those matrices which have only finitely many
non-zero entries. 
Transposition also makes $S$ into a ring with involution such that the
quotient map $C \twoheadrightarrow S$ is a map of rings with involution. 
For another description of the suspension ring $S$, consider the matrices $e_n
\in C$, $n\in \N$, with entries $(e_n)_{i,j}=1$ for $i=j\geq n$ 
and zero otherwise.
They are symmetric idempotents, \ie ${^te_n}=e_n=e_n^2$, and
they form a multiplicative subset of $C$ which satisfies the {\O}re condition,
that is, the multiplicative subset satisfies the axioms for a calculus of
fractions. 
One checks that the quotient map $C \twoheadrightarrow S$ identifies
the suspension ring $S$ with the 
localization of the cone ring $C$ with respect to the elements $e_n\in C$,
$n\in \N$. 
In particular, the suspension ring $S$ is also a flat $\Z$-module.

Let $X$ be a quasi-compact and quasi-separated scheme.
For a quasi-coherent sheaf $\A_X$ of $O_X$-algebras, write $C\A_X$ and $S\A_X$
for the quasi-coherent sheaves of $O_X$-algebras associated with the
presheaves $C\otimes_{\Z}\A_X$ and $S\otimes_Z\A_X$.
On quasi-compact open subsets $U \subset X$, we have
$(C\A_X)(U)=C\otimes_{\Z}A_X(U)$ and $S\A_X=S\otimes_{\Z}A_X(U)$, by flatness
of $C$ and $S$.
If $\A_X$ is a sheaf of algebras with involutions, then the involutions on $C$
and on $S$ make $C\A_X$ and $S\A_X$ into sheaves of $O_X$-algebras with
involution.

Let $\meps =1-e_1 \in C$ be the symmetric idempotent with entries $1$ at
$(0,0)$ and zero otherwise. 
The image $C\meps$ of the right multiplication map $\times \meps: C \to C$
is a finitely generated projective left $C$-module.
It is equipped with a symmetric form $\ffi:C\meps
\otimes_{\Z}(C\meps)^{op} \to C :x\otimes y^{op} \mapsto x\cdot {^ty}$.
The idempotent $\meps$ makes
$(C\meps,\ffi)$ into a direct factor of the unit symmetric form $(C,\mu)$, see
\ref{subsec:basicMffiOtimes} \ref{subsec:basicMffiOtimesA}.
Therefore, tensor product $(C\meps,\ffi)\otimes_{\Z}?$ defines a non-singular
exact form functor
$$\iota: (\sPerf(\A_X),\sharp_L^n) \to (\sPerf(C\A_X),\sharp_L^n): V \mapsto C\meps \otimes_{\Z}V.$$
Since $S$ is a flat $C$-algebra, the quotient map $C \to S$ induces an exact
functor $\rho:C\A_X\Lmod \to S\A_X\Lmod:M \mapsto S\otimes_CM$ on categories
of modules which sends vector bundles to vector bundles.
The two functors $\iota$ and $\rho$ yield a sequence of non-singular exact
form functors
\begin{equation}
\label{eqn:iotaRho}
(\sPerf(\A_X),\sharp_L^n) \stackrel{\iota}{\longrightarrow}
(\sPerf(C\A_X),\sharp_L^n) 
\stackrel{\rho}{\longrightarrow} (\sPerf(S\A_X),\sharp_L^n).
\end{equation}
The functors satisfy $\rho\circ\iota=0$ because $S\otimes_CC\meps = \im(\times
\meps:S \to S)=0$ as $0=\meps\in S$.

The following theorem will allow us to extend the results of
\S\ref{sec:LocnZar} to negative Grothendieck-Witt groups.

\subsection{Theorem}
\label{thm:deloop}
{\it
Let $X$ be a scheme with an ample family of line-bundles, let $Z \subset X$ be
a closed subset with quasi-compact open complement $X-Z$, and let $\A_X$ be a
quasi-coherent $O_X$-algebra with involution.
Then for any line bundle $L$ on $X$, and any integer $n \in \Z$, the sequence 
(\ref{eqn:iotaRho})
induces a homotopy fibration of Grothendieck-Witt spaces with contractible
total space 
$$GW^n(\A_X\phantom{i}on\phantom{i}Z,L) \to
GW^n(C\A_X\phantom{i}on\phantom{i}Z, L) \to 
GW^n(S\A_X\phantom{i}on\phantom{i}Z, L).$$
}
\vspace{2ex}

The proof of theorem \ref{thm:deloop} will occupy sections
\ref{lem:iotalFullFaith} to \ref{dfn:GWspct}.

\subsection{Lemma}
\label{lem:iotalFullFaith}
{\it
The functor $\iota$ in (\ref{eqn:iotaRho}) is fully faithful.
}

\subsection*{\it Proof}
The image $\meps C$ of the left multiplication map $\meps \times: C \to
C$ is a right $C$-module. 
We have a $\Z$-bimodule isomorphism $\eta: \Z \to \meps C \otimes_C C\meps: 1
\mapsto \meps \otimes_C \meps = 1 \otimes \meps = \meps \otimes 1$ and a
$C$-bimodule map $\mu: C\meps \otimes_{\Z} \meps C \to C: A\meps \otimes \meps
B \mapsto A\meps B$ such that the compositions
$$C\meps \cong C\meps \otimes_{\Z}\Z \stackrel{id\otimes \eta}{\longrightarrow}
C\meps \otimes_{\Z}\meps C \otimes_C C\meps \stackrel{\mu\otimes
  id}{\longrightarrow} C \otimes_C C\meps \cong C\meps\hspace{3ex}{\rm and}
$$
$$\meps C \cong \Z \otimes_{\Z} \meps C \stackrel{\eta \otimes
  id}{\longrightarrow} \meps C \otimes_C C \meps \otimes_{\Z} \meps C 
\stackrel{id \otimes \mu}{\longrightarrow} \meps C \otimes_C C \cong \meps C\phantom{\hspace{3ex}{\rm and}}$$
are the identity maps.
It follows that $\eta$ and $\mu$ define unit and counit of an adjunction
between the functors 
$\A_X\Lmod \to C\A_X\Lmod:M\mapsto C\meps\otimes_{\Z}M$ and 
$C\A_X\Lmod \to \A_X\Lmod:N\mapsto \meps C\otimes_CN$.
Since the unit $\eta$ is an isomorphism, the first functor is fully faithful.
In particular, $\iota$ is fully faithful.
\Qed

\subsection{Proposition}
\label{prop:SuspExSeq}
{\it
The sequence of triangulated categories
$$D^b\Vect(\A_X) \stackrel{\iota}{\longrightarrow} D^b\Vect(C\A_X)
\stackrel{\rho}{\longrightarrow} D^b\Vect(S\A_X)$$
is exact up to direct factors.
}

\subsection*{\it Proof}
The multiplication map $\mu:C\meps \otimes_{\Z}\meps C \to C$ factors through
$M_{\infty} \subset C$ and induces an isomorphism $\mu: C\meps
\otimes_{\Z}\meps C \to 
M_{\infty}$ (it is a filtered colimit of isomorphisms of finitely generated
free $\Z$-modules). 
The exact functors 
$$\Qcoh(\A_X) \stackrel{\iota}{\longrightarrow} \Qcoh(C\A_X)
\stackrel{\rho}{\longrightarrow} \Qcoh(S\A_X)$$
have exact right adjoints $\kappa: \Qcoh(C\A_X) \to \Qcoh(\A_X): M \mapsto
\meps C \otimes_C M$ and $\sigma: \Qcoh(S\A_X) \to \Qcoh(C\A_X): N \mapsto N$
such that for a left $C\A_X$-module $M$ the adjuntion maps $\iota\kappa \to
id$ and $id \to \sigma\rho$ are part of a functorial exact sequence 
\begin{equation}
\label{eqn:prop:SuspExSeq}
0 \to \iota\kappa M \to M \to \sigma\rho M \to 0
\end{equation}
which is the tensor product (over $C$) of $M$ with the exact sequence of flat
$C$-modules $0 \to M_{\infty} \to C \to S \to 0$.
It follows that the sequence of triangulated categories
$$D\Qcoh(\A_X) \stackrel{\iota}{\longrightarrow} D\Qcoh(C\A_X)
\stackrel{\rho}{\longrightarrow} D\Qcoh(S\A_X)$$
is exact as $\kappa$ and $\rho$ induce right adjoint functors on derived
categories, and (\ref{eqn:prop:SuspExSeq}) induces a functorial distinguished
triangle for every object of $D\Qcoh(C\A_X)$. 
By proposition \ref{prop:DXgenByLineBdls}, 
these triangulated categories are compactly generated.
Since $\iota$
and $\rho$ preserve compact objects, the associated sequence of compact
objects -- which is the sequence in proposition \ref{prop:SuspExSeq} --
is exact up to factors \ref{thm:neeman}.
\Qed
\vspace{2ex}

Let $\sPerf_S(C\A_X) \subset \sPerf(C\A_X)$ be the full subcategory of those
complexes $V$ for which $S\otimes_CV$ is acyclic.
This subcategory is closed under the involution $\sharp_L^n$, so that 
$\sPerf_S(C\A_X)$ inherits the structure of an exact category with weak
equivalences and duality from $\sPerf(C\A_X)$.
Since $\rho\iota=0$, $\iota$ induces a non-singular exact form functor
$\iota:\sPerf(\A_X) \to \sPerf_S(C\A_X)$.

\subsection{Proposition}
\label{prop:sPerfA=sPerfCA}
{\it
For any line bundle $L$ on $X$, and any $n\in \Z$, the functor 
$\iota$ induces a homotopy equivalence 
$$GW^n(\sPerf(\A_X),\quis,L) \stackrel{\sim}{\longrightarrow} GW^n(\sPerf_S(C\A_X),\quis,L).$$
}

\subsection*{\it Proof}
The proof is a consequence of theorem \ref{thm:ApprFunctorialCalcFrac} (or of
lemma \ref{lem:resolution}).
Since $\iota$ is fully faithful, conditions (e) and (f) are satisfied.
Since $\iota$ induces a fully faithful functor on derived categories, by
proposition \ref{prop:SuspExSeq}, condition (b) is also satisfied.
The only non-trivial condition to check is (c).
By lemma \ref{lem:C'replacesC}, we only need to
show that for every $M \in \sPerf_S(C\A_X)$ there is an $A \in \sPerf(\A_X)$
and a quasi-isomorphism $C\meps \otimes_{\Z} A \to M$.
Let $M$ be a strictly perfect complex of $C\A_X$ modules with $S\otimes_CM$
acyclic. 
By proposition \ref{prop:SuspExSeq}, there is a zigzag of quasi-isomorphisms
$C\meps \otimes_{\Z} B \leftarrow N \to M$ in $\sPerf_S(C\A_X)$ with $B \in
\sPerf(\A_X)$. 
Since $N \in \sPerf_S(C\A_X)$, proposition \ref{prop:SuspExSeq} implies that
the counit of adjunction  
$C\meps \otimes_{\Z}\meps C \otimes_C N \to N$ is a quasi-isomorphism.
We apply lemma \ref{lem:UquisAppr} with $C\A_X$ in place of $\A_X$ and $U=X$,
$\A = \sPerf(\A_X)$ to the quasi-isomorphism 
$\meps C \otimes_C N \to \meps C \otimes_C C\meps \otimes_{\Z} B \cong B$,
and obtain a strictly perfect complex $A$ of $\A_X$-modules and a
quasi-isomorphism $A \to \meps C \otimes_C N$.
Finally, the composition
$C\meps \otimes_{\Z} A \to C\meps \otimes_{\Z} \meps C \otimes_C N =
M_{\infty} \otimes_C N \to N \to M$ is a quasi-isomorphism.
\Qed
\vspace{2ex}

For a quasi-coherent $O_X$-algebra $\A_X$, call an $\A_X$-module $M$ {\em
  quasi-free} if it is isomorphic to a 
finite direct sum
$\bigoplus_i \A_X\otimes L_i$ of $A_X$-modules of 
the form $\A_X\otimes L_i$ for some line bundles $L_i$ on $X$.
Note that a quasi-free $\A_X$-module is a vector bundle.

\subsection{Lemma}
\label{lem:CAMcalcFrac}
{\it
Let $X$ be a quasi-compact and quasi-separated scheme, and let $\A_X$ be a
quasi-coherent sheaf of $O_X$-algebras.
Let $A,M$ be quasi-coherent $C\A_X$-modules with $A$ quasi-free.
Then the following hold.
\begin{enumerate}
\item
\label{lem:CAMcalcFracA}
For every map $f: \rho A \to \rho M$ of $S\A_X$-modules, 
there are maps $s: B \to A$ and $g:B \to M$ of $C\A_X$-modules
with $B$ quasi-free such that $f\circ 
\rho(s)=\rho(g)$ and $\rho(s)$ an isomorphism.
\item
\label{lem:CAMcalcFracB}
For any two maps $f,g:A \to M$ of $C\A_X$-modules such that
$\rho(f)= \rho(g)$ there is a map $s:B \to A$ of of
quasi-free $C\A_X$-modules such that $f\circ s = g\circ s$ and $\rho(s)$ is an 
isomorphism.
\end{enumerate}
}

\subsection*{\it Proof}
The proof reduces to 
$A=C\A_X\otimes L$ with $L$ a line-bundle on $X$.
For such an $A$, the map $Hom_{C\A_X}(A,M) \to Hom_{S\A_X}(\rho S,\rho
M)$ can be identified with the map on global sections 
$\Gamma(X,M\otimes L^{-1}) \to \Gamma(X,S\otimes_C M \otimes
L^{-1})=S\otimes_C \Gamma(X,M\otimes L^{-1})$.
This map is surjective since $C \to S$ is, proving \ref{lem:CAMcalcFracA}.
The map is also a localization by a calculus of fractions with respect to the
set of elements $e_n \in C$, $n\in \N$, of \ref{subsec:ConeAndSusp}.
This shows that \ref{lem:CAMcalcFracB} also holds.
\Qed

\subsection{Lemma}
\label{lem:EilenbergSwindle}
{\it
Let $X$ be a quasi-compact and quasi-separated scheme, $L$ a line bundle on
$X$, and let $\A_X$ be a
quasi-coherent sheaf of $O_X$-algebras with involution.
Then for every $n\in \Z$, the Grothendieck-Witt space 
$$GW^n(C\A_X,L)$$
is contractible.
}

\subsection*{\it Proof \rm (compare \cite{karoubi:heidelberg})}
We will define a $C$-bimodule $M$, which is finitely generated
projective as left 
$C$-module, together with a symmetric form $\ffi:
M\otimes_CM^{op} \to C$ in $C\Bimod$ whose adjoint $M \to [M^{op},C]_C$ is an
isomorphism.
Furthermore, we will construct an isometry
$(C,\mu) \perp (M,\ffi) \cong (M,\ffi)$ of symmetric forms in $C\Bimod$.
Therefore, tensor product $(M,\ffi)\otimes_C?$ defines a non-singular exact
form functor 
$(F,\ffi):(\sPerf(C\A_X),\quis,\sharp_L^n) \to
(\sPerf(C\A_X),\quis,\sharp_L^n)$ which satisfies $id \perp (F,\ffi) \cong
(F,\ffi)$, so that on 
higher Grothendieck-Witt groups we have
$GW_i^n(id)+GW^n_i(F,\ffi)=GW^n_i(F,\ffi)$ which implies $GW^n_i(id)=0$, that
is, $GW^n_i(C\A_X,L)=0$, hence $GW^n(C\A_X,L)$ is contractible.

To construct $(M,\ffi)$ and the bimodule isometry  $(C,\mu) \perp (M,\ffi) \cong (M,\ffi)$
we choose a bijection $\sigma:\N \stackrel{\cong}{\to} \N\times \N:n\mapsto
(\sigma_1(n),\sigma_2(n))$
 and define a homomorphism of rings with involutions 
$$\sigma: C \to C: a \mapsto \sigma(a) \phantom{qw}{\rm with}\phantom{qw}
\sigma(a)_{ij} = \left\{\begin{array}{cl}
                         a_{\sigma_1(i),\sigma_1(j)}& {\rm
                           if}\phantom{ab}\sigma_2(i)=\sigma_2(j)\\
                         0 & {\rm otherwise.}
                       \end{array}
                 \right.
$$
The $C$-bimodule $M$ is $C$ as a left module, and has right multiplication
defined by $M \times C \to M: (x,a) \mapsto x\cdot I(a)$.
The symmetric form $\ffi$ is the $C$-bimodule map $M\otimes_CM^{op} \to
C:x\otimes y^{op} \mapsto x\cdot {^ty}$.
Since, as a left $C$-module, $(M,\ffi)$ is just the unit symmetric form
$(C,\mu)$ (see \ref{subsec:basicMffiOtimes}
\ref{subsec:basicMffiOtimesA}), the adjoint $M \to [M^{op},C]_C$ of $\ffi$ is
an isomorphism.

In order to define the bimodule isometry $(C,\mu) \perp (M,\ffi) \cong
(M,\ffi)$, consider the 
elements $\gamma, \delta \in C$ defined by
$$
\gamma_{ij} = \left\{\begin{array}{cl}
                         1& {\rm
                           if}\phantom{ab}\sigma(j)=(i,0)\\
                         0 & {\rm otherwise,}
                       \end{array}
                 \right.
\phantom{qw}{\rm and}\phantom{qw}
\delta_{ij} = \left\{\begin{array}{cl}
                         1& {\rm
                           if}\phantom{ab}\sigma(j)=\sigma(i)+(0,1)\\
                         0 & {\rm otherwise.}
                       \end{array}
                 \right.
$$
The homomorphism $I$ and the elements $\gamma, \delta \in C$ are related by the
following identities 
$$\begin{array}{lll}
\delta \cdot {^t\gamma}=0, &\hspace{2ex}\gamma \cdot {^t\gamma} = \delta \cdot
{^t\delta} = 1, & \hspace{2ex}{^t\gamma}\cdot\gamma + {^t\delta}\cdot \delta =
1, \\ 
a\cdot \gamma =\gamma\cdot I(a), & \hspace{2ex}I(a)\cdot \delta = \delta \cdot
I(a)& 
\end{array}
$$
for all $a\in C$.
Therefore, the map $C\oplus M \to M:(a,x)\mapsto a\cdot \gamma + x\cdot
\delta$ is a $C$-bimodule isomorphism with inverse the map $M \to C\oplus M: x
\mapsto (x\cdot {^t\gamma}, x \cdot {^t\delta})$.
It preserves forms because $(a\gamma + x\delta)\cdot {^t(b\gamma + y\delta)} =
a\cdot{^tb}+ x\cdot{^ty}$. 
\Qed
\vspace{2ex}

Write $\sPerf^0(\A_X) \subset \sPerf(\A_X)$ for the full subcategory of those
strictly perfect complexes of $\A_X$-modules which are degree-wise quasi-free.
Note that this category is closed under the duality $\sharp_L^n$.
We equip the category $\sPerf^0(\A_X)$ with the degree-wise split exact
structure.
Together with the set of quasi-isomorphisms of complexes of $\A_X$ vector
bundles, it becomes a category of complexes in the sense of Definition
\ref{dfn:CatOfCxs}.

\subsection{Lemma}
\label{lem:sPerf0sPerf}
{\it
Let $X$ be a scheme with an ample family of line bundles.
The inclusion of quasi-free modules into the category of vector bundles 
induces a (fully faithful) cofinal triangle
functor 
$$\D(\, \sPerf^0(\A_X),\, \quis)\, \subset\, 
\D(\, \sPerf(\A_X),\, \quis).$$
Moreover, for every strictly perfect complex $M$ of $\A_X$-modules with class
$[M]$ in the image of the map $K_0(\, \sPerf^0(\A_X),\, \quis) \to 
K_0(\, \sPerf(\A_X),\, \quis)$ there is a quasi-isomorphism $A \to M$ of
complexes of $\A_X$-modules with $A$ a bounded complex of quasi-free modules.
}

\subsection*{\it Proof}
The triangle functor in the lemma is fully faithful, by lemma
\ref{lem:UquisAppr} with $U=X$ and $\A=\sPerf^0(\A_X)$.
It is cofinal, by Neeman's theorem \ref{thm:neeman} \ref{thm:neemanA} and proposition \ref{prop:DXgenByLineBdls}.
Let $\sPerf_{K_0}(\A_X) \subset \sPerf(\A_X)$ be the full subcategory of
those strictly perfect complexes of $\A_X$-modules $M$ whose class $[M]$ is
in the image of the map $K_0(\, \sPerf^0(\A_X),\, \quis) \to 
K_0(\, \sPerf(\A_X),\, \quis)$.
Then the inclusion $(\, \sPerf^0(\A_X),\, \quis) \subset
(\sPerf_{K_0}(\A_X),\, \quis)$ of exact categories with weak equivalences
induces an equivalence of derived categories, so 
that another application of lemma \ref{lem:UquisAppr} with $U=X$ and
$\A=\sPerf^0(\A_X)$ finishes the proof of the claim.
\Qed

\subsection*{\it Proof of theorem \ref{thm:deloop}}
By theorem \ref{thm:Locn1}, we only need to treat the case $Z=X$.
In this case, the total spaces are contractible, by lemma
\ref{lem:EilenbergSwindle}.

Let $\sPerf_{K_0}(S\A_X) \subset \sPerf(S\A_X)$ be the full subcategory of
those strictly perfect complexes of $S\A_X$-modules $E$ whose class $[E]$ is
zero in the Grothendieck group $K_0(\sPerf(S\A_X),\quis)$ of $S\A_X$-vector
bundles. 
Furthermore, call a map $f$ of strictly perfect complexes of
$C\A_X$-modules an $S$-quasi-isomorphism if $\rho(f)$ is a quasi-isomorphism
of complexes of $S\A_X$-modules.
The set of $S$-quasi-isomorphisms is denoted by $S\Quis$.
Consider the commutative diagram of exact categories with weak equivalences
and duality $\sharp_L^n$ induced by inclusions and the map of rings
with involution $C \to S$
$$\xymatrix{
(\sPerf^0(C\A_X),\quis)
 \ar[r] \ar[d] &
(\sPerf^0(C\A_X),S\Quis)
 \ar[r]^{\rho} \ar[d] &
(\sPerf^0(S\A_X),\quis)  \ar[d] \\
(\sPerf(C\A_X),\quis) \ar[r] &
(\sPerf(C\A_X),S\Quis) \ar[r]^{\rho} &
(\sPerf_{K_0}(S\A_X),\quis).
}$$
Note that $K_0$ of all categories with weak equivalences in the diagram is $0$.
For the two left hand categories, this follows from lemma
\ref{lem:EilenbergSwindle} and lemma \ref{lem:sPerf0sPerf}, since for cofinal
triangle functors $\T^0 \subset \T$, the map $K_0(\T^0) \to K_0(\T)$ is
injective.
Since the left horizontal maps are surjective on $K_0$, the middle two
categories with weak equivalences have trivial Grothendieck group $K_0$.
For the upper right corner, vanishing of $K_0$ follows 
moreover from the fact that its $K_0$ is generated by classes of complexes
concentrated in degree $0$
and the fact that
every quasi-free $S\A_X$-module is the
image of a quasi-free $C\A_X$-module, so that the upper horizontal map is
surjective on $K_0$.
Therefore, 
the right vertical and the lower right horizontal functors, which -- a
priori -- have images in $\sPerf(S\A_X)$,  have indeed
image in $\sPerf_{K_0}(S\A_X)$.
We will show that the upper right horizontal and middle and right vertical
functors induce 
equivalences of Grothendieck-Witt spaces (for any duality
$\sharp^n_L$), so that the lower right horizontal functor will induce an
equivalence, too.

The upper right horizontal functor is a localization by a calculus of right
fractions, by lemma \ref{lem:CAMcalcFrac} and
\ref{lem:extendingCalcFrac} \ref{enum:lem:extendingCalcFrac3}.
Therefore, theorem \ref{thm:ApprCalcOfFrac2} shows that it
induces a homotopy equivalence of
Grothendieck-Witt spaces.
For the right vertical functor, the resolution lemma \ref{lem:resolution}
(which we can apply because of lemma \ref{lem:sPerf0sPerf}) shows that it
induces an equivalence of Grothendieck-Witt spaces.
Similarly, by lemma \ref{lem:sPerf0sPerf}, for every strictly perfect complex
of $C\A_X$-modules $M$, there is a bounded complex $A$ of quasi-free
$C\A_X$-modules and a quasi-isomorphism $A \to M$.
A quasi-isomorphism of complexes of $C\A_X$-modules is, a fortiori, an
$S$-quasi-isomorphism, so that the resolution lemma applies to show that the
middle vertical functor induces an equivalence of Grothendieck-Witt spaces.
Summarizing, we have shown that the lower right horizontal functor induces an
equivalence of Grothendieck-Witt spaces.

By the ``change of weak equivalence theorem'' (theorem \ref{thm:ChgOfWkEq}),
the sequence of exact categories with weak equivalences and duality
$\sharp_L^n$ 
$$(\, \sPerf_S(C\A_X),\, \quis\, ) \to (\, \sPerf(C\A_X),\, \quis\, ) \to (\,
\sPerf(C\A_X),\, S\Quis\, )$$
induces a homotopy fibration of Grothendieck-Witt spaces.
Using proposition \ref{prop:sPerfA=sPerfCA} we can replace the left hand term
with $(\, \sPerf(\A_X),\, \quis\, )$.
Using the equivalence of Grothendieck-Witt spaces of the lower right
horizontal functor above and cofinality (theorem
\ref{thm:cofinal}) applied to the inclusion of exact categories with weak
equivalences and duality 
$(\, \sPerf_{K_0}(S\A_X),\, \quis\, ) \subset \
(\, \sPerf(S\A_X),\, \quis\, )$, we can replace the right hand term in the
sequence by $(\, \sPerf(S\A_X),\, \quis\, )$.
\Qed
\vspace{2ex}

Since the total space in the fibration of theorem \ref{thm:deloop} is
contractible, we obtain a homotopy equivalence of spaces
\begin{equation}
\label{eqn:deloop}
GW^n(\A_X\phantom{i}on\phantom{i}Z,L)\, \stackrel{\simeq}{\longrightarrow}\, 
\Omega\ GW^n(S\A_X\phantom{i}on\phantom{i}Z,L).
\end{equation}

\subsection{Definition}
\label{dfn:GWspct}
Let $X$ be a scheme with an ample family of line-bundles, $\A_X$ be a
quasi-coherent sheaf of $O_X$-algebras with 
involution, $L$ a line bundle on $X$, $Z\subset X$ a closed subscheme with
quasi-compact open complement $X-Z$ and $n
\in \Z$ an integer.
The {\em Grothendieck-Witt spectrum}
$$\GW^n(\A_X\phantom{i}on\phantom{i}Z,\phantom{i}L)$$
of symmetric spaces over $\A_X$ with coefficients in the $n$-th shifted
line bundle $L[n]$ and support in $Z$ is
the sequence
$$ GW^n(S^k\A_X\phantom{i}on\phantom{i}Z,\phantom{i}L),\phantom{qw} k\in\N, $$
of Grothendieck-Witt spaces
together with the bonding maps given by the homotopy equivalence
(\ref{eqn:deloop}). 
As usual, if $Z=X$, $n=0$, $\A_X=O_X$ or $L=O_X$, we omit the label
corresponding to $Z$, $n$, $\A$, or $L$, respectively.
 
By construction, we have
$$\pi_i\GW^n(\A_X\phantom{i}on\phantom{i}Z,\phantom{i}L) = 
\left\{\begin{array}{lcl}
GW_i^n(\A_X\phantom{i}on\phantom{i}Z,\phantom{i}L) & {\rm for} & i\geq 0\\
GW_0^n(S^{-i}\A_X\phantom{i}on\phantom{i}Z,\phantom{i}L) & {\rm for} & i\leq 0.
\end{array}
\right.
$$

\subsection{Remark}
By proposition \ref{prop:Periodicity}, there are natural homotopy equivalences
of spectra  
$\GW^n(\A_X\phantom{i}on\phantom{i}Z,\phantom{i}L) \simeq \GW^{n+4}(\A_X\phantom{i}on\phantom{i}Z,\phantom{i}L)$.
\vspace{2ex}

Finally, we are in position to prove the main theorems of this article.

\subsection{Theorem \rm (Localization)}
\label{thm:LocnSp1}
{\it
Let $X$ be a scheme with an ample family of line-bundles, let $Z \subset X$
be a closed subscheme with quasi-compact open complement
 $j:U \subset X$, and let $L$ a line bundle on $X$.
Let $\A_X$ be a quasi-coherent sheaf of $O_X$-algebras with involution.
Then for every $n\in \Z$, the following sequence is a homotopy fibration of
Grothendieck-Witt spectra 
$$\GW^n(\A_X\phantom{i}on\phantom{i}Z,\phantom{i}L) \longrightarrow
\GW^n(\A_X,\phantom{i}L) \longrightarrow \GW^n(\A_U,\phantom{i}j^*L).$$ 
}

\subsection*{\it Proof}
This is because the sequences
$$GW^n(S^i\A_X\phantom{i}on\phantom{i}Z,\phantom{i}L) \longrightarrow
GW^n(S^i\A_X,\phantom{i}L) \longrightarrow GW^n(S^i\A_U,\phantom{i}j^*L)$$ 
are homotopy fibrations for $i \in \N$, by theorem \ref{thm:Locn1}.
\Qed

\subsection{Theorem \rm (Zariski-excision)}
\label{thm:ZarExcSp1}
{\it
Let $X$ be a scheme with an ample family of line-bundles, let $Z \subset X$
be a closed subscheme with quasi-compact open complement,
let $L$ be a line bundle on $X$ and
let $\A_X$ be a quasi-coherent sheaf of $O_X$-algebras with involution.
Then for every $n\in \Z$ and every quasi-compact open subscheme $j:V \subset X$
containing $Z$,
restriction of vector-bundles induces a homotopy equivalence of
Grothendieck-Witt spectra
$$\GW^n(\A_X\phantom{i}on\phantom{i}Z,\phantom{i}L)
\stackrel{\sim}{\longrightarrow}
\GW^n(\A_V\phantom{i}on\phantom{i}Z,\phantom{i}j^*L).$$
}

\subsection*{\it Proof}
This is because the maps
$$GW^n(S^i\A_X\phantom{i}on\phantom{i}Z,\phantom{i}L)
\longrightarrow
GW^n(S^i\A_V\phantom{i}on\phantom{i}Z,\phantom{i}j^*L).$$
are homotopy equivalences for $i \in \N$, by theorem \ref{thm:ZarExc1}.
\Qed

\subsection{Corollary \rm (Mayer-Vietoris for open covers)}
\label{cor:fullMV}
{\it
Let $X=U\cup V$ be a scheme with an ample family of line-bundles
which is covered by two open quasi-compact subschemes $U,V \subset X$.
Let $\A_X$ be a quasi-coherent $O_X$-module with involution.
Let $L$ be a line-bundle on $X$, and $n \in \Z$.
Then restriction of vector bundles induces a homotopy cartesian square of
Grothendieck-Witt spectra
$$\xymatrix{
\GW^n(\A_X,L) \ar[r] \ar[d] & \GW^n(\A_U,L) \ar[d]\\
\GW^n(\A_V,L) \ar[r] & \GW^n(\A_{U\cap V},L).
}$$
}

\subsection*{\it Proof}
The map on vertical homotopy fibres is an equivalence, by theorems
\ref{thm:LocnSp1} and \ref{thm:ZarExcSp1}.
\Qed

\end{document}